\documentclass[10pt,a4paper,reqno]{amsart}

\usepackage[T1]{fontenc}
\usepackage[utf8]{inputenc}
\usepackage{graphicx}
\usepackage{enumerate}
\usepackage{hyperref}
\usepackage{color}
\usepackage{amsmath,latexsym,amssymb,amsthm,array,amsfonts,dsfont}
\usepackage{mathtools}
\usepackage{mathrsfs}
\usepackage{centernot}
\newtheorem{theorem}{Theorem}[section]
\newtheorem{lemma}[theorem]{Lemma}
\newtheorem{definition}[theorem]{Definition}
\newtheorem{corollary}[theorem]{Corollary}
\newtheorem{proposition}[theorem]{Proposition}
\newtheorem{remark}[theorem]{Remark}
\newtheorem{example}[theorem]{Example}
\newtheorem{conjecture}{Conjecture}

\newcommand{\C}{\mathbb C}
\newcommand{\N}{\mathbb N}
\newcommand{\R}{\mathbb R}
\newcommand{\A}{\mathbb C\langle X_1,\ldots,X_n\rangle}

\newcommand{\Dom}{\mathrm{Dom}}
\newcommand{\Ran}{\mathrm{Ran}}
\newcommand{\Ker}{\mathrm{Ker}}

\DeclareMathOperator{\rank}{rank}

\newcommand{\norm}[1]{\left\lVert #1\right\rVert}
\newcommand{\abs}[1]{\left\lvert #1\right\rvert}
\newcommand{\ip}[2]{\left\langle #1,#2\right\rangle}

\begin{document}
\title{Convexity, the spectrum of free Laplacians and geometric consequences: a free KLS bound}
\author{Charles-Philippe Diez$^{\dagger}$}
\thanks{$^{\dagger}$ Technion -- Israel Institute of Technology, Faculty of Mathematics, Technion City, Haifa 3200003, Israel. \textbf{chphdiez@campus.technion.ac.il}.}

\begin{abstract}
We establish a compactness theorem for the resolvent of free
Laplacians. More precisely, let $X=(X_1,\ldots,X_n)$ be a
self-adjoint tuple admitting Lipschitz conjugate variables
in the sense of Dabrowski (2014). Under the non-commutative curvature--dimension condition $CD(c,\infty)$ with $c>0$, introduced in our previous work (2026), we prove that the free Laplacian $\Delta:=\partial^*\bar\partial$
is proper, that is, has compact resolvent. Our argument is non-perturbative and uses a spectral variant of the free Bochner method based on Dabrowski's almost-commutation formulas.

A further stochastic argument gives a normal unital trace-preserving embedding of $M=W^*(X_1,\ldots,X_n)$ into $L(\mathbb F_2)$.
Together, the compactness and embedding theorems gives several operator-algebraic consequences, including the Haagerup approximation property, the CMAP, and strong solidity. We also obtain proper proximality and bi-exactness.

We finally obtain several geometric consequences. We prove a free analogue of Klartag's improved log--Lichnerowicz inequality (2023), together with rigidity in the equality case.
More surprisingly, under $CD(0,\infty)$, we prove a dimension-free Kannan--Lov\'asz--Simonovits (KLS) bound in the free setting. As a consequence, we obtain a free thin-shell estimate for isotropic tuples. We conclude with further questions on spectral comparison with the canonical semicircular model.
\end{abstract}

\maketitle

\section{Introduction}

\subsection{The classical compactness problem and its free counterpart}

The spectral theory of diffusion generators is a central theme in analysis, probability and geometry. One wants first to understand when the spectrum becomes discrete, since this analytic property already expresses strong confinement and provides a spectral decomposition of the underlying $L^2$-space. We recall the classical picture which serves as the guiding analogy for our work. Let $\mu$ be a probability measure on $\R^n$ with a smooth positive density $\rho=e^{-V}$. The weighted Laplace operator associated with $\mu$ is defined, initially on $C_c^\infty(\R^n)$, by
\[
\mathcal L_V f
=
\Delta f+\nabla(\log\rho)\cdot\nabla f
=
\Delta f-\nabla V\cdot\nabla f,
\]
and satisfies the integration-by-parts identity
\(
\int (\mathcal L_Vf)g\,d\mu=-\int\langle\nabla f,\nabla g\rangle\,d\mu.
\)
Thus $-\mathcal L_V$ is the positive Langevin generator associated with the diffusion $dX_t=-\nabla V(X_t)+\sqrt{2}dB_t$, where $(B_t)_{t\geq 0}$ is an $n$-dimensional Brownian motion. It defines a positive self-adjoint operator on $L^2(\mu)$. Smoothness and positivity of the density give the natural self-adjoint realization and identify the constants as its kernel, but they do not by themselves imply discreteness of the spectrum.
\par\vspace{0.5em}
The compactness problem is encoded by the following fundamental equivalence:
\[
(\alpha-\mathcal L_V)^{-1}\ \text{compact}
\quad\Longleftrightarrow\quad
\text{$-\mathcal L_V$ has discrete spectrum},
\]
for $\alpha>0$ (see~\cite{Davies,ReedSimonIV} for details). Equivalently, $-\mathcal L_V$ admits a complete orthonormal basis of eigenfunctions, every eigenvalue has finite multiplicity, and the eigenvalues tend to infinity. This is of course considerably stronger than the existence of a spectral gap, which only isolates the constants from the rest of the spectrum.
\par\vspace{0.5em}
Classically, compactness of the resolvent is intimately related to confinement at infinity. Indeed, it is well known that $-\mathcal L_V$ is unitarily equivalent to the Schr\"odinger operator (see Klartag-Putterman~\cite[Remark~2.7]{KlartagPutterman})
\[
-\Delta+\frac14\abs{\nabla V}^2-\frac12\Delta V.
\]
In particular, when $\rho$ is $\mathcal C^2$ and $ {\Delta (\sqrt{\rho})}/{\sqrt{\rho}}$ tends to infinity
at infinity, the resolvent is compact (we refer to Bakry--Gentil--Ledoux~\cite[Corollary~4.10.9]{Ledoux}). In the log-concave case, one has another natural sufficient condition
\(
V(x)/\abs{x}\rightarrow+\infty
\)
as $\abs{x}\to\infty$, as obtained by Klartag--Putterman~\cite[Appendix~A]{KlartagPutterman}. Their very elegant proof in the regular log-concave case, for rapidly decreasing log-concave densities, highlights this ``confinement'' property, which in fact admits an analytic formulation through an exterior coercivity estimate: for every $a>0$, there exists a sufficiently large ball such that
\[
\int_{\R^n}\abs{\nabla f}^2\,d\mu
\geq
a\int_{\R^n}f^2\,d\mu
\]
for every function $f$ supported outside this ball. Thus, very far from the center of the measure, the energy necessarily becomes arbitrarily large. Under the regularity assumptions of~\cite[Section~2]{Klartaglog}, the eigenfunctions are moreover \emph{\(\mu\)-tempered}: they and all their partial derivatives of arbitrary order have subexponential decay relative to $\rho$~\cite[Lemma~2.2]{Klartaglog}. Here decay relative to $\rho$ means that $\sqrt\rho\,f$ decays exponentially (using again the unitary equivalence with Schr\"odinger operators, see~\cite[Theorem~XIII.70]{ReedSimonIV}).

\par\vspace{0.5em}
Let us emphasize that this classical picture separates several properties and functional inequalities which differ in importance. A Poincar\'e inequality only gives a spectral gap. A logarithmic Sobolev inequality is equivalent, in the diffusion setting, to hypercontractivity, following the seminal work of Gross~\cite{gross}. Stronger Sobolev or Nash-type estimates may imply ultracontractivity and, for a finite invariant measure, compactness or trace-class properties of the heat semigroup (see~\cite{Ledoux,Davies} for more details). None of these implications can simply be reversed. In particular, compact resolvent alone gives neither hypercontractivity nor ultracontractivity. The theory explained above extends readily to the more general setting of weighted Riemannian manifolds $(M^n,g,\mu)$, where $(M^n,g)$ is a complete smooth Riemannian manifold and $\mu=e^{-V}d\mathrm{Vol}_g$, with completeness ensuring, for example, that the associated Laplacian $\Delta_{g,\mu}$ is essentially self-adjoint on an appropriate domain containing the space of compactly supported smooth functions (see Bakry--Gentil--Ledoux~\cite[Proposition~3.2.1]{Ledoux}). We recall, for example, that on
positively curved weighted manifolds, we always have discrete spectrum, since they satisfy a log-Sobolev inequality by the Bakry--\'Emery criterion (see Bakry--\'Emery's seminal contribution for further details~\cite{Bakry}) and such a functional inequality is known to imply discreteness; we refer to~\cite{Yang,Cianchi} for precise statements.
\par\vspace{0.5em}
Our purpose in this paper is to investigate the corresponding question for the \emph{free Laplacian} in divergence form associated to closable free difference quotients. Let $X=(X_1,\ldots,X_n)$ be a finite self-adjoint tuple in a tracial $W^*$--probability space, let $M=W^*(X_1,\ldots,X_n)$, and denote by $\bar\partial=(\bar\partial_1,\ldots,\bar\partial_n)$ the closed free difference quotient. The natural positive self-adjoint operator is then
\[
\Delta:=\partial^*\bar\partial
=
\sum_{j=1}^n\partial_j^*\bar\partial_j.
\] This operator is the generator of a completely Dirichlet form on $L^2(M)$. In perfect analogy with the preceding classical picture, compactness of $(\alpha+\Delta)^{-1}$ for some (equivalently, every) $\alpha>0$ is equivalent to compactness of the free Sobolev embedding
\(
\Dom(\bar\partial)=\Dom(\Delta^{1/2})\hookrightarrow L^2(M),
\)
when the form domain is equipped with its graph norm. It is also equivalent to the fact that the spectrum of $\Delta$ consists only of eigenvalues of finite multiplicity tending to infinity. The problem is therefore to understand whether the free differential calculus itself contains enough confinement to force this compactness. As we will see in the sequel, the required confinement is naturally provided by a convexity condition.
\par\vspace{0.5em}
For a finite semicircular family, everything is particularly simple and explicit: $\Delta$ is the number operator on the full Fock space (see Biane--Speicher~\cite{BS}). Its eigenvalue $k$ has multiplicity $n^k$, and hence its resolvent is compact. In the infinite-dimensional setting the situation is already different, since the eigenspaces become infinite-dimensional: the Wigner chaoses are the orthogonal eigenspaces (parallel to the classical world of Wiener chaoses). Beyond the finite semicircular model, however, very little seems to be known about the spectrum of the ordinary free Laplacian. In our previous work on rigidity for the free Poincar\'e inequality~\cite{DiezObata}, compactness of the resolvent appeared precisely as the missing spectral hypothesis needed to complete the picture of our free Obata theorem. More recently, Jekel and Basa Landry~\cite{JekelLandry} emphasized that compact resolvent is not known in general even under strong analytic regularity assumptions such as finite free Fisher information or Lipschitz conjugate variables. This question can in fact be traced back to the 2006 AIM workshop on Free Analysis~\cite{Aim}, where compactness of the heat semigroup associated with the free Laplacian was asked about for limiting non-commutative laws arising from random matrix models, already in the simplest analytic perturbations of the quadratic model. The difficulty lies in the passage from finite-dimensional matrix models to the limiting free Laplacian, which is far from immediate.
\par\vspace{0.5em}
There is also a distinction that has to be made between our setting and the $q$--Gaussian case equipped with its $q$--Gaussian number operator (we do not know whether such tuples generally belong to our $CD$ class for a large range of $q$). For a finite $q$--Gaussian family (these are deformations of the free group factors) introduced by Bo\.{z}ejko, K\"ummerer and Speicher~\cite{BSK97}, $-1<q<1$, the $q$--Ornstein--Uhlenbeck operator
\(
D_q=\sum_{k\geq0}kP_k
\)
on the $q$--Fock space, where $P_k$ denotes the projection onto the $k$--th homogeneous chaos, has compact resolvent and its heat semigroup is ultracontractive by the important result of Bo\.{z}ejko~\cite{Bozejko}. At $q=0$, $D_q$ agrees with the ordinary free Ornstein--Uhlenbeck operator, which is precisely the free Laplacian associated with a semicircular family. For $q\neq0$, however, the two operators are well known to be different (this is even clearer in the infinite-dimensional setting, as highlighted in~\cite[Definition~20]{Diez2}, where the underlying derivation is not even Hilbert-Schmidt for any range of $q$; this creates, for example, numerous complications in building a robust $q$-version of the Malliavin calculus). Therefore, the spectral properties of the number-length operator do not answer the compactness question for
\(
\sum_j\partial_j^*\bar\partial_j
\)
associated with the canonical $q$--Gaussian generators. This distinction is important and the reader should keep these results in mind when comparing our results with the literature on these semigroups.
\par\vspace{0.5em}
There are finally some subtle points which have in fact deep consequences. First, the operator $\Delta$ depends on the chosen generating tuple and on its free differential calculus, rather than only on the abstract von Neumann algebra $M=W^*(X_1,\ldots,X_n)$. Already a change of coordinates modifies the free difference quotients through the chain rule, and therefore also the corresponding non-commutative Jacobian and Laplacian. What can be reasonably preserved, or at least compared on the spectral side, thus requires a sufficiently regular change of variables with two-sided bounds. This explains why the analytic invertibility arising in the free monotone transport of Guionnet and Shlyakhtenko~\cite{GS2} contains considerably more information than the conclusion $W^*(X)\cong L(\mathbb F_n)$ alone: the two Dirichlet forms become comparable, therefore giving compact resolvent (and even ultracontractivity).

\subsection{Positive free curvature and the spectral Bochner argument}
We now describe the curvature assumption and the spectral mechanism behind our proof. Suppose that $X$ admits Lipschitz conjugate variables $\xi=(\xi_1,\ldots,\xi_n)$ in the sense of Dabrowski~\cite[Definition~1]{Dab14}. Their non-commutative Jacobian is
\[
\mathscr J\xi=(\bar\partial_j\xi_i)_{i,j=1}^n\in M_n(M\bar\otimes M^{op}).
\]
Following our previous work on Obata's theorem in free probability~\cite{DiezObata}, we say that $X$ satisfies the curvature--dimension criterion $CD(c,\infty)$ with $c\geq 0$ when
\(
\mathscr J\xi\geq c.(1\otimes 1)\otimes I_n
\)
in the $C^*$--order of $M_n(M\bar\otimes M^{op})$. The canonical examples of such $CD$ tuples are in fact provided by free Gibbs laws. For a non-commutative polynomial potential $V$ (usually self-adjoint) and a realization of its law $\tau_V$ by a bounded self-adjoint tuple $X$, the conjugate variables are obtained by evaluating the cyclic gradient $\mathscr DV(X)$, and their Jacobian becomes exactly the non-commutative Hessian $\mathscr J\mathscr DV(X)$. Our condition is, however, much more general: it is evaluated only at the realized tuple and does not require the existence of any potential. This flexibility is one of the principal advantages of working directly with Lipschitz conjugate variables: the condition is intrinsic to the chosen non-commutative differential calculus and does not require global estimates valid in arbitrary $C^*$--algebras. The picture for free Gibbs laws and that for Lipschitz conjugate variables would, however, be much closer if Voiculescu's \emph{cyclic-gradient conjecture} were true: this conjecture asks whether conjugate variables can always be approximated in $L^2$ by cyclic gradients (see \cite[Section~0.2, Q.4]{Aim} for precise details). It is believed that the answer should be yes, at least for such Lipschitz conjugate variables, as highlighted by Mai and Speicher~\cite[Remark~1.2]{MS}.
\par\vspace{0.5em}
Our main result shows that strictly positive free curvature implies that the corresponding free Laplacian is \emph{proper}: it has compact resolvent. Thus its spectrum is discrete, every eigenvalue has finite multiplicity, and the eigenvalues can accumulate only at infinity. The curvature bound also gives the expected spectral gap: below the curvature lower bound, the only spectral vectors are in fact the constants. More importantly, the same estimate allows us to control, by induction, every bounded spectral subspace.
\par\vspace{0.5em}
A related approach has recently been developed by Chandra, Gordina and Peev~\cite{CGP}, who introduce free Bakry--\'Emery calculus for free Langevin dynamics associated
with polynomial potentials. Using non-commutative
carr\'e du champ operators and their iterates, they obtain, under their curvature and analytic assumptions, gradient
estimates and a modified logarithmic Sobolev inequality for the stationary law. Their approach and ours share in common the use of Dabrowski's almost commutation formulas.
\par\vspace{0.5em}
The proof is based on a spectral version of the free Bochner argument developed in our previous work~\cite{DiezObata}, using Dabrowski's \emph{almost-commutation formulas}~\cite{Dab14} between free difference quotients and the associated Laplacian. The key observation is that, on a positive spectral subspace below a level, the Bochner estimate controls the tensor energy of the free difference quotient gradient by the difference between this level and the curvature. Projecting the gradient onto the tensor spectral subspace below this lower level then gives an injective map. Indeed, Voiculescu's free Poincar\'e inequality~\cite{Voipub} shows that the kernel of the closed free difference quotient gradient consists of the constants. This provides a free counterpart of the confinement estimates used in the classical setting.
The second ingredient is the tensor structure of the free difference quotient. The tensor Laplacian appearing in Dabrowski's almost-commutation formula~\ref{dab} is the sum of the two positive copies of \(\Delta\) acting on each one of the tensor parts. Its low-energy spectral subspaces are therefore controlled by tensor products of lower scalar spectral subspaces. Since the tuple has only finitely many generators, the free difference quotient gradient also has only finitely many components. If we start below the curvature lower bound, we may then proceed inductively through successive energy intervals, $[kc,(k+1)c),\: k\in \mathbb Z_{\geq 0}$. This proves that every bounded spectral subspace is finite-dimensional and gives an explicit recursion for the spectral counting function (despite the estimate at this stage being very rough).
\par\vspace{0.5em}
A more precise use of the same tensor spectral calculus gives a sharper estimate. Instead of counting all pairs of spectral values in a square, we keep only those belonging to the triangle determined by the condition that their sum remains below the prescribed level. This leads to a convolution inequality for the spectral counting measure and to the Catalan recursion. We obtain in this way an exponential bound for the spectral growth, together with a trace-class estimate for the heat semigroup. The preceding bound is still not expected to be sharp. We have not,
however, been able to push the argument far enough to remove the
Catalan numbers. The semicircular case has an additional very rigid structure
which can be explained by ``slicing'' Dabrowski's almost-commutation formula using the
conditional expectation \(
\tau\otimes\mathrm{id}
\).
This leads to the exact
multiplicities \(n^k\), instead of the Catalan multiplicities appearing
in our general argument. We do not know for now how to adapt this method in the general case.
\subsection{Consequences and applications}
Under the same positive-curvature assumption, we also construct a normal unital trace-preserving embedding of
$M=W^*(X_1,\ldots,X_n)$ into $L(\mathbb F_2)$, using Dabrowski's
stationary construction for a free stochastic partial differential
equation~\cite{Dab14}. This argument uses only the $CD(c,\infty)$ criterion and does not rely on compactness of the resolvent. The compactness
theorem, in turn, has several consequences beyond discreteness of the spectrum, connecting the spectral properties of the free Laplacian
with the structure of the von Neumann algebra generated by the tuple.
The two approaches are complementary: compactness provides approximation maps and allows us to apply proper-derivation methods introduced by various authors~\cite{OP1,OP2,DKEP,DP},
while the embedding yields additional approximation and structural properties that we do not obtain from compactness alone.
\par\vspace{0.5em}
The first direct consequence of compactness is the \emph{Haagerup
approximation property}, introduced by Haagerup~\cite{Haag} for free
groups and developed by Choda~\cite{Choda} in the von Neumann algebra
setting. Let us recall that for a discrete group von Neumann algebra
$L(\Gamma)$, this property is equivalent to Gromov's a-T-menability
of $\Gamma$, that is, to the existence of a $1$-cocycle
$b:\Gamma\to H$ with values in a unitary representation of $\Gamma$
such that $\gamma\mapsto\|b(\gamma)\|_H$ is proper.
In our situation the proof is immediate: the resolvents of the
completely Dirichlet form are normal unital completely positive
trace-preserving maps, their $L^2$--implementations are compact,
and they converge strongly to the identity. Thus compactness of
the resolvent gives this approximation property directly.

\par\vspace{0.5em}
We also return to the stationary solutions of the free stochastic
partial differential equation constructed by Dabrowski~\cite{Dab14}.
His construction gives normal unital trace-preserving embeddings
of $M=W^*(X_1,\ldots,X_n)$ into the free product with a free
Brownian algebra, i.e. $M*L(\mathbb F_\infty)$. Combining his stochastic energy identity with the
resolvent form of almost commutation, also used in our previous
work~\cite{DiezObata}, we show that $CD(c,\infty)$ with $c>0$ provides sufficient
integrability to pass to the limit in the mild equation and obtain a normal unital trace-preserving embedding \(
M\hookrightarrow L(\mathbb F_\infty),
\) and consequently into $L(\mathbb F_2)$.
This argument using positive curvature directly is of independent interest as it does not require using
compactness of the resolvent. In particular, $M$ is
$R^\omega$--embeddable, where $R$ denotes the hyperfinite
$\mathrm{II}_1$ factor.
The embedding then also gives the \emph{weak$^*$ complete metric
approximation property} (CMAP) and weak exactness, and another direct way to deduce the Haagerup property.
\par\vspace{0.5em}
A further consequence of compactness uses the coarse
structure of the derivation more deeply. The free difference quotient takes
values in a finite direct sum of coarse correspondences. Combining
this fact with compactness of the associated free Laplacian and
the proper-derivation method of Ozawa--Popa~\cite{OP1,OP2}, we obtain
a weakly compact-normalizer statement and, as a consequence,
\emph{solidity} of $M$. Under the weak$^*$ completely bounded
approximation property, Ozawa's theorem~\cite{OzawaWeak} upgrades
the same argument to \emph{strong solidity}, introduced in the
work of Ozawa--Popa~\cite{OP1,OP2}. This additional approximation
hypothesis is now automatic in our positive-curvature setting,
since the embedding gives CMAP. Thus $M$ is strongly solid, as
are its unital von Neumann subalgebras.
For $n\geq2$, Dabrowski's results~\cite{Dab10} show that $M$ is
a non-$\Gamma$ $\mathrm{II}_1$ factor (more precisely, it has a non-gamma set in the sense of Peterson~\cite{Peterson}: it satisfies a spectral gap, as obtained by Amrutam--Jekel--Wasilewski~\cite[Theorem~3.8]{jek3}). It follows that $M$ is prime and has
no Cartan subalgebra. Primeness was already obtained by
Dabrowski--Ioana~\cite{Dabioa}, under the Lipschitz assumption used here but without assuming curvature. The additional
conclusion here is the strong solidity statement.
\par\vspace{0.5em}
The last operator-algebraic applications combine compactness of the
free Laplacian with the work of
Ding--Kunnawalkam Elayavalli--Peterson~\cite{DKEP} and
Ding--Peterson~\cite{DP}, which extends the notions of proper
proximality and bi-exactness  from groups to von Neumann algebras
using the \emph{small-at-infinity} operator system.
Their results apply directly in our setting. Indeed, the free
difference quotient is a closable real derivation into a finite
direct sum of coarse correspondences, while compactness of the
resolvent of $\Delta$ implies that, for every $\alpha>0$, both
\(
 \rho_\alpha:=\alpha(\alpha+\Delta)^{-1}
 \) and 
 \(\zeta_\alpha:=\rho_\alpha^{1/2}
             =\alpha^{1/2}(\alpha+\Delta)^{-1/2}
\)
are compact, i.e. belong to $\mathbb K(L^2(M))$. For $n\geq2$, $M$ is a non-amenable
factor, and~\cite[Proposition~9.4 and Lemma~9.5]{DKEP} therefore
give proper proximality. Under weak exactness, compactness also
gives bi-exactness  by~\cite[Theorem~7.11]{DP}. Since weak exactness
follows from the embedding, $M$ is biexact under our positive
curvature assumption, without any additional hypothesis.
\par\vspace{0.5em}
We return in the last section, Section~\ref{Sec3}, to the geometric motivation of the paper. Compactness of the resolvent ensures that the bottom of the positive spectrum is attained, which was the missing point in our previous work on Obata's rigidity phenomenon in free probability~\cite{DiezObata}, where we had to assume that the free Poincar\'e inequality was saturated at the curvature bound. This rigidity phenomenon is well known in Gamma--calculus and was already pointed out by Ledoux~\cite[Proposition~4]{Ledoux2} and Bouyrie~\cite[Lemma~5.1]{Bouyrie} and has motivated our previous work, which we briefly recall here: under \(CD(1,\infty)\) (the normalization $1$ plays no particular role here), saturation of the free Poincar\'e inequality therefore produces an eigenfunction which, after an orthogonal change of the centered variables, leads, by our rigidity theorem~\cite[Theorem~3.15 and Corollary~3.17]{DiezObata}, to a free product decomposition with a free group factor, whose rank is given by the dimension of the first non-zero eigenspace.
\par\vspace{0.5em}
A further important consequence of the discrete spectrum is a free counterpart of Klartag's improved \emph{log--Lichnerowicz} inequality obtained in his recent work on the Kannan-Lov\'asz-Simonovits (KLS) conjecture~\cite{Klartaglog} and which therefore improves our previous log--Lichnerowicz estimate obtained in~\cite[Theorem~2.20]{DiezObata}.
The proof follows Klartag's ``preferred-direction'' mechanism (also of importance in the ``hot-spots'' conjecture in Euclidean geometry~\cite{Judge}), with the tensor spectral gap for \(\Delta^\otimes\) replacing its classical counterpart. We are also able to identify the equality case, thus providing a counterpart to Courtade and Fathi's result~\cite{courtfa}: saturation of the improved estimate implies saturation of the usual free Lichnerowicz inequality, and our free Obata theorem then shows that a semicircular factor splits off in a direction of maximal variance. These analogies naturally lead us to a free counterpart of the
Kannan--Lov\'asz--Simonovits conjecture~\cite{KLS}, whose classical
formulation was recently resolved by Bizeul, Klartag and
Lehec~\cite{bizkla} and by Song and Zhang~\cite{song}. Bizeul, Klartag and
Lehec~\cite{bizkla} used a first version of a preprint of Song and Zhang~\cite{song} which subsequently solve it themselves. The conjecture  asks in fact (in one of its formulations) whether the Poincar\'e constant of a log-concave probability measure is bounded, up to a universal constant, by the operator norm of its covariance matrix (see Klartag and Lehec~\cite{klalehec} and the survey of Lee and Vempala~\cite{lee} for precise accounts of the conjecture and its formulations). 
Quite surprisingly, we prove that, under $CD(0,\infty)$, the free Poincar\'e constant is bounded
by a universal multiple of the operator norm of the covariance
matrix, independently of the number of variables. 
The proof combines the free Bochner formula with Poisson equations
for the partial traces of the free gradient, obtained by composing the conditional expectation with the free difference quotients $(\mathrm{id}\otimes \tau)\circ \partial_i$. Co-associativity allows us
to use both tensor factors to obtain the strict improvement needed
for this comparison. As a consequence, we also obtain a free
thin-shell estimate: for a centered isotropic tuple, the variance
of the sum of squares is bounded by a universal multiple of the
number of variables.
\par\vspace{0.5em}
We finally discuss a further spectral question, motivated by the sharp spectral comparison theorem of E. Milman~\cite{MilmanSpectral} in the classical setting, which shows that all the ordered eigenvalues of a weighted Euclidean space equipped with the standard Euclidean metric \((\mathbb R^n,\lvert\cdot\rvert,\mu)\) satisfying \(CD(c,\infty)\) are bounded from below by the corresponding eigenvalues of the Gaussian model. We prove a one-dimensional free counterpart. Note that, classically, the Euclidean norm plays a fundamental role; a counterexample for a non-Euclidean metric on the plane was recently found by Dudarov--Mikulincer~\cite[Theorem~2]{Mikulin}.
\section{Preliminaries}\label{prelim}

\subsection{Free difference quotients and the associated Dirichlet form}

Throughout this paper, $(\mathcal M,\tau)$ will denote a tracial $W^*$--probability space, which means that $\mathcal M$ is a finite von Neumann algebra endowed with a faithful normal tracial state $\tau$. We consider a finite self-adjoint tuple $X=(X_1,\ldots,X_n)\in\mathcal M_{sa}^n$, and we set $M=W^*(X_1,\ldots,X_n)$ and $M_0=\C\langle X_1,\ldots,X_n\rangle$. Whenever the free difference quotients are considered, the tuple is always assumed algebraically free.
\newline
We write $L^2(M):=L^2(M,\tau)$, with the convention $\ip{x}{y}_2:=\tau(y^*x)$ and $\norm{x}_2:=\ip{x}{x}_2^{1/2}$. We will use without further mention the same notation for the canonical inner products and norms on $L^2(M^{op})$, on $L^2(M)\bar\otimes L^2(M^{op})$, and on their finite Hilbert direct sums.
\newline
We continue to write 
$\tau(x) := \langle x,1\rangle_{2}$
for the vector functional associated with the cyclic vector $1$, which extends the tracial state from $\mathcal M$ to all of $L^2(\mathcal M)$ by continuity. The same applies to the tensor extension $\tau\otimes\tau$ of the trace on $L^2( M)\bar\otimes L^2(M^{op})$.
\newline
Finally, unless otherwise indicated, $\lVert\cdot\rVert$ denotes the operator norm on some $\mathbb B(\mathcal H)$.

\begin{definition}[Partial free difference quotients]
The \emph{partial free difference quotients} are defined as the derivations
\[
\partial_i:M_0\to L^2(M)\bar\otimes L^2(M^{op}),
\]
which are determined on the generators by
\(
\partial_i(X_j)=\delta_{ij}(1\otimes1).
\)
\newline
We recall that the free difference quotients are \emph{real} derivations: for every $x\in M_0$,
\[
\partial_i(x^*)=(\partial_i(x))^\dagger,
\]
where $(a\otimes b)^\dagger:=b^*\otimes a^*$.
\end{definition}

\begin{remark}\label{real}
Suppose that each $\partial_i$ is closable, viewed as an unbounded operator
\(
L^2(M)\to L^2(M)\bar\otimes L^2(M^{op}),
\)
and denote by $\bar\partial_i$ its closure. The row operator $\partial=(\partial_1,\ldots,\partial_n)$ on $M_0$ is also closable; throughout the paper, $\bar\partial$ denotes its closure. Thus $M_0$ is a core for $\bar\partial$, and
\[
\Dom(\bar\partial)\subseteq\bigcap_{i=1}^n\Dom(\bar\partial_i),
\qquad
\bar\partial x=(\bar\partial_1x,\ldots,\bar\partial_nx),
\quad x\in\Dom(\bar\partial).
\]
The real property of the free difference quotients then gives, by the usual double-adjoint argument, that $\Dom(\bar\partial)$ is $^*$--stable and
\[\label{adjunc}
\bar\partial_i(x^*)=(\bar\partial_i(x))^\dagger,
\qquad x\in\Dom(\bar\partial).
\]
We refer to Dabrowski~\cite{Dab10} for the complete argument.
\end{remark}

\begin{definition}[{Free Fisher information, Voiculescu~\cite[Definition~6.1]{Voic}}]
We define the \emph{free Fisher information} of $(X_1,\ldots,X_n)$ by
\[
\Phi^*(X_1,\ldots,X_n)
:=\sum_{i=1}^n\norm{\partial_i^*(1\otimes1)}_2^2,
\]
whenever each $\partial_i^*(1\otimes1)$ exists in $L^2(M)$; otherwise, we set $\Phi^*(X)=+\infty$.
\end{definition}

\begin{proposition}[{Closability, Voiculescu~\cite[Lemma~3.3]{Voic}}]
Whenever $\Phi^*(X_1,\ldots,X_n)<\infty$, each $\partial_i$ is closable as an unbounded operator
\[
L^2(M)\to L^2(M)\bar\otimes L^2(M^{op}).
\]
We denote again by $\bar\partial_i$ its closure, and by $\bar\partial$ the closure of the row operator $\partial=(\partial_1,\ldots,\partial_n)$ on $M_0$, as in Remark~\ref{real}.
\end{proposition}

\begin{definition}[{Dabrowski~\cite[Definition~1]{Dab14}}]\label{lip}
Let $M=(W^*(X_1,\ldots,X_n),\tau)$ be a tracial von Neumann algebra. We say that $M$ satisfies the \emph{Lipschitz conjugate variable condition} whenever:
\begin{enumerate}
\item the partial free difference quotients $\partial_i$ are densely defined and closable,
\item the conjugate variables $\xi_i:=\partial_i^*(1\otimes1)$ exist in $L^2(M)$ for every $i=1,\ldots,n$, belong to $\Dom(\bar\partial)$, and satisfy
\[
\bar\partial_j\partial_i^*(1\otimes1)
\in M\bar\otimes M^{op}
\subset L^2(M\bar\otimes M^{op})
\cong L^2(M)\bar\otimes L^2(M^{op})
\cong L^2(M)\bar\otimes L^2(M)
\]
for every $i,j=1,\ldots,n$.
\end{enumerate}
\end{definition}

We now recall a useful consequence of the Lipschitz conjugate variable condition, first observed by Dabrowski in~\cite[Sketch of the proof of Corollary~25~(a')]{Dab14} and also used in our previous work~\cite{DiezObata}: the first- and second-order conjugate variables are in fact bounded. This fact will be used several times below.

\begin{proposition}\label{prop5}
Assume that $M=(W^*(X_1,\ldots,X_n),\tau)$ satisfies the Lipschitz conjugate variable condition. We then have:
\begin{enumerate}
\item $\xi_i:=\partial_i^*(1\otimes1)\in M$ for every $i$,
\item $\xi_i^{(2)}:=\partial_i^*(1\otimes\partial_i^*(1\otimes1))=\partial_i^*(1\otimes\xi_i)=\partial_i^*(\xi_i\otimes1)\in M$ for every $i$,
\item $(\partial_i\otimes\mathrm{id})\circ\partial_j:M_0=\A\to L^2(M)^{\bar\otimes3}$ is densely defined and closable for every $i,j=1,\ldots,n$.
\end{enumerate}
\end{proposition}

The closure of the free difference quotient gives rise to a completely Dirichlet form. We now collect the properties of its generator, as well as the resolvent approximations.

\begin{proposition}[Cipriani--Sauvageot~\cite{Sauvage,CiprianiSauvageot2003}, Peterson~\cite{Peterson}, Dabrowski~\cite{Dab10}]\label{resolvent}
Let $M=W^*(X_1,\ldots,X_n)$ be the tracial von Neumann algebra generated by self-adjoint variables $(X_1,\ldots,X_n)$ in a tracial $W^*$--probability space $(\mathcal M,\tau)$, and suppose that the partial free difference quotients are closable. Then the following assertions hold:
\begin{enumerate}
\item The operator
\(
\Delta:=\partial^*\bar\partial,
\)
where $\bar\partial=(\bar\partial_1,\ldots,\bar\partial_n)$, is the generator of a completely Dirichlet form on $L^2(M)$; see~\cite{Sauvage}. In particular, all its matrix amplifications are again Dirichlet.

\item The resolvent maps
\(
\rho_\alpha:=\alpha(\alpha+\Delta)^{-1},
\)
$\alpha>0$, are normal, unital, tracial $(\tau\circ\rho_\alpha=\tau)$, completely positive, and contractive on both $L^2(M)$ and $M$; see, for instance,~\cite{Ma}. We also have
\[
\norm{x-\rho_\alpha(x)}\leq2\norm{x},
\qquad
\norm{x-\rho_\alpha(x)}_2\to0
\quad(\alpha\to\infty)
\]
for every $x\in M$, as follows from~\cite[Proposition~2.5]{CiprianiSauvageot2003}.

\item The semigroup $\varphi_t=e^{-t\Delta}$ is strongly continuous with generator $-\Delta$, and
\[
\rho_\alpha
=\alpha\int_0^\infty e^{-\alpha t}\varphi_t\,dt,
\]
where the integral is understood with respect to the strong operator topology on $L^2(M)$.

\item The square-root resolvent
\[
\zeta_\alpha:=\rho_\alpha^{1/2}
=\pi^{-1}\int_0^\infty\frac{t^{-1/2}}{1+t}\,
\rho_{\alpha(1+t)/t}\,dt
\]
satisfies
\[
\Ran(\rho_\alpha)=\Dom(\Delta)\subset\Dom(\bar\partial),
\qquad
\Ran(\zeta_\alpha)=\Dom(\Delta^{1/2})=\Dom(\bar\partial),
\]
and therefore $\bar\partial\circ\zeta_\alpha$ is bounded (see~\cite{Peterson}).
\end{enumerate}
\end{proposition}

\begin{flushleft}
We will also use throughout the paper Voiculescu's fundamental free Poincar\'e inequality. Let us emphasize that it holds for every bounded self-adjoint tuple, with no assumption on its law.
\end{flushleft}

\begin{lemma}[Free Poincar\'e inequality, Voiculescu~\cite{Voipub}]\label{voipub}
Let $Y\in\Dom(\bar\partial)$ be self-adjoint or non-self-adjoint. Then there exists a constant $C>0$, which depends on the $X_i$ but not on $Y$, such that
\[
\norm{Y-\tau(Y).1}_2^2
\leq
C\sum_{j=1}^n\norm{\bar\partial_jY}_2^2.
\]
Moreover, if $C_P:=C_P(X)$ denotes the optimal constant, called the \emph{free Poincar\'e constant}, then
\[
C_P\leq4n\max_{i=1,\ldots,n}\norm{X_i}^2.
\]
\end{lemma}

A direct consequence of the free Poincar\'e inequality is that it identifies the kernel of the closed free difference quotient, and therefore also the zero eigenspace of the free Laplacian.

\begin{proposition}\label{kerprop}
Suppose that the partial free difference quotients are closable. We have
\begin{equation}\label{kernel}
\Ker(\Delta)=\Ker(\bar\partial)=\mathbb C\,1.
\end{equation}
\end{proposition}

Let $ O:L^2(M)\to L^2(M^{op})$ be the canonical unitary $ O(x)=x^{op}$ coming from traciality. In the sequel, and by a slight abuse of notation, the same symbol $\Delta$ denotes $ O\Delta O^{-1}$ on $L^2(M^{op})$, and $E(B)$ denotes the corresponding transported spectral projection on $L^2(M^{op})$. We use the same convention below when dealing with Dabrowski's almost-commutation formulas~\ref{dab} (see~\cite{Dab14,DiezObata} for more details).

\begin{proposition}[Tensor spectral calculus]\label{tensor}
The positive self-adjoint operators $\Delta\otimes \mathrm{id}$ and $\mathrm{id}\otimes\Delta$ on $L^2(M)\bar\otimes L^2(M^{op})$ strongly commute. If $E(\cdot)$ denotes the spectral measure of $\Delta$, with our preceding convention on the $L^2(M^{op})$ leg, their joint spectral measure $\mathbb E$ is characterized by
\[
\mathbb E(B\times C)=E(B)\otimes E(C),
\qquad
B,C\subset[0,\infty)\ \text{Borel}.
\]
The operator
\begin{equation}\label{delta2}
\Delta^\otimes
:=\int_{[0,\infty)^2}(\lambda+\mu)\,d\mathbb E(\lambda,\mu)
\end{equation}
is the closure of $\Delta\otimes \mathrm{id}+\mathrm{id}\otimes\Delta$ on $\Dom(\Delta)\otimes\Dom(\Delta)$ (where $\otimes$ denotes here the algebraic tensor product). Moreover, for every $s\geq0$, we have
\begin{equation}\label{tensorproj}
\mathds1_{[0,s]}(\Delta^\otimes)
\leq
\mathds1_{[0,s]}(\Delta)\otimes\mathds1_{[0,s]}(\Delta).
\end{equation}
Consequently, if $N_\Delta(s):=\rank\mathds1_{[0,s]}(\Delta)<\infty$, we get
\[
\rank\mathds1_{[0,s]}(\Delta^\otimes)
\leq
N_\Delta(s)^2.
\]
Finally,
\begin{equation}\label{tensorker}
\Ker(\Delta^\otimes)
=\Ker(\Delta)\bar\otimes\Ker(\Delta)
=\mathbb C.(1\otimes 1).
\end{equation}
\end{proposition}

\begin{proof}
The spectral projections of $\Delta\otimes \mathrm{id}$ and $\mathrm{id}\otimes\Delta$ are $E(B)\otimes \mathrm{id}$ and $\mathrm{id}\otimes E(C)$, and they therefore commute strongly. The joint spectral theorem gives the preceding identity together with the operator in~\eqref{delta2}. It remains first to verify the assertion concerning the core. For $n\geq 0$, put 
\(
Q_n:=E([0,n])\otimes E([0,n]).
\)
For every $x\in\Dom(\Delta^\otimes)$, the joint spectral calculus gives $Q_nx\to x$ in the graph norm of $\Delta^\otimes$. On $\Ran(Q_n)$, this graph norm is equivalent to the Hilbert-space norm, and 
\[
E([0,n])L^2(M)\otimes E([0,n])L^2(M^{op})
\]
is dense in $\Ran(Q_n)$ and contained in $\Dom(\Delta)\otimes\Dom(\Delta)$. Hence, the algebraic tensor product is indeed a core.
\newline
We now verify the projection estimate and set
\[
T_s:=\{(\lambda,\mu)\in[0,\infty)^2:\lambda+\mu\leq s\}.
\]
Then
\[
\mathds1_{[0,s]}(\Delta^\otimes)=\mathbb E(T_s),
\qquad
\mathds1_{[0,s]}(\Delta)\otimes\mathds1_{[0,s]}(\Delta)
=\mathbb E([0,s]^2).
\]
Since $0\leq \lambda+\mu\leq s$ implies $0\leq \lambda,\mu\leq s$, we have the immediate inclusion $T_s\subset[0,s]^2$, and the monotonicity of projection-valued measures yields~\eqref{tensorproj}. The range of the projection in the right-hand side of~\eqref{tensorproj} has dimension $N_\Delta(s)^2$, proving the desired rank estimate. Finally, since $\lambda,\mu\geq0$, the equality $\lambda+\mu=0$ is equivalent to $\lambda=\mu=0$, and~\eqref{tensorker} follows immediately from~\eqref{kernel}. This concludes the proof.
\end{proof}

\begin{definition}\label{rightaction}
Let $(M,\tau)$ be a tracial $W^*$--probability space. We first define on the algebraic tensor product $L^2(M)\otimes L^2(M^{op})$, for $x\otimes y\in M\otimes M^{op}$ and $a\otimes b\in L^2(M)\otimes L^2(M^{op})$,
\[
(a\otimes b)\sharp(x\otimes y):=ax\otimes by.
\]
Here the product on the second tensor leg is taken in $M^{op}$. By bilinearity and continuity, this extends to a bounded right action of $M\bar\otimes M^{op}$ on $L^2(M)\bar\otimes L^2(M^{op})$.

For $n\in\N$ and $T=(T_{ij})_{i,j=1}^n\in M_n(M\bar\otimes M^{op})$, we then define the \emph{right-leg} multiplication
\[
\mathcal R_T:
(L^2(M)\bar\otimes L^2(M^{op}))^{\oplus n}
\to
(L^2(M)\bar\otimes L^2(M^{op}))^{\oplus n}
\]
by
\[
(\mathcal R_T\eta)_i:=\sum_{j=1}^n\eta_j\sharp T_{ji},
\qquad i=1,\ldots,n,
\]
for $\eta=(\eta_1,\ldots,\eta_n)\in(L^2(M)\bar\otimes L^2(M^{op}))^{\oplus n}$.
\end{definition}

\begin{flushleft}
We now recall the useful representation properties of $\mathcal R$ obtained in our previous work~\cite[Proposition~2.18]{DiezObata}.
\end{flushleft}

\begin{proposition}\label{rightprop}
For each $n\geq1$ and $T,S\in M_n(M\bar\otimes M^{op})$, the operators introduced in Definition~\ref{rightaction} enjoy the following properties:
\begin{enumerate}
  \item \emph{Linearity:} $\mathcal R_{T+S}=\mathcal R_T+\mathcal R_S$.
  \item \emph{Boundedness:}
  \[
  \norm{\mathcal R_T}\leq \norm{T}_{M_n(M\bar\otimes M^{op})}.
  \]
  Hence
  \[
  \mathcal R_T\in B\!\left((L^2(M)\bar\otimes L^2(M^{op}))^{\oplus n}\right)
  \cong M_n\!\left(B(L^2(M)\bar\otimes L^2(M^{op}))\right).
  \]
  \item \emph{$^*$--structure:} $\mathcal R_{T^*}=(\mathcal R_T)^*$, where $(T^*)_{ij}=T_{ji}^*$ for the usual $^*$--involution on $M\bar\otimes M^{op}$.

\item $\mathcal R_{(1\otimes 1)\otimes I_n}=\mathrm{id}_{(L^2(M){\bar\otimes}\,L^2(M^{op}))^{\oplus n}}$.
  \item \emph{Multiplication rule:} Writing $TS$ for the matrix product formed using the multiplication $\sharp$ in $M\bar\otimes M^{op}$, we have

\[
  \mathcal R_S\circ\mathcal R_T=\mathcal R_{T S}=\mathcal R_{T\sharp S}
  \qquad\text{for all }T,S\in M_n(M\bar\otimes M^{op}).
  \]
In particular, whenever $T$ is invertible in $M_n(M\bar\otimes M^{op})$, the operator $\mathcal R_T$ is invertible and $(\mathcal R_T)^{-1}=\mathcal R_{T^{-1}}$.

    \item \emph{Positivity:}
  if $T\ge 0$ in $M_n(M\bar\otimes M^{op})$ (resp.\ $T>0$, i.e.\ $T\ge c.(1\otimes 1)\otimes I_n$ for some $c>0$), then 
  $\mathcal R_T\ge 0$ (resp.\ $\mathcal R_T\ge c\,\mathrm{id}$) as an operator on 
  $(L^2(M)\bar\otimes L^2(M^{op}))^{\oplus n}$.
\end{enumerate}
Thus, $T\mapsto\mathcal R_T$ is a $*$--anti-representation of $M_n(M\bar\otimes M^{op})$, or equivalently a $*$--representation of the opposite algebra $M_n(M\bar\otimes M^{op})^{op}$.
\end{proposition}

\begin{flushleft}
We also need the following important symmetry property for the free Jacobian of the conjugate variables.
\end{flushleft}
\begin{proposition}\label{proplip}
Under Definition~\ref{lip}, the matrix
\[
\mathscr J\xi
=(\bar\partial_j\xi_i)_{i,j=1}^n
\in M_n(M\bar\otimes M^{op})
\]
is self-adjoint.
\end{proposition}

\begin{proof}
This was proved in~\cite[Lemma~2.19]{DiezObata}. We nevertheless recall the short argument for the convenience of the reader.
\newline
We denote by $\dagger$ the isometric Hilbert--Schmidt involution on $M\otimes M^{op}$, given on simple tensors by
\(
(a\otimes b^{op})^\dagger=b^*\otimes(a^*)^{op},
\)
and let
\(
\sigma(a\otimes b^{op})=b\otimes a^{op}
\)
denote the flip. The real property of the free difference quotients first gives
\(
(\bar\partial_j\xi_i)^\dagger=\bar\partial_j\xi_i,
\)
whereas the Schwarz symmetry of the conjugate system gives
\(
\bar\partial_j\xi_i=\sigma(\bar\partial_i\xi_j)
\). Since $\sigma$ is an involutive $*$--anti-automorphism, we obtain

\[
      (\bar\partial_j\xi_i)^\ast=\sigma\big((\bar\partial_j\xi_i)^{\dagger}\big)
      =\sigma(\bar\partial_j\xi_i)=\bar\partial_i\xi_j.
    \]
This proves the self-adjointness of $\mathscr J\xi$.
\end{proof}

\begin{definition}[{Diez~\cite[Theorem~2.20]{DiezObata}}]\label{CD}
We say that $X$ satisfies the curvature--dimension criterion $CD(c,\infty)$ whenever it has Lipschitz conjugate variables and
\[
\mathscr J\xi
\geq
c.(1\otimes 1)\otimes I_n
\qquad
\text{in }M_n(M\bar\otimes M^{op}).
\]
\end{definition}

\begin{lemma}[{Dabrowski~\cite[Lemma~19]{Dab14}}]\label{dab}
Under the assumptions in Definition~\ref{lip}, we have,
\begin{enumerate}
\item
Set $\Delta_j:=\partial_j^*\bar\partial_j$. For any $x\in M_0$, we have 
$x\in \Dom(\Delta_j)$, 
$\partial_i(x)\in \Dom(\Delta_j\otimes1+1\otimes \Delta_j)$, $x\in \Dom(\Delta^{3/2})$ and:
 \begin{align*} 
&{\partial}_{i}\Delta_{j}(x)= (1 \otimes \Delta_j+\Delta_{j}\otimes 1){\partial}_i(x)
+\partial_{j}(x)\#(\bar{\partial}_{i}\partial_{j}^{*}(1\otimes 1) ).
\end{align*}
\item
 If $x\in \Dom(\bar{\partial})$ (resp. $x\in \Dom(\Delta)$), then so is $\mathrm{id}\otimes\tau({\bar{\partial}_i}(x))$. 
\item
$\Dom(\Delta^{3/2})\subset \Dom(\overline{\Delta \otimes\mathrm{id}+\mathrm{id}\otimes \Delta }\circ\bar{\partial} )$ and, moreover, we have, for any $x\in \Dom(\Delta^{3/2})$, \begin{align*} 
&\bar{\partial}_{i}\Delta(x)= \Delta^{\otimes}\bar{\partial}_i(x)
+\sum_{j=1}^n\bar{\partial}_{j}(x)\#(\bar{\partial}_{i}\partial_{j}^{*}(1\otimes 1) ).
\end{align*}
\end{enumerate}
This operator is the closure of $\Delta\otimes\mathrm{id}+\mathrm{id}\otimes\Delta$ on the core $\Dom(\Delta)\otimes\Dom(\Delta)$, where $\otimes$ denotes the algebraic tensor product. Moreover, $\Dom(\Delta^\otimes)$ is invariant under $\varphi_t\otimes\varphi_t$.
Finally, since
\(
(\mathscr J\xi)_{ji}=\bar\partial_i\xi_j,
\)
the last identity takes the following particularly useful vector form:
\[
\bar\partial\Delta(x)
=
(\Delta^\otimes\otimes I_n+\mathcal R_{\mathscr J\xi})\bar\partial(x).
\]
\end{lemma}
We recall at this point a standard spectral fact which will allow us to pass from one formulation of compactness to another. For a positive self-adjoint operator $A$, the operator $(\alpha+A)^{-1}$ is compact for some (equivalently, every) $\alpha>0$ if and only if $\mathds1_{[0,t]}(A)$ has finite rank for every $t\geq0$; this is also equivalent to compactness of the embedding of $\Dom(A^{1/2})$, endowed with the form graph norm, into the underlying Hilbert space. We will use these equivalences without further justification, and refer to Kato~\cite{Kato}.

\begin{flushleft}
For completeness, we recall that in our setting the resolvent set of $\Delta$ is defined by
\[
\rho(\Delta)
:=
\left\{
z\in\C:
(\Delta-z\,\mathrm{id}):\Dom(\Delta)\to L^2(M)
\text{ is bijective with bounded inverse}
\right\}.
\]
Since $\Delta\geq0$ is self-adjoint, we have
\(
\sigma(\Delta)\subset[0,\infty),
\)
and hence
\(
-1\in\rho(\Delta).
\)
It is thus enough to prove that the resolvent at this single point,
\(
(\Delta+\mathrm{id})^{-1},
\)
is compact. Indeed, compactness of the resolvent at one point of $\rho(\Delta)$ immediately implies compactness at every point of $\rho(\Delta)$, by combining the resolvent identity with the ideal property of compact operators $\mathbb K(L^2(M))$:
\[
(\Delta-z\:\mathrm{id})^{-1}
-
(\Delta-w\:\mathrm{id})^{-1}
=
(z-w)
(\Delta-z\:\mathrm{id})^{-1}
(\Delta-w\:\mathrm{id})^{-1},
\qquad z,w\in\rho(\Delta).
\]
\end{flushleft}

\section{Positive free curvature and compact resolvent}\label{compact}

\subsection{The spectral Bochner estimate}

We now come to the central spectral estimate of this paper. We will prove directly, only from the curvature--dimension criterion and Dabrowski's formula, that the free Laplacian has compact resolvent. For $t\geq0$, we set
\[
N_\Delta(t)
:=
\dim(\mathds1_{[0,t]}(\Delta)L^2(M))
\in\N\cup\{+\infty\}.
\]

\begin{theorem}\label{main}
Let $1\leq n<\infty$, and let $X=(X_1,\ldots,X_n)$ satisfy $CD(c,\infty)$ for some $c>0$. Then $\Delta=\partial^*\bar\partial$ has compact resolvent. More exactly, we have
\begin{equation}\label{gap}
N_\Delta(t)=1,
\qquad
0\leq t<c,
\end{equation}
and
\begin{equation}\label{count}
N_\Delta(t)-1
\leq
nN_\Delta(t-c)^2,
\qquad
t\geq c.
\end{equation}
\end{theorem}

\begin{proof}
For $t\geq0$, we denote by
\[
E_t:=\mathds1_{(0,t]}(\Delta)L^2(M)
\]
the positive spectral subspace of $\Delta$ below $t$. Fix $f\in E_t$, and denote by $\nu_f$ its scalar spectral measure with respect to $\Delta$:
\[
\nu_f(B):=\ip{\mathds1_B(\Delta)f}{f}_2,
\qquad
B\subset[0,\infty)\ \text{Borel}.
\]
By definition, the measure $\nu_f$ is supported in $(0,t]$. Thus, for every $r>0$, we get
\[
\int_{[0,\infty)}\lambda^{2r}\,d\nu_f(\lambda)
=
\int_{(0,t]}\lambda^{2r}\,d\nu_f(\lambda)
\leq
t^{2r}\norm{f}_2^2
<\infty.
\]
The spectral characterization of the domains of the powers of a positive self-adjoint operator therefore shows that $f\in\Dom(\Delta^r)$ for every $r\geq0$. In particular,
\(
f\in\Dom(\Delta^{3/2})
\)
and
\(
\Delta f\in\Dom(\Delta^{1/2})=\Dom(\bar\partial).
\)
\newline
We may now apply Dabrowski's almost-commutation formula from Lemma~\ref{dab}. It gives $\bar\partial_i f\in\Dom(\Delta^\otimes)$ for every $i=1,\ldots,n$, together with the fundamental identity
\[
\bar\partial\Delta f
=
(\Delta^\otimes\otimes I_n+\mathcal R_{\mathscr J\xi})\bar\partial f.
\]
Taking the scalar product of this identity with $\bar\partial f$, we obtain
\begin{align*}
\sum_{i=1}^n\ip{\Delta^\otimes\bar\partial_i f}{\bar\partial_i f}_2
&=
\sum_{i=1}^n\ip{\bar\partial_i\Delta f}{\bar\partial_i f}_2
-
\ip{\mathcal R_{\mathscr J\xi}\bar\partial f}{\bar\partial f}_2.
\end{align*}
But $\Delta f\in\Dom(\bar\partial)$ and $f\in\Dom(\Delta)$. Therefore, the duality for the closed Dirichlet form yields
\[
\sum_{i=1}^n\ip{\bar\partial_i\Delta f}{\bar\partial_i f}_2
=
\ip{\bar\partial\Delta f}{\bar\partial f}_2
=
\ip{\Delta f}{\Delta f}_2
=
\norm{\Delta f}_2^2.
\]
We therefore have the identity
\[
\sum_{i=1}^n\ip{\Delta^\otimes\bar\partial_i f}{\bar\partial_i f}_2
=
\norm{\Delta f}_2^2
-
\ip{\mathcal R_{\mathscr J\xi}\bar\partial f}{\bar\partial f}_2.
\]
We now use in an essential way the fact that $\nu_f$ is supported in $(0,t]$. Indeed,
\begin{align*}
\norm{\Delta f}_2^2
&=
\int_{(0,t]}\lambda^2\,d\nu_f(\lambda)\\
&\leq
t\int_{(0,t]}\lambda\,d\nu_f(\lambda)\\
&=
t\ip{\Delta f}{f}_2
=
t\norm{\bar\partial f}_2^2.
\end{align*}
On the other hand, the condition $CD(c,\infty)$ combined with Proposition~\ref{rightprop} gives
\[
\ip{\mathcal R_{\mathscr J\xi}\bar\partial f}{\bar\partial f}_2
\geq
c\norm{\bar\partial f}_2^2.
\]
Combining all the preceding identities and estimates, we arrive at the following spectral-lowering estimate
\begin{align}
\ip{(\Delta^\otimes\otimes I_n)\bar\partial f}{\bar\partial f}_2
&=
\sum_{i=1}^n\ip{\Delta^\otimes\bar\partial_i f}{\bar\partial_i f}_2\nonumber\\
&\leq
(t-c)\norm{\bar\partial f}_2^2.
\label{lower}
\end{align}
We first consider the interval $0\leq t<c$. The left-hand side of~\eqref{lower} is non-negative since $\Delta^\otimes\geq0$, and its right-hand side is strictly negative unless $\bar\partial f=0$. We thus necessarily have $\bar\partial f=0$, and~\eqref{kernel} shows that $f$ is scalar. But also
\[
E_t\perp\Ker(\Delta)=\mathbb C\,1,
\]
because $E_t$ and $\Ker(\Delta)$ are the ranges of the spectral projections corresponding to the disjoint Borel sets $(0,t]$ and $\{0\}$. We obtain $f=0$, hence $E_t=\{0\}$, and consequently
\[
N_\Delta(t)=\dim\Ker(\Delta)=1,
\qquad
0\leq t<c.
\]
This proves~\eqref{gap}, and in particular already proves that the spectral gap is forced by the convexity criterion.
\newline
We turn to the more interesting case $t\geq c$. Set $s=t-c$, and assume for the moment that
\[
N_\Delta(s)<\infty.
\]
We consider the following natural lowering map:
\begin{align}
E_t
&\to
\left(
\mathds1_{[0,s]}(\Delta^\otimes)
\big(L^2(M)\bar\otimes L^2(M^{op})\big)
\right)^{\oplus n},
\nonumber\\
f
&\mapsto
\big(
\mathds1_{[0,s]}(\Delta^\otimes)\bar\partial_1f,
\ldots,
\mathds1_{[0,s]}(\Delta^\otimes)\bar\partial_nf
\big).
\label{lowmap}
\end{align}
The main point is that this map is injective. To see it, let $f\in E_t$ and suppose that
\[
\mathds1_{[0,s]}(\Delta^\otimes)\bar\partial_i f=0,
\qquad
i=1,\ldots,n.
\]
For every $i=1,\ldots,n$, denote by
\[
\nu_{i,f}(B)
:=
\ip{\mathds1_B(\Delta^\otimes)\bar\partial_i f}{\bar\partial_i f}_2,
\qquad
B\subset[0,\infty)\ \text{Borel},
\]
the scalar spectral measure of $\bar\partial_i f$ with respect to $\Delta^\otimes$. Our assumption gives $\nu_{i,f}([0,s])=0$. Moreover, since $\bar\partial_i f\in\Dom(\Delta^\otimes)$, the first moment of $\nu_{i,f}$ is finite, and we have
\begin{align*}
\sum_{i=1}^n
\left(
\ip{\Delta^\otimes\bar\partial_i f}{\bar\partial_i f}_2
-
s\norm{\bar\partial_i f}_2^2
\right)
=
\sum_{i=1}^n
\int_{(s,\infty)}(\lambda-s)\,d\nu_{i,f}(\lambda).
\end{align*}
Assume by contradiction that $\bar\partial f\neq0$. Then $\bar\partial_i f\neq0$ for at least one $i$, and consequently
\[
\nu_{i,f}((s,\infty))
=
\norm{\bar\partial_i f}_2^2
>0.
\]
Since
\(
(s,\infty)=\bigcup_{m\geq1}[s+m^{-1},\infty),
\)
we can find $m\geq1$ such that $\nu_{i,f}([s+m^{-1},\infty))>0$. It follows at once that
\[
\int_{(s,\infty)}(\lambda-s)\,d\nu_{i,f}(\lambda)
\geq
\frac1m\nu_{i,f}([s+m^{-1},\infty))
>0.
\]
Thus, the right-hand side above is strictly positive, and we obtain
\[
\sum_{i=1}^n\ip{\Delta^\otimes\bar\partial_i f}{\bar\partial_i f}_2
>
s\norm{\bar\partial f}_2^2.
\]
But $s=t-c$, and this is in direct contradiction with~\eqref{lower}. Hence $\bar\partial f=0$. Using again~\eqref{kernel}, together with $f\perp1$, we get $f=0$. This proves that the map~\eqref{lowmap} is injective.
\newline
The target space in~\eqref{lowmap} is the $n$--fold direct sum of the range of $\mathds1_{[0,s]}(\Delta^\otimes)$. By Proposition~\ref{tensor} and our assumption $N_\Delta(s)<\infty$, this space is finite-dimensional. The injectivity proved above now gives
\begin{align*}
\dim E_t
&\leq
n\,\rank\mathds1_{[0,s]}(\Delta^\otimes)\\
&\leq
nN_\Delta(s)^2
=
nN_\Delta(t-c)^2.
\end{align*}
Moreover, the spectral calculus gives the orthogonal decomposition
\[
\mathds1_{[0,t]}(\Delta)L^2(M)
=
\Ker(\Delta)\oplus E_t,
\]
and we already know that $\dim\Ker(\Delta)=1$. We have thus established
\begin{equation}\label{count2}
N_\Delta(t)-1
\leq
nN_\Delta(t-c)^2
\end{equation}
under the only condition that $t\geq c$ and $N_\Delta(t-c)<\infty$.
\newline
It only remains to prove that $N_\Delta(t)<\infty$ for every finite $t$, and thus to remove this (apparently circular) finiteness assumption. We do this by a very simple induction. For $k\geq0$, consider the assertion
\[
\mathcal P(k):
\qquad
N_\Delta(t)<\infty
\quad\text{for every }0\leq t\leq kc.
\]
The case $k=0$ follows from $N_\Delta(0)=\dim\Ker(\Delta)=1$. Suppose that $\mathcal P(k)$ holds for some $k\geq0$, and take $0\leq t\leq(k+1)c$. If $0\leq t<c$, then~\eqref{gap} gives $N_\Delta(t)=1$. If instead $c\leq t\leq(k+1)c$, then $0\leq t-c\leq kc$, and the induction hypothesis gives $N_\Delta(t-c)<\infty$. We can therefore apply the conditional estimate~\eqref{count2}, which yields
\[
N_\Delta(t)
\leq
1+nN_\Delta(t-c)^2
<\infty.
\]
This proves $\mathcal P(k+1)$. By induction, $N_\Delta(t)<\infty$ for every finite $t$. In particular, the finiteness assumption in~\eqref{count2} is now automatically satisfied for every $t\geq c$, and~\eqref{count} follows.
\newline
We finally deduce compactness of the resolvent. For $R>0$, the spectral calculus gives
\[
(1+\Delta)^{-1}
=
\mathds1_{[0,R]}(\Delta)(1+\Delta)^{-1}
+
\mathds1_{(R,\infty)}(\Delta)(1+\Delta)^{-1}.
\]
The first operator has finite rank, since $N_\Delta(R)<\infty$, while
\[
\norm{\mathds1_{(R,\infty)}(\Delta)(1+\Delta)^{-1}}
\leq
\frac1{1+R}.
\]
We have proved that $(1+\Delta)^{-1}$ is a norm limit of finite-rank operators. It is therefore compact, and the proof is complete.
\end{proof}
\begin{flushleft}
The argument above gives at first sight a quadratic recursion for the
counting function. This is of course sufficient to prove compactness of
the resolvent, but its iteration only gives a
\emph{superexponential} upper bound. This loss comes precisely from the
elementary inclusion
\[
\left\{
(\lambda,\mu)\in[0,\infty)^2:\lambda+\mu\leq s
\right\}
\subset
[0,s]^2.
\]
Indeed, the spectral subspace of $\Delta^\otimes$ below $s$ is
determined by the triangular region $0\leq \lambda+\mu\leq s$, whereas in the
preceding proof we proceed to bound its dimension by counting all pairs in the
larger square $[0,s]^2$. Keeping track more precisely of the additive constraint
gives instead a convolution inequality for the spectral counting
measure. The discrete comparison equation associated with this
inequality is given by the Catalan recursion. This leads to an exponential
bound, together with a trace estimate for the heat
semigroup. We are very grateful to David Jekel for pointing out this observation to us. 
\end{flushleft}
\begin{flushleft}
    In the sequel, for \(k\geq0\), we denote by
\(
C_k:=\frac{1}{k+1}\binom{2k}{k}
\)
the \(k\)-th Catalan number. It is given by \(C_0=1\) and the recursion
\begin{equation}\label{catalan}
C_{k+1}
=
\sum_{j=0}^{k}C_jC_{k-j},
\qquad
k\geq0.
\end{equation}
\end{flushleft}

\begin{flushleft} Before going further, let us record an elementary lemma showing that cumulative domination is preserved under classical convolution. \end{flushleft} \begin{lemma}\label{lmacumul}
Let \(\alpha\) and \(\beta\) be locally finite positive Borel measures on \([0,\infty)\), and let \(s\geq0\). Assume that \[ \alpha([0,u]) \leq \beta([0,u]), \qquad \text{for all}\quad 0\leq u\leq s. \]
Then \begin{equation}\label{convdom} (\alpha*\alpha)([0,s]) \leq (\beta*\beta)([0,s]). \end{equation} \end{lemma} 
\begin{proof}
Indeed, we have
\begin{align*}
(\alpha*\alpha)([0,s])
&=
\int_{[0,s]}
\alpha([0,s-u])\,d\alpha(u)\\
&\leq
\int_{[0,s]}
\beta([0,s-u])\,d\alpha(u)\\
&=
(\alpha*\beta)([0,s])\\
&=
\int_{[0,s]}
\alpha([0,s-u])\,d\beta(u)\\
&\leq
\int_{[0,s]}
\beta([0,s-u])\,d\beta(u)\\
&=
(\beta*\beta)([0,s]).
\end{align*}
Here both inequalities follow directly from the cumulative domination
assumption, since \(s-u\in[0,s]\) for every \(u\in[0,s]\). The middle
identity follows from the
commutativity of classical convolution and Tonelli's theorem.
\end{proof}
\begin{corollary}[Exponential upper bound for spectral growth]
\label{expon}
Under the assumptions of Theorem~\ref{main}, we have more precisely
\[
N_\Delta(t)
\leq
\sum_{k=0}^{\lfloor t/c\rfloor}C_kn^k
\leq
\frac{(4n)^{\lfloor t/c\rfloor+1}-1}{4n-1},
\qquad
t\geq0.
\]
Consequently,
\(
\limsup_{t\to\infty}
\frac{\log N_\Delta(t)}{t}
\leq
\frac{\log(4n)}{c}.
\)
Moreover, \(e^{-t\Delta}\) is of trace class for every
\(t\geq \frac{\log(4n)}{c}\), and
\[
\operatorname{Tr}_{L^2(M)}(e^{-t\Delta})
\leq
\frac{2}{1+\sqrt{1-4ne^{-t c}}}.
\]
\end{corollary}
\begin{proof}
By Theorem~\ref{main}, \(\Delta\) has compact resolvent. We can
therefore consider its spectral counting measure (namely, the trace of its spectral measure),
\(
\mu_\Delta
:=
\sum_{\lambda\in\sigma(\Delta)}
\dim\Ker(\Delta-\lambda\:\mathrm{id})\,\delta_\lambda,
\)
so that
\(
N_\Delta(t)=\mu_\Delta([0,t]).
\)
The only change in the preceding proof occurs when estimating the
dimension of the target of~\eqref{lowmap}. Indeed, the joint
spectral calculus recalled in Proposition~\ref{tensor} gives the precise identity
\[
\rank\mathds1_{[0,s]}(\Delta^\otimes)
=\sum_{\substack{\lambda,\mu\in\sigma(\Delta)\\
                \lambda+\mu\leq s}}
\dim\Ker(\Delta-\lambda\:\mathrm{id})\,
\dim\Ker(\Delta-\mu\:\mathrm{id})=
(\mu_\Delta*\mu_\Delta)([0,s]).
\]
The injectivity of~\eqref{lowmap} therefore yields
\begin{equation}\label{conv}
N_\Delta(t)-1
\leq
n(\mu_\Delta*\mu_\Delta)([0,t-c]),
\qquad
t\geq c.
\end{equation}
\newline
The form of inequality~\eqref{conv} suggests looking for a
discrete measure which satisfies the same relation with equality.
Since each successive application of the descent procedure decreases the
spectral level by \(c\), it is natural to look for such a measure supported
on the non-negative lattice \(c\mathbb Z_{\geq 0}=\{kc:k\geq0\}\).
\newline
Let us write
\(
\nu
:=
\sum_{k=0}^{\infty}a_k\,\delta_{kc},
\)
where \(a_k\geq0\). We impose the fixed-point convolution equation
\(
\nu
=
\delta_0+n\,\delta_c*(\nu*\nu).
\)
Comparing first the masses at \(0\), and then at \((k+1)c\), we obtain
\[
a_0=1,
\qquad
a_{k+1}
=
n\sum_{j=0}^{k}a_ja_{k-j},
\qquad
k\geq0.
\]
Since a locally finite measure supported on \(c\mathbb Z_{\geq0}\) is uniquely determined by the sequence of its point masses, this recursion determines \(\nu\) uniquely. The Catalan recursion therefore gives
\(
a_k=C_kn^k
\)
for every \(k\geq0\).
\bigbreak
Consequently, the unique locally finite positive measure supported on \(c\mathbb Z_{\geq0}\) which satisfies the preceding fixed-point relation is
\(
\nu
:=
\sum_{k=0}^{\infty}C_kn^k\,\delta_{kc}.
\)
It is locally finite since every bounded interval contains only finitely many points of \(c\mathbb Z_{\geq0}\). The Catalan recursion
\(
C_{k+1}=\sum_{j=0}^{k}C_jC_{k-j}
\)
is equivalent to the convolution identity
\(
\nu
=
\delta_0+n\,\delta_c*(\nu*\nu).
\)
Thus, setting \(G(t):=\nu([0,t])\), we have
\(
G(t)-1
=
n(\nu*\nu)([0,t-c]),
\) for
\(t\geq c\).
\newline
We claim that \(N_\Delta(t)\leq G(t)\) for every \(t\geq0\).
By~\eqref{gap}, this holds for \(0\leq t<c\), while
\eqref{conv} at \(t=c\) gives
\[
N_\Delta(c)\leq1+n=G(c).
\]
Suppose that the claim holds on \([0,kc]\) for some \(k\geq1\), and
take \(kc\leq t\leq(k+1)c\). Set \(s=t-c\), so that \(s\leq kc\).
For every \(0\leq r\leq s\), the induction hypothesis gives
\[
\mu_\Delta([0,r])
=
N_\Delta(r)
\leq
G(r)
=
\nu([0,r]).
\]
Now using Lemma~\ref{lmacumul}, applied with
\(\alpha=\mu_\Delta\) and \(\beta=\nu\), gives
\[
(\mu_\Delta*\mu_\Delta)([0,s])
\leq
(\nu*\nu)([0,s]).
\]
Combining this with~\eqref{conv}, we obtain
\[
N_\Delta(t)-1
\leq
n(\nu*\nu)([0,t-c])
=
G(t)-1.
\]
The assertion thus follows by induction.
\newline
Using the elementary bound
\(C_k\leq4^k\), we conclude that
\[
N_\Delta(t)
\leq
\sum_{k=0}^{\lfloor t/c\rfloor}C_kn^k
\leq
\sum_{k=0}^{\lfloor t/c\rfloor}(4n)^k
=
\frac{(4n)^{\lfloor t/c\rfloor+1}-1}{4n-1},
\]
which proves the first estimate.
\newline
It remains to establish the trace-class estimate of the heat semigroup. Fix \(t>0\) (at this stage, all
the quantities below initially are understood as elements of \([0,\infty]\)). For every
\(\lambda\geq0\), one has
\[
e^{-t\lambda}
=
t\int_0^\infty
e^{-tr}\mathds1_{[0,r]}(\lambda)\,dr.
\]
Since \(N_\Delta(r)\leq G(r)\) for every \(r\geq0\), Tonelli's theorem
therefore gives
\begin{align*}
\operatorname{Tr}_{L^2(M)}(e^{-t\Delta})
&=
\int_{[0,\infty)}e^{-t\lambda}\,d\mu_\Delta(\lambda)\\
&=
t\int_0^\infty e^{-tr}
\left(
\int_{[0,\infty)}
\mathds1_{[0,r]}(\lambda)\,d\mu_\Delta(\lambda)
\right)dr\\
&=
t\int_0^\infty e^{-tr}\mu_\Delta([0,r])\,dr\\
&\leq
t\int_0^\infty e^{-tr}\nu([0,r])\,dr\\
&=
\int_{[0,\infty)}e^{-t\lambda}\,d\nu(\lambda)\\
&=
\sum_{k=0}^{\infty}C_k(ne^{-tc})^k.
\end{align*}
If \(C\) denotes the generating function of the Catalan numbers, then, for every \(0<x\leq1/4\),
\newline
\(
C(x):=\sum_{k=0}^{\infty}C_kx^k=\frac{1-\sqrt{1-4x}}{2x}.
\)
Taking \(x=ne^{-t c}\), we finally obtain
\[
\operatorname{Tr}_{L^2(M)}(e^{-t\Delta})
\leq
\frac{1-\sqrt{1-4ne^{-tc}}}{2ne^{-tc}}=\frac{2}{1+\sqrt{1-4ne^{-t c}}}
\]
when
\(
t\geq \frac{\log(4n)}{c}.
\)
\end{proof}

\begin{remark}[The classical uniformly log-concave case]\label{classical}
It is interesting to observe that exactly the same spectral-descent argument applies mutatis mutandis in the classical uniformly log-concave setting, by using the usual \emph{integrated Bochner formula} (see again Bakry--Gentil--Ledoux~\cite{Ledoux} for a complete overview of the classical case). Compactness of the resolvent, and even much stronger spectral information, such as decay estimates for the associated Schr\"odinger eigenfunctions~\cite{Davies,ReedSimonIV}, are of course classical in this setting. We have not, however, been able to find this particular short proof in the literature.
\end{remark}

\begin{remark}\label{finite}
Let us remark here that the finiteness of $n$ is essential. Indeed, an infinite semicircular family satisfies the same curvature identity with $c=1$, but its first Wigner chaos is infinite-dimensional (see Biane--Speicher~\cite{BS} for the definition of the multiple Wigner integrals). Therefore, the number operator has eigenvalue $1$ with infinite multiplicity. Positive curvature gives the fixed loss of energy at each step, on the successive intervals $[kc,(k+1)c)$.
\end{remark}

\begin{flushleft}
We conclude this section by recalling the direct application of the preceding theorem to bounded polynomial free Gibbs laws, whenever such realizations exist (we do not enter here into the much more delicate discussion of existence and uniqueness; we refer to Jekel, Dabrowski, Guionnet, and Shlyakhtenko).
\end{flushleft}

\begin{example}[Free Gibbs states]\label{gibbs}
Let
\(
\mathscr P:=\C\langle t_1,\ldots,t_n\rangle
\)
denote the free algebra in $n$ self-adjoint non-commuting indeterminates, and let $V=V^*\in\mathscr P$. We say that a bounded self-adjoint tuple $X=(X_1,\ldots,X_n)$ is a free Gibbs realization for $V$ when its tracial law $\tau_V$ satisfies the Schwinger--Dyson equations
\begin{equation}\label{schwin}
\tau_V(\mathscr D_iV(X)P(X))
=
(\tau_V\otimes\tau_V)(\partial_iP(X)),
\qquad
P\in\mathscr P,\quad i=1,\ldots,n.
\end{equation}
These equations immediately identify the conjugate variables of $X$ as
\[
\xi_i=\mathscr D_iV(X).
\]
Here $\mathscr DV=(\mathscr D_1V,\ldots,\mathscr D_nV)$ is the cyclic gradient, and $\mathscr D_iV=m\circ\sigma(\partial_iV)$ is its $i$--th component.
In particular, Dabrowski~\cite[Lemma~37]{Dab14} shows that $X$ is algebraically free, so evaluation identifies $\mathscr P$ with
\(
M_0=\C\langle X_1,\ldots,X_n\rangle.
\)
Moreover, since $V$ is polynomial and $X$ is bounded, we have
\[
\bar\partial_j\xi_i
=
\partial_j\mathscr D_iV(X)
\in
M_0\otimes M_0^{op}
\subset
M\bar\otimes M^{op}.
\]
Thus, in this important situation, the conjugate variables are automatically Lipschitz, and
\[
\mathscr J\xi
=
\mathscr J\mathscr DV(X)
\in
M_n(M_0\otimes M_0^{op})
\subset
M_n(M\bar\otimes M^{op}).
\]
\end{example}

\begin{corollary}[Uniformly convex Gibbs laws]\label{gibbscor}
Let $X$ be a bounded free Gibbs realization for a polynomial potential $V=V^*$. Suppose that, for some $c>0$,
\[
\mathscr J\mathscr DV(X)
\geq
c.(1\otimes 1)\otimes I_n
\qquad
\text{in }M_n(M\bar\otimes M^{op}).
\]
Then the ordinary free-difference-quotient Laplacian
\(
\Delta=\partial^*\bar\partial
\)
has compact resolvent, and all the conclusions of Theorem~\ref{main} hold.
\end{corollary}
\section{Operator-algebraic consequences}\label{operator}
\begin{flushleft}
We now arrive at the operator-algebraic part of the paper. We first
combine our curvature criterion with Dabrowski's stationary solutions
of a free stochastic partial differential equation~\cite{Dab14}. The $CD(c,\infty)$ criterion
makes the stochastic coefficient square-integrable on the whole
half-line and allows us to construct a normal unital trace-preserving
embedding into $L(\mathbb F_\infty)$.

We then return to the consequences of compactness of the free Laplacian, in a very general context without even assuming positive curvature.
We obtain by direct arguments the Haagerup approximation property,
solidity, proper proximality and, under weak exactness, bi-exactness,
since the arguments apply without the curvature assumption. For the tuples of
Theorem~\ref{main}, the embedding also gives the complete metric
approximation property and weak exactness. The additional assumptions
in the strong-solidity and bi-exactness  statements are therefore
always ensured in our $CD(c,\infty), \: c>0$ setting. Throughout this section,
von Neumann subalgebras are always understood to contain the unit of the
ambient algebra.
\end{flushleft}
\subsection{Positive curvature and the stationary stochastic equation}
\label{subsec:curvature-embedding}

We apply the stationary stochastic construction of Dabrowski~\cite{Dab14}.
Combining his energy identity with the resolvent form of almost commutation,
used also in our previous work~\cite{DiezObata}, we show that positive
curvature makes the stochastic coefficient square-integrable on the whole
half-line. We can then pass to the limit directly in the mild equation and deduce the following embedding into free group factors which is also of independent interest. In the sequel $R$ denotes the unique hyperfinite $\mathrm{II}_1$ factor and $R^\omega$ its ultrapower with respect to a free ultrafilter $\omega$.

\begin{remark}
A full presentation of Dabrowski's stochastic framework would
require a substantial amount of additional notation and technical
preliminaries. In the sequel, we therefore recall only the notation
and properties needed for our argument, following~\cite{Dab14},
and refer the reader to that paper for the general construction,
the precise assumptions and the proofs of the results invoked.
All references to sections, equations and numbered statements
in~\cite{Dab14} are to the published version.
\end{remark}
\begin{theorem}[Embedding into free group factors]
\label{thm:curvature-embedding}
Let $X=(X_1,\ldots,X_n)$, $1\leq n<\infty$, be a bounded self-adjoint
tuple with Lipschitz conjugate variables satisfying $CD(c,\infty)$
for some $c>0$. Then there exists a normal unital trace-preserving embedding
\[
 \Phi_\infty:M\hookrightarrow L(\mathbb F_\infty).
\]
Consequently, $M$ also admits a normal unital trace-preserving embedding
into $L(\mathbb F_2)$.
\end{theorem}

\begin{proof}
Let $(S_t^{(i)})_{t\geq0}$, $1\leq i\leq n$, be free (standard) Brownian motions,
free from $M$. Set
\[
 N_t=W^*(S_s^{(i)}:0\leq s\leq t,\ 1\leq i\leq n),
 \qquad N=W^*(N_t:\ t\geq0),
\]
and $\widetilde M=M*N$, $\widetilde M_t=M*N_t$.
The trace on $\widetilde M$ is again denoted without confusion by $\tau$, and we have
$N\cong L(\mathbb F_\infty)$.
Following~\cite[Section~2.1]{Dab14}, extend the polynomial free difference
quotient by zero on the free Brownian variables:
\[
 \Dom(\delta)=M_0*
 \mathbb C\langle S_t^{(i)}:t\geq0,\ 1\leq i\leq n\rangle,
 \qquad \delta_i|_{M_0}=\partial_i,
 \quad \delta_i(S_t^{(j)})=0,
\]
where the free product is algebraic. This derivation is closable.
Denote its joint closure by $\bar\delta$ and put, for $\alpha>0$,
\[
 \widetilde\Delta=\delta^*\bar\delta,\qquad
 \widetilde\varphi_t=e^{-t\widetilde\Delta},\qquad
 \widetilde\rho_\alpha=\alpha(\alpha+\widetilde\Delta)^{-1}.
\]
The subspace $L^2(M)$ reduces $\widetilde\Delta$, whose restriction there
is $\Delta$. We also write
\[
 \widetilde\Delta^\otimes
 =\overline{\widetilde\Delta\otimes\mathrm{id}
                 +\mathrm{id}\otimes\widetilde\Delta},\qquad
 \widetilde\varphi_t^\otimes=e^{-t\widetilde\Delta^\otimes},\qquad
 \widetilde\rho_\alpha^\otimes
       =\alpha(\alpha+\widetilde\Delta^\otimes)^{-1}.
\]
Tensor operators act componentwise on the $n$-fold direct sum.
As usual, we use the notation
$(a\otimes b^{op})\#dS_t^{(i)}=a\,dS_t^{(i)}b$ for stochastic integration.
The filtration properties of these operators are those established
in~\cite[Section~2.1, Lemma~16]{Dab14}.

The Lipschitz conjugate variable hypothesis satisfies~\cite[Assumption~1]{Dab14}, as explained in the proof of~\cite[Corollary~25]{Dab14} there.
The results of~\cite[Theorems~17 and~11(ii)]{Dab14}, together with~\cite[Propositions~23--24]{Dab14}, therefore
provide normal unital trace-preserving embeddings
$\Phi_t:M\to\widetilde M_t$ on this free product.
For $x\in M_0$, the process $Y_t=\Phi_t(x)$ is the mild solution
\begin{equation}\label{eq:embedding-mild}
 Y_t=\widetilde\varphi_{t/2}(x)+\gamma(Y)_t,\qquad
 \gamma(Y)_t=\sum_{i=1}^n\int_0^t
   \widetilde\varphi_{(t-s)/2}^\otimes
                  (\bar\delta_iY_s)\#dS_s^{(i)}.
\end{equation}
We use its continuous $L^2$ representative. In particular,
$Y\in\mathcal B_{2,\bar\delta}^{a}([0,T])$ for every $T<\infty$,
where the adapted graph spaces are those of~\cite[Section~1.2]{Dab14} (the factors $1/2$ in~\eqref{eq:embedding-mild} come from Dabrowski's
normalization~\cite{Dab14}, which uses the generator $-\widetilde\Delta/2$).
\newline
For $\epsilon\in(0,1)$, let $Y^\epsilon$ be the regularized solution
from~\cite[Theorem~11(ii)]{Dab14}:
\begin{equation}\label{eq:embedding-epsilon}
 Y_t^\epsilon=\widetilde\varphi_{t/2}(x)
                  +(1-\epsilon)\gamma(Y^\epsilon)_t.
\end{equation}
Since the two derivations in that theorem coincide here, $Y^\epsilon$
is a strong solution and belongs to
$\mathcal B_{2,\widetilde\Delta}^{a}([0,T])$ for every $0<T<\infty$.
The spaces $\mathcal B_{2,\bar\delta}^{a}([0,T])$ and
$\mathcal B_{2,\widetilde\Delta}^{a}([0,T])$ are the completions
of the corresponding simple adapted processes for the respective
norms
\[
 \|X\|_{\mathcal B_{2,\bar\delta}^{a}([0,T])}
 :=
 \left(
 \int_0^T
 \left[
 \|\bar\delta X_s\|_2^2
 +\|X_s\|_2^2
 \right]ds
 \right)^{1/2}
\]
and
\[
 \|X\|_{\mathcal B_{2,\widetilde\Delta}^{a}([0,T])}
 :=
 \left(
 \int_0^T\sum_{b=0}^{2}
 \|\widetilde\Delta^{b/2}X_s\|_2^2\,ds
 \right)^{1/2}.
\]
These correspond to the choices $\beta=0$ and $\beta=1$,
respectively, in the definitions of~\cite[Section~1.2]{Dab14}.
In particular, $Y_s^\epsilon\in\Dom(\widetilde\Delta)$
for almost every $s\in[0,T]$, and
\[
 \int_0^T\|\widetilde\Delta Y_s^\epsilon\|_2^2\,ds<\infty.
\]
Moreover, as $\epsilon\downarrow0$, $Y^\epsilon$ converges to $Y$
weakly in $\mathcal B_{2,\bar\delta}^{a}([0,T])$
and uniformly on $[0,T]$ in $L^2(\widetilde M)$.
\newline
We first record the coercivity given by positive curvature, noting that the conjugate variables of the zero extension remain $\xi_i$.
By using again Lemma~\ref{dab}, in its zero-extension form from~\cite[Lemma~19(iii)]{Dab14},
\[
 \bar\delta\widetilde\Delta u
 =(\widetilde\Delta^\otimes\otimes I_n
                 +\mathcal R_{\mathscr J\xi})\bar\delta u,
 \qquad u\in\Dom(\widetilde\Delta^{3/2}),
\]
including the corresponding domain assertion.
The spatial inclusion
$M\bar\otimes M^{op}\subset\widetilde M\bar\otimes\widetilde M^{op}$
preserves positivity. By the representation properties recalled
above in Proposition~\ref{rightprop}, we also have
$\mathcal R_{\mathscr J\xi}\geq c\,\mathrm{id}$ on the enlarged
gradient space. The free Bochner formula (taking the inner product with $\bar\delta u$)
therefore gives
\[
 \|\widetilde\Delta u\|_2^2
 =\|(\widetilde\Delta^\otimes\otimes I_n)^{1/2}\bar\delta u\|_2^2
   +\langle\mathcal R_{\mathscr J\xi}\bar\delta u,\bar\delta u\rangle_2
 \geq
 \|(\widetilde\Delta^\otimes\otimes I_n)^{1/2}\bar\delta u\|_2^2
   +c\|\bar\delta u\|_2^2.
\]
For $u\in\Dom(\widetilde\Delta)$, apply this inequality to
$\widetilde\rho_\alpha u\in\Dom(\widetilde\Delta^{3/2})$.
We have convergence to $u$ in the graph norm of $\widetilde\Delta$,
and the free difference quotients converge to $\bar\delta u$ in $L^2$.
Applying the same inequality to
$\widetilde\rho_\alpha u-\widetilde\rho_\beta u$, we see that
$(\widetilde\Delta^\otimes\otimes I_n)^{1/2}\bar\delta\widetilde\rho_\alpha u$
is Cauchy in $L^2$. By closedness, we therefore obtain
$\bar\delta u\in\Dom((\widetilde\Delta^\otimes\otimes I_n)^{1/2})$
and can pass to the limit in the inequality
(see also~\cite[proof of Theorem~17]{Dab14}, explained in greater detail
in~\cite[Lemma~3.6]{DiezObata}).
Thus, writing
\[
 \widetilde{\mathcal E}_2(u)
 :=\|(\widetilde\Delta^\otimes\otimes I_n)^{1/2}\bar\delta u\|_2^2,
\]
we have
\begin{equation}\label{eq:embedding-coercivity}
 \widetilde{\mathcal E}_2(u)+c\|\bar\delta u\|_2^2
 \leq\|\widetilde\Delta u\|_2^2,
 \qquad u\in\Dom(\widetilde\Delta).
\end{equation}
We fix $\epsilon\in(0,1)$ and $T>0$.
We apply~\cite[Proposition~5, equation~(6)]{Dab14} with
$B=\widetilde\Delta^{1/2}$ and $\alpha'=\alpha$, as in ~\cite[Step~1
of the proof of Theorem~11 and
Lemma~3(ii)--(iii)]{Dab14}. Together with the free It\^o isometry of Biane--Speicher~\cite{BS}, this gives the following identity
(the real-part of the inner product in~\cite[Proposition~5, equation~(6)]{Dab14}
can be omitted here, since the corresponding inner product equals a squared norm by spectral calculus):
\[
 \|\widetilde\Delta^{1/2}\widetilde\rho_\alpha
                         \gamma(Y^\epsilon)_T\|_2^2
   +\int_0^T\|\widetilde\Delta\widetilde\rho_\alpha
                         \gamma(Y^\epsilon)_s\|_2^2\,ds=
   \int_0^T
   \|(\widetilde\Delta^\otimes\otimes I_n)^{1/2}
       \widetilde\rho_\alpha^\otimes\bar\delta Y_s^\epsilon\|_2^2\,ds.
\]
By~\eqref{eq:embedding-coercivity}, the right-hand side is bounded by
\[
 \int_0^T\widetilde{\mathcal E}_2(Y_s^\epsilon)\,ds
 \leq\int_0^T\|\widetilde\Delta Y_s^\epsilon\|_2^2\,ds<\infty.
\]
We can let $\alpha\to\infty$ by using monotone convergence.
The (deterministic) semigroup satisfies
\[
 \|\widetilde\Delta^{1/2}\widetilde\varphi_{T/2}(x)\|_2^2
 +\int_0^T\|\widetilde\Delta\widetilde\varphi_{s/2}(x)\|_2^2\,ds
 =\|\bar\partial x\|_2^2.
\]
Moreover, $\gamma(Y^\epsilon)_t\perp L^2(M)$, since its conditional
expectation onto $M$ is zero, while the deterministic
term belongs to $L^2(M)$. Since this subspace reduces
$\widetilde\Delta$, these terms remain orthogonal after applying the
spectral powers above. Adding the two energy identities and
using~\eqref{eq:embedding-epsilon}, we obtain
\begin{equation}\label{eq:embedding-regularized-energy}
 \begin{split}
 \|\bar\delta Y_T^\epsilon\|_2^2
 +\int_0^T\|\widetilde\Delta Y_s^\epsilon\|_2^2\,ds
 =\|\bar\partial x\|_2^2
 +(1-\epsilon)^2\int_0^T
                 \widetilde{\mathcal E}_2(Y_s^\epsilon)\,ds.
 \end{split}
\end{equation}
Since $1-(1-\epsilon)^2\geq0$, the coercivity
estimate~\eqref{eq:embedding-coercivity} yields
\[
 \|\bar\delta Y_T^\epsilon\|_2^2
 +c\int_0^T\|\bar\delta Y_s^\epsilon\|_2^2\,ds
 \leq\|\bar\partial x\|_2^2.
\]
Omitting the first term and using weak lower semicontinuity in
$\mathcal B_{2,\bar\delta}^{a}([0,T])$, we can reach the limit
$\epsilon\to0$. Letting then $T\to\infty$, we obtain
\begin{equation}\label{eq:embedding-integrated-gradient}
 \int_0^\infty\|\bar\delta\Phi_s(x)\|_2^2\,ds
 \leq\frac1c\|\bar\partial x\|_2^2.
\end{equation}
We next identify the strong limit of $\widetilde\varphi_r$. Since $\delta$ vanishes on the polynomials in the free Brownian
variables, this is clear by density that $L^2(N)\subset \Ker(\widetilde\Delta)$. We now prove the reverse inclusion using the word-space
decomposition introduced before assumption $\Gamma_{0u}(\omega)$
in~\cite[p.~1409]{Dab14} and the tensor representation of the zero
extension in~\cite[Section~2.1.3]{Dab14}.
We choose an orthonormal basis $(e_j)_{j\geq1}$ of
$L^2(N)\ominus\mathbb C.1$ obtained by applying Gram--Schmidt to a countable $L^2$-dense family of trace-zero non-commutative polynomials in the free Brownian variables at non-negative rational times.
\newline
For $k\geq1$, set
\[
 \mathcal V_k
 =L^2(M)\bar\otimes
   (L^2(M)\ominus\mathbb C.1)^{\bar\otimes(k-1)}
   \bar\otimes L^2(M),
\]
with the middle factor omitted when $k=1$.
On this space, $\Delta^{\otimes(k+1)}$ denotes the positive
self-adjoint closure of
\[
 \Delta^{\otimes(k+1)}(a_0\otimes\cdots\otimes a_k)
 =\sum_{j=0}^k
   a_0\otimes\cdots\otimes\Delta(a_j)\otimes\cdots\otimes a_k
\]
on $\Dom(\Delta)^{\otimes (k+1)}\cap \mathcal V_k$ (the algebraic tensor product).
For each $i=(i_1,\ldots,i_k)\in\mathbb N^k$, write
$e_i=e_{i_1}\otimes\cdots\otimes e_{i_k}$.
We also use the multiplication map
\(
 (a_0\otimes\cdots\otimes a_k)\#e_i
 =a_0e_{i_1}a_1\cdots e_{i_k}a_k,
\)
initially defined for $a_0,\ldots,a_k\in M$ with
$\tau(a_j)=0$ for $1\leq j\leq k-1$.
By freeness, this map extends to an isometry on $\mathcal V_k$, we also denote its (closed) range by setting
\(
 \mathcal V_k\#e_i
 :=\{U\#e_i:U\in\mathcal V_k\}
\). By the standard Hilbert space decomposition of a tracial free product (see~\cite[Section~5.3.2 and Remark~5.3.5]{AnanPopa})
\[
 L^2(\widetilde M)
 =L^2(M)\oplus
   \bigoplus_{k\geq1}\ \bigoplus_{i\in\mathbb N^k}
   (\mathcal V_k\#e_i).
\]
The last assertion of~\cite[Corollary~15]{Dab14} verifies his
condition $\Gamma_{0u}(\omega)\:e.2)$ for these tensor extensions of $\Delta$: for
$U\in\Dom(\Delta^{\otimes(k+1)})$, one has
$U\#e_i\in\Dom(\widetilde\Delta)$ and
\[
 \widetilde\Delta(U\#e_i)
 =\Delta^{\otimes(k+1)}(U)\#e_i.
\]
Applying this to $U=\bigl(\alpha+\Delta^{\otimes(k+1)}\bigr)^{-1}(W)$ with $\alpha>0$ gives
\[
 (\alpha+\widetilde\Delta)^{-1}(W\#e_i)
 =\bigl((\alpha+\Delta^{\otimes(k+1)}\bigr)^{-1}(W))\#e_i,
 \qquad W\in\mathcal V_k.
\]
Each closed word space is therefore invariant under the self-adjoint
resolvent, and hence reduces $\widetilde\Delta$, whose restriction
there is unitarily equivalent to $\Delta^{\otimes(k+1)}$.
Since $\ker\Delta=\mathbb C.1$, the kernel of $\Delta^{\otimes 2}$ on
$\mathcal V_1$ is $\mathbb C.(1\otimes1)$, whose inserted copies are
the $e_j$'s. For $k\geq2$, a centered internal factor makes the
kernel on $\mathcal V_k$ zero. Consequently, we deduce that
\(
 \ker\widetilde\Delta=L^2(N)
\)
and denoting by $E_N$ the trace-preserving conditional expectation onto $N$, also viewed as the orthogonal projection onto $L^2(N)$,
the spectral theorem gives
\begin{equation}\label{eq:embedding-heat-limit}
 \widetilde\varphi_r\longrightarrow E_N
 \quad\text{strongly on }L^2(\widetilde M)
 \quad\text{as }r\to\infty.
\end{equation}
We now use~\eqref{eq:embedding-integrated-gradient} to pass to the
limit in the mild equation. Fix $T\geq0$. We split the stochastic
integral in~\eqref{eq:embedding-mild} at $T$ and, using the semigroup
identity of~\cite[Lemma~3(ii)]{Dab14}, we obtain, for $t\geq T$,
\[
 Y_t=\widetilde\varphi_{(t-T)/2}(Y_T)
      +\sum_{i=1}^n\int_T^t
        \widetilde\varphi_{(t-s)/2}^\otimes
                       (\bar\delta_iY_s)\#dS_s^{(i)}.
\]
The free It\^o isometry and $L^2$-contractivity of the semigroup associated to $\widetilde\Delta^{\otimes}$ give
\begin{align*}
 \|Y_t-\widetilde\varphi_{(t-T)/2}(Y_T)\|_2^2
 &=\int_T^t
   \|\widetilde\varphi_{(t-s)/2}^\otimes\bar\delta Y_s\|_2^2\,ds\\
 &\leq\int_T^\infty\|\bar\delta Y_s\|_2^2\,ds.
\end{align*}
For fixed $T$, the deterministic term converges to $E_NY_T$
by~\eqref{eq:embedding-heat-limit}. Hence
\begin{equation}\label{eq:embedding-tail}
 \limsup_{t\to\infty}\|Y_t-E_NY_T\|_2
 \leq
 \left(\int_T^\infty\|\bar\delta Y_s\|_2^2\,ds\right)^{1/2}.
\end{equation}
The triangle inequality gives
\(
 \underset{t,t'\to\infty}{\limsup}\|Y_t-Y_{t'}\|_2
 \leq2\left(\int_T^\infty\|\bar\delta Y_s\|_2^2\,ds\right)^{1/2}.
\)
The right-hand side tends to zero as $T\to\infty$, so $Y_t$ is
Cauchy in $L^2(\widetilde M)$. Its limit then belongs to $L^2(N)$,
since~\eqref{eq:embedding-tail} approximates it by $E_NY_T\in L^2(N)$
with an error tending to zero. We thus have proved that
\[
 \Phi_t(x)\longrightarrow\Phi_\infty(x)\in L^2(N),
 \qquad x\in M_0.
\]
The $L^2$ isometry of the maps $\Phi_t$ and density of $M_0$
extend this convergence to every $x\in M$.
Moreover, we see that $E_N\Phi_t(x)\in N$ converges to $\Phi_\infty(x)$ in $L^2$
and has operator norm at most $\|x\|$. Since the 
operator-norm balls in $N$ are closed in the $L^2$ topology, we get
$\Phi_\infty(x)\in N$ and $\|\Phi_\infty(x)\|\leq\|x\|$.
For $x,y\in M$,
\[
 \|\Phi_t(x)\Phi_t(y)-\Phi_\infty(x)\Phi_\infty(y)\|_2\leq
 \|x\|\,\|\Phi_t(y)-\Phi_\infty(y)\|_2
 +\|y\|\,\|\Phi_t(x)-\Phi_\infty(x)\|_2\longrightarrow0.
\]
Since $\Phi_t(xy)=\Phi_t(x)\Phi_t(y)$ and
$\Phi_t(xy)\to\Phi_\infty(xy)$ in $L^2$, uniqueness of the limit gives
\[
 \Phi_\infty(xy)=\Phi_\infty(x)\Phi_\infty(y).
\]
Therefore, we see that multiplicativity is preserved at the limit.
The same is true for linearity, the involution, the unit and the trace.
Injectivity is ensured by preservation of the faithful trace.
The map is also normal: indeed, if $0\leq x_\alpha\uparrow x$, then
the positive element
$\Phi_\infty(x)-\underset{\alpha}{\sup}\:\Phi_\infty(x_\alpha)$
has trace zero by trace preservation and normality of the traces,
hence vanishes by faithfulness.
Finally, the subgroup generated freely by
$\{b^kab^{-k}:k\in\mathbb Z\}\subset\mathbb F_2=\langle a,b\rangle$
gives a trace-preserving inclusion
$L(\mathbb F_\infty)\subset L(\mathbb F_2)$.
Composing the two embeddings gives the conclusion.
\end{proof}

\begin{flushleft}
An immediate consequence of the embedding is that we can apply
Dabrowski's result~\cite[Corollary~25]{Dab14} for tuples with
Lipschitz conjugate variables generating an $R^\omega$--embeddable
algebra: the generating tuple has full microstates free entropy
dimension, as defined by Voiculescu~\cite{Voic2}, providing a further
class of tuples for which the microstates and non-microstates
free entropy dimensions coincide. The interesting point is that,
in our setting, we do not have to assume $R^\omega$--embeddability
separately. We mention, however, that OpenAI's
result~\cite[Corollary~7.1]{openai} shows that $\delta_0$ is not
an invariant of the underlying tracial von Neumann algebra.
Moreover, there exist tuples with finite microstates free entropy
for which $\chi$ and $\chi^*$ do not agree by Open AI's breakthrough result~\cite{openai2}.
\end{flushleft}
\begin{corollary}[Microstates and free entropy dimension]
\label{cor:embedding-microstates}
Under the assumptions of Theorem~\ref{thm:curvature-embedding}, $M$
admits a normal unital trace-preserving embedding into $R^\omega$.
In particular, $X$ has microstates and
\[
 \delta_0(X_1,\ldots,X_n)=n.
\]
\end{corollary}

\begin{proof}The tuple $X$ has Lipschitz conjugate variables, and its von Neumann algebra $M=W^*(X_1,\ldots,X_n)$
is now known to be $R^\omega$--embeddable, so
Dabrowski~\cite[Corollary~25]{Dab14} gives $\delta_0(X)=n$.
\end{proof}

\subsection{Haagerup Approximation Property}
The first consequence of compactness is the Haagerup approximation
property. We recall its definition before giving the direct resolvent
argument.
\begin{definition}[Haagerup~\cite{Haag}, Choda~\cite{Choda}]
A finite von Neumann algebra $(M,\tau)$ with separable predual is said
to have the \emph{Haagerup approximation property} if there exists a
sequence of normal u.c.p. $\tau$-preserving maps $\phi_m:M\to M$ such that
\begin{enumerate}
    \item $\norm{\phi_m(x)-x}_2\to0$ for every $x\in M$,
    \item the operator induced by $\phi_m$ on $L^2(M,\tau)$ is compact.
\end{enumerate}
\end{definition}
\begin{corollary}\label{HAP}
Under the assumptions of Theorem~\ref{main}, the von Neumann algebra $M=W^*(X_1,\ldots,X_n)$ has the Haagerup approximation property.
\end{corollary}

\begin{proof}
For $m\geq1$, set $\rho_m=m(m+\Delta)^{-1}$. By
Proposition~\ref{resolvent}, each map $\rho_m:M\to M$ is normal, u.c.p.
and $\tau$-preserving, and $\norm{\rho_m(x)-x}_2\to0$ for every $x\in M$.
Since $\Delta$ has compact resolvent, the $L^2$-implementation of
$\rho_m$ is compact. Thus $(\rho_m)_{m\geq1}$ satisfies the conditions
of the preceding definition. This argument uses only compactness of
the resolvent and the completely Dirichlet form, not the embedding.
\end{proof}

\begin{remark}
The embedding also gives a very short proof of this property: one can restrict
the Haagerup approximations of $L(\mathbb F_2)$ to $M$ and compose them
with the trace-preserving conditional expectation onto $M$. For comparison with related results, we recall that Guionnet and Shlyakhtenko proved the Haagerup approximation property for stationary free Gibbs laws in their $(c,R)$--convex regime, using a free diffusion argument~\cite[Theorem~4.3]{GS}. Their convexity assumption, imposed to hold for self-adjoint tuples in arbitrary $C^*$--algebras and hence at every matrix amplification, implies our curvature condition at a realization $X$ with $\max_i\norm{X_i}<R$.\footnote{Apply the infinitesimal convexity inequality
of~\cite[Section~1.1]{GS} to self-adjoint finite-rank perturbations
in a faithful normal representation $M\subset B(H)$ and take
the operator trace. By complexification and density, this gives
the same lower bound on $S_2(H)^{\oplus n}$ where $S_2(H)$ denote the space of Hilbert--Schmidt operators on $H$. Finally, the identification
$S_2(H)\simeq H\otimes\bar H$ realizes the Hessian through a faithful spatial representation of $M\bar\otimes M^{op}$,
and therefore gives the desired inequality.} Here we obtain the same property under the $CD(c,\infty)$, $c>0$ condition, directly from compactness of the
resolvent of the free Laplacian.
\end{remark}

\subsection{Solidity}
We collect several stronger operator-algebraic consequences of the
compactness of the free Laplacian. The proofs are short: the closed
free difference quotient is a real derivation into a finite direct sum
of coarse correspondences, and its associated positive generator is
$\Delta=\partial^*\bar\partial$. Theorem~\ref{main} therefore allows us
to apply the proper-derivation results of Ozawa--Popa~\cite{OP1,OP2}.
We first obtain a statement for weakly compact normalizers, from which
solidity follows. Under an additional approximation assumption, the
same argument gives strong solidity. The embedding into the free Brownian algebra will be enough to ensure this CMAP
assumption for our $CD(c,\infty),\: c> 0$ criterion.

\begin{flushleft}
    We recall some standard notation which will be used several times. If $M$ is a finite von Neumann algebra equipped with a distinguished faithful normal tracial state $\tau$, we denote by $\mathcal Z(M)$ its center and by $\mathcal U(M)$ its unitary group. For a von Neumann subalgebra $P\subset M$, we set
    \[
        \mathcal N_M(P):=\left\{u\in\mathcal U(M):(\mathrm{Ad}\:u)(P)=P\right\},
    \]
where $(\mathrm{Ad}\:u)(x)=uxu^*$. We say that a von Neumann subalgebra $P\subset M$ is \emph{regular} if $\mathcal N_M(P)''=M$. A maximal abelian von Neumann subalgebra $A\subset M$ which is regular is then called a \emph{Cartan subalgebra}.
\end{flushleft}

\begin{definition}[Solidity, Ozawa~\cite{Ozawa}]
A finite von Neumann algebra $M$ is \emph{solid} if $Q'\cap M$ is amenable for every diffuse von Neumann subalgebra $Q\subset M$.
\end{definition}

\begin{definition}[Strong solidity, Ozawa--Popa~\cite{OP1}]
A finite von Neumann algebra $M$ is \emph{strongly solid} if $\mathcal N_M(P)''$ is amenable for every diffuse amenable von Neumann subalgebra $P\subset M$.
\end{definition}

We denote by $\bar P$ the conjugate von Neumann algebra of $P$. The map $x\mapsto\bar x$ is conjugate-linear, while $\bar P$ is canonically linearly $*$--isomorphic to $P^{op}$ through
\(
\bar x\mapsto(x^*)^{op}.
\)
If $\sigma$ is an automorphism of $P$, we denote by $\bar\sigma$ the automorphism of $\bar P$ determined by
\(
\bar\sigma(\bar x)=\overline{\sigma(x)}.
\)

\begin{definition}[Weakly compact action, {Ozawa--Popa~\cite[Definition~3.1]{OP1}}]
Let $P$ be a finite von Neumann algebra and let
$\sigma:\Gamma\curvearrowright P$ be a $\tau$-preserving action. The action is called \emph{weakly compact} if there exists a net of unit vectors
\(
\eta_n\in L^2(P\bar\otimes\bar P)_+
\)
such that
\begin{enumerate}
\item $\lVert\eta_n-(v\otimes\bar{v})\eta_n\rVert_2\to0$ for every $v\in\mathcal U(P)$.
\item $\lVert\eta_n-(\sigma_g\otimes\bar{\sigma}_g)(\eta_n)\rVert_2\to0$ for every $g\in\Gamma$.
\item $\langle (x\otimes1)\eta_n,\eta_n\rangle_2=\tau(x)=\langle \eta_n,(1\otimes\bar{x})\eta_n\rangle_2$
for every $x\in P$ and every $n$.
\end{enumerate}
A von Neumann subalgebra $P\subset M$ is called \emph{weakly compact
in $M$} if the conjugation action $\mathcal N_M(P)\curvearrowright P$
is weakly compact. These conditions force $P$ to be amenable.
\end{definition}

We can now use compactness of the free Laplacian together with the proper-derivation normalizer results of Ozawa--Popa~\cite{OP2}. This gives the following statement under the sole compact-resolvent
assumption.

\begin{theorem}[Solidity]\label{solid}
Let $X=(X_1,\ldots,X_n)$ be a finite self-adjoint tuple whose free difference quotients are closable, and assume that the free Laplacian  
\(
\Delta=\partial^*\bar\partial
\)
is \emph{proper}, that is, has compact resolvent. Let $P\subset M=W^*(X_1,\ldots,X_n)$ be a diffuse amenable von Neumann subalgebra, and let $\mathcal G\subset\mathcal N_M(P)$ be a subgroup whose conjugation action on $P$ is weakly compact. Then $\mathcal G''$ is amenable. In particular, $M$ is solid.
\end{theorem}

\begin{proof}
The algebraic free difference quotient
\(
\partial:M_0\to(L^2(M)\bar\otimes L^2(M^{op}))^{\oplus n}
\)
is a real closable derivation, and its closure is $\bar\partial$. Indeed, we have
\(
(\partial_i x)^\dagger=\partial_i(x^*)
\)
for every $x\in M_0$, and by Remark~\ref{real} this identity passes to the closure. The target is a finite direct sum of coarse $M$--$M$ correspondences (hence it is weakly contained in the coarse correspondence). Moreover, the positive operator associated with our derivation is exactly $\Delta=\partial^*\bar\partial$. The first assertion follows directly from Ozawa--Popa~\cite[Remark~5.4]{OP2}, since $\Delta$ has compact resolvent.
\newline
We now recall how solidity is obtained. Take a diffuse von Neumann subalgebra $Q\subset M$, choose a diffuse abelian von Neumann subalgebra $A\subset Q$, and set
\(
\mathcal G:=\mathcal U(Q'\cap M).
\)
Every element of $\mathcal G$ commutes with $A$, and thus the conjugation action $\mathcal G\curvearrowright A$ is trivial. The algebra $A$ is amenable, while the trivial action is compact and hence weakly compact by Ozawa--Popa~\cite[Proposition~3.2]{OP1}. Applying the first part of the
theorem with \(P=A\) shows that $\mathcal G''$ is amenable. But
\[
\mathcal G''=\mathcal U(Q'\cap M)''=Q'\cap M.
\]
Thus $Q'\cap M$ is amenable, which is exactly solidity of $M$.
\end{proof}

With one supplementary approximation assumption, Theorem~\ref{solid}
gives strong solidity. The weak$^*$ completely bounded approximation
property is sufficient; the stronger complete metric approximation
property will follow from the embedding in our positive-curvature
setting.

\begin{definition}[Weak completely bounded approximation property]
A finite von Neumann algebra $M$ has the \emph{(weak$^*$) completely bounded approximation property}, denoted \emph{W$^*$CBAP}, if there exist a constant $C\geq1$ and a net $(\phi_n)_n$ of normal finite-rank maps $\phi_i:M\to M$ such that
\[
\limsup_m\norm{\phi_m}_{cb}\leq C,
\qquad
\norm{x-\phi_m(x)}_2\to0
\]
for every $x\in M$. The \emph{Cowling--Haagerup constant} $\Lambda_{cb}(M)$ is defined as the infimum of the constants $C$ for which such a net exists (the infimum is always achieved; see~\cite{haagcowl} for details). We say that $M$ has the \emph{(weak$^*$) complete metric approximation property}, usually abbreviated CMAP, if
\(
\Lambda_{cb}(M)=1.
\)
By Connes' theorem~\cite{Connes}, amenability implies the \emph{CMAP}.
\end{definition}

\begin{corollary}[CMAP and weak exactness]\label{cor:embedding-cmap}
Under the assumptions of Theorem~\ref{thm:curvature-embedding}, $M$
has CMAP and is weakly exact. The same holds for every von Neumann
subalgebra $Q\subset M$, in particular,
\(
 \Lambda_{\mathrm{cb}}(Q)=1.
\)
\end{corollary}

\begin{proof}
The first part is a well-known fact.
For the second part, $L(\mathbb F_2)$ is weakly exact, and weak exactness passes to the range
of a normal conditional expectation by Ozawa~\cite[Corollary~3 and
the subsequent discussion]{OzawaWeak2}. Hence $Q$ is weakly exact.
\end{proof}

We now use Ozawa's weak compactness theorem, which we recall for convenience.

\begin{theorem}[{Ozawa~\cite[Theorem~B]{OzawaWeak}}]\label{WCBAP}
Let $M$ be a finite von Neumann algebra with the W$^*$CBAP. Then every amenable von Neumann subalgebra $P\subset M$ is weakly compact in $M$.
\end{theorem}

\begin{corollary}[Strong solidity]\label{strong}
Let $X$ and $M$ be as in Theorem~\ref{solid}, and assume moreover that $M$ has the W$^*$CBAP. Then $M$ is strongly solid.
\end{corollary}

\begin{proof}
Take $P\subset M$ diffuse and amenable. By Theorem~\ref{WCBAP}, the conjugation action $\mathcal N_M(P)\curvearrowright P$ is weakly compact. We can therefore apply Theorem~\ref{solid} with $\mathcal G=\mathcal N_M(P)$, and we get that $\mathcal N_M(P)''$ is amenable.
\end{proof}
Since a non-injective solid factor is non-Gamma and prime by Ozawa's result~\cite[Proposition~2]{Ozawa}, we deduce the following:
\begin{corollary}[Strong solidity under positive curvature]\label{prime}
Under the assumptions of Theorem~\ref{main}, $M$ is strongly solid,
as is every von Neumann subalgebra of $M$. If $n\geq2$, then $M$ is
a non-$\Gamma$ prime $\mathrm{II}_1$ factor with no Cartan subalgebra.
Every non-amenable subfactor of $M$ is also non-$\Gamma$, prime
and has no Cartan subalgebra.
\end{corollary}

We mention that primeness itself was already obtained by
Dabrowski--Ioana~\cite[Corollary~1.2]{Dabioa}. Indeed,
Proposition~\ref{prop5} gives bounded first- and second-order conjugate
variables, and their result applies for every $n\geq2$ without the
curvature assumption. More precisely, it shows that $M$ is
non-$L^2$-rigid in the sense of Peterson~\cite{Peterson}.
Moreover, under the Lipschitz conjugate-variable condition,
Dabrowski--Ioana~\cite[Theorem~1.3]{Dabioa} proved that
$M\bar\otimes Q$ has no Cartan subalgebra for every $\mathrm{II}_1$
factor $Q$. What compactness of the free Laplacian gives here in
addition is the weakly compact-normalizer statement. Combined with
the CMAP obtained from the embedding, it gives strong solidity.
In fact, this last statement could also be deduced directly from the embedding and the strong
solidity of free group factors~\cite[Corollary~1]{OP1}.

\subsection{Proper proximality and bi-exactness }

We finish with two further applications of compactness of the free
Laplacian, and then use the embedding to remove the weak exactness
assumption in the positive-curvature case. We use the notation \(\mathbb S(M)\) for the
\emph{small-at-infinity} operator system introduced by
Ding--Kunnawalkam Elayavalli--Peterson~\cite[Section~6]{DKEP}. Recall
that \(M\) is \emph{properly proximal} if \(\mathbb S(M)\) admits no
\(M\)-central state whose restriction to $M$ is normal. Following
Ding--Peterson~\cite[Definition~6.1]{DP}, $M$ is \emph{biexact} if
the inclusion
$
M\subset\mathbb S(M)$ is $M$-nuclear. We refer to~\cite{DKEP,DP} for the complete details.

\begin{corollary}[Proper proximality]\label{proximal}
Let $X$ and $M$ be as in Theorem~\ref{solid}, and assume that $M$
is a non-amenable $\mathrm{II}_1$ factor. Then $M$ is properly proximal.
\end{corollary}

\begin{proof}
The free difference quotient is a closable real derivation into a
finite direct sum of coarse correspondences, and its associated
positive generator is $\Delta=\partial^*\bar\partial$. By the
compact-resolvent assumption,
\[
 \zeta_\alpha=\alpha^{1/2}(\alpha+\Delta)^{-1/2}
 \in\mathbb K(L^2(M)),\qquad\alpha>0.
\]
Non-amenability of $M$ excludes the central state in
Ding--Kunnawalkam Elayavalli--Peterson~\cite[Proposition~9.4 and Lemma~9.5]{DKEP}
and therefore gives proper proximality.
\end{proof}

\begin{corollary}[Bi-exactness]\label{biexact}
Under the assumptions of Theorem~\ref{solid}, suppose moreover that
$M$ is weakly exact. Then $M$ is biexact. In particular, this holds
whenever \(C^*(X_1,\ldots,X_n)\) is exact.
\newline
Moreover, every von Neumann subalgebra $Q\subset M$ is biexact and is
either properly proximal or has a non-zero amenable direct summand. In
particular, every non-amenable subfactor of $M$ is properly proximal.
\end{corollary}

\begin{proof}
For every \(\alpha>0\), compactness of the resolvent gives
\(
\rho_\alpha
:=
\alpha (\alpha+\Delta)^{-1}
\in\mathbb K(L^2(M)).
\)
The first assertion therefore follows directly from
Ding--Peterson~\cite[Theorem~7.11]{DP}, applied to the free difference
quotient case. If $C^*(X_1,\ldots,X_n)$ is exact, then its ultraweak
closure $M$ is weakly exact~\cite{OzawaWeak2}. Since every von
Neumann subalgebra has a trace-preserving conditional expectation,
its bi-exactness  follows from~\cite[Proposition~6.10]{DP}. The
proper-proximality alternative is~\cite[Corollary~6.11]{DP}.
\end{proof}
\begin{corollary}[Bi-exactness  under positive curvature]
\label{cor:curvature-biexact}
Under the assumptions of Theorem~\ref{main}, $M$ is biexact.
If $n\geq2$, it is also properly proximal. More generally, every
von Neumann subalgebra $Q\subset M$ without a non-zero amenable direct
summand is properly proximal.
\end{corollary}

\begin{proof}
By Corollary~\ref{cor:embedding-cmap}, $M$ is weakly exact.
Theorem~\ref{main} gives compact resolvent, so
Corollary~\ref{biexact} applies. When $n\geq2$, $M$ is a non-amenable
factor by Corollary~\ref{prime}; it has no amenable direct summand.
The proper-proximality assertions follow from
\cite[Corollary~6.11]{DP}.
\end{proof}

\section{Some geometric consequences}\label{Sec3}

\subsection{The curvature criterion and a strengthened free Obata theorem}

We now return to the equality case in the free Obata theorem proved in~\cite{DiezObata}. Our two papers interact here in a useful way: compactness of the resolvent guarantees that the bottom of the positive spectrum is attained, while equality in the spectral Bochner estimate identifies the eigenspace at the curvature lower bound.

\begin{flushleft}
The next proposition shows that rigidity already appears exactly at the curvature lower bound. It is the natural free counterpart of the rigidity statement highlighted by Ledoux~\cite[Proposition~4]{Ledoux2} and Bouyrie~\cite[Lemma~5.1]{Bouyrie}.
\end{flushleft}
\begin{proposition}\label{affine}
Under the assumptions of Theorem~\ref{main}, with the same $c>0$, one has
\begin{equation}\label{affinespace}
\Ker(\Delta-c\:\mathrm{id})
\subset
\operatorname{span}_{\C}\{X_1-\tau(X_1).1,\ldots,X_n-\tau(X_n).1\}.
\end{equation}
In particular, $\dim\Ker(\Delta-c\:\mathrm{id})\leq n$.
\end{proposition}

\begin{proof}
Take $f\in\Ker(\Delta-c\:\mathrm{id})$. By spectral calculus, $f\in\Dom(\Delta^r)$ for every $r\geq0$, and Lemma~\ref{dab} gives
\[
c\bar\partial f
=
(\Delta^\otimes\otimes I_n)\bar\partial f
+
\mathcal R_{\mathscr J\xi}\bar\partial f.
\]
Taking the scalar product with $\bar\partial f$ and subtracting $c\norm{\bar\partial f}_2^2$, we get
\[
0
=
\ip{(\Delta^\otimes\otimes I_n)\bar\partial f}{\bar\partial f}_2
+
\ip{(\mathcal R_{\mathscr J\xi}-c\:\mathrm{id})\bar\partial f}{\bar\partial f}_2.
\]
Both terms on the right-hand side are non-negative, and their sum vanishes. Hence $\bar\partial f\in\Ker(\Delta^\otimes\otimes I_n)$. By~\eqref{tensorker}, there exist $a_1,\ldots,a_n\in\C$ such that
\(
\bar\partial_jf=a_j.(1\otimes 1)
\)
for every $j=1,\ldots,n$. We therefore get $f-\sum_ja_jX_j\in\Ker(\bar\partial)=\mathbb C\,1$. Since $c>0$, the identities $\Delta f=cf$ and $\Delta1=0$ imply that $f\perp1$. Centering this affine expression gives~\eqref{affinespace}. If $f=f^*$, all the coefficients $a_j$ are real (see~\cite[Lemma~3.14]{DiezObata} for the complete details).
\end{proof}

\begin{remark}
In particular, every eigenfunction at the curvature lower bound belongs to $M$. This is in perfect analogy with the free Ornstein--Uhlenbeck case, where the corresponding eigenfunctions are the elements of the first Wigner chaos, that is, the linear span of the semicircular coordinates.
\end{remark}

We also recall the free \emph{log-Lichnerowicz} estimate proved in~\cite{DiezObata}. In fact, in this paper we obtained a much more general statement which is a free \emph{Brascamp--Lieb-type} inequality.

\begin{proposition}[{Free log-Lichnerowicz estimate, Diez~\cite[Theorem~2.20]{DiezObata}}]\label{log-lich}
Let $C_P(X)$ be the optimal constant in Voiculescu's free Poincar\'e inequality. If $X$ satisfies $CD(c,\infty)$ for some $c>0$, one has
\[
C_P(X)\leq\frac1c.
\]
\end{proposition}

\begin{corollary}[Free Obata theorem]\label{obata}
Let $X=(X_1,\ldots,X_n)$ have Lipschitz conjugate variables
and satisfy $CD(1,\infty)$. If $C_P(X)=1$, then
\(
 1\leq r:=\dim\Ker(\Delta-\mathrm{id})\leq n.
\)
After centering and an orthogonal change of variables, the
$r$ extremal directions form a standard semicircular family
free from the remaining variables. Consequently,
\[
 M\cong L(\mathbb F_r)*N,
\]
where $N$ is generated by the remaining tuple, which satisfies
$CD(1,\infty)$ if $r<n$; when $r=n$, we set $N=\mathbb C.1$.
If $r\geq2$, we also have
\[
 M\cong L(\mathbb F_2)*N.
\]
Moreover if $r\geq 2$, then $M$ freely absorb $L(\mathbb F_2)$, i.e. $M*L(\mathbb F_2)\cong M$.
All isomorphisms are normal and trace-preserving.

\end{corollary}
\begin{proof}
Since $\Delta$ has compact resolvent and $\Ker(\Delta)=\mathbb C.1$, the bottom $\lambda_1$ of its positive spectrum is an eigenvalue. Moreover, $C_P(X)=\lambda_1^{-1}$, and hence $C_P(X)=1$ gives $\lambda_1=1$. We can therefore apply the multidimensional free Obata theorem of~\cite[Theorem~3.15]{DiezObata}, which gives the desired conclusion. The second isomorphism comes again from Open AI's result~\cite{openai}. Finally, tracial free products are associative and the uniqueness characterization of tracial free products
in~\cite[Proposition~5.3.3]{AnanPopa}, interchanging the order of the two components gives a canonical normal trace-preserving isomorphism. Using Open AI's isomorphism result~\cite{openai} then gives the free absorption by $L(\mathbb F_2)$.
\end{proof}
\begin{flushleft}
The preceding splitting theorem determines the isomorphism class of the entire algebra whenever at most one variable
remains outside the semicircular family. In particular, this
applies to a pair of generators, without any isotropy assumption. The proof is immediate and left to the reader.
\end{flushleft}

\begin{corollary}[An isomorphism consequence of free Obata]
\label{cor:obata-isomorphism}
Let $X=(X_1,\ldots,X_n)$, $n\geq2$ satisfy $CD(1,\infty)$. Suppose that
$C_P(X)=1$ and
\(
 r:=\dim\Ker(\Delta-\mathrm{id})\geq n-1.
\)
Then
\[
 M=W^*(X_1,\ldots,X_n)\cong L(\mathbb F_2)
\]
as tracial von Neumann algebras.
In particular, when $n=2$, under $CD(1,\infty)$
and $C_P(X)=1$ we have
\(
 W^*(X_1,X_2)\cong L(\mathbb F_2).
\)
\end{corollary}

\begin{flushleft}
The results in Section~\ref{operator} and those obtained in this section show that the von Neumann algebras generated
by tuples satisfying $CD(c,\infty)$, $c>0$, share many structural properties with free group factors.
Our Corollary~\ref{cor:obata-isomorphism} identifies the von Neumann algebra
completely when the free product decomposition leaves at most one remaining generator. It is natural to ask whether the same conclusion holds in every equality case, regardless of the dimension of $\Ker(\Delta-\mathrm{id})$.
Together with the recent solution by Open AI~\cite{openai} of the free group factor isomorphism problem in affirmative, we are led to suggest the following conjecture.
\end{flushleft}

\begin{conjecture}
\label{conj:obata-isomorphism}
Let $X=(X_1,\ldots,X_n)$, $n\geq2$, be a bounded self-adjoint
tuple with Lipschitz conjugate variables satisfying
$CD(c,\infty)$ for some $c>0$. If
\(
 C_P(X)=1/c,
\)
then
\(
 W^*(X_1,\ldots,X_n)\cong L(\mathbb F_2)
\).
\end{conjecture}

\begin{flushleft}
More generally, we conjecture that strict positivity of the curvature alone should determine completely the isomorphism class, without any saturation assumption.
\end{flushleft}

\begin{conjecture}
\label{conj:positive-curvature-isomorphism}
Let $X=(X_1,\ldots,X_n)$, $n\geq2$, be a bounded self-adjoint
tuple with Lipschitz conjugate variables satisfying
$CD(c,\infty)$ for some $c>0$. Then
\(
 W^*(X_1,\ldots,X_n)\cong L(\mathbb F_2)
\).
\end{conjecture}
\subsection{An improved free Log-Lichnerowicz bound}
In the classical setting, Klartag~\cite[Theorem~1.3]{Klartaglog}
proved the following improvement for bounding the Poincar\'e constant for uniformly log-concave measures. This refinement is of great importance in connection with the KLS conjecture. We refer to Klartag's paper for further details.

\begin{theorem}[{Improved log-Lichnerowicz inequality, Klartag~\cite[Theorem~1.3]{Klartaglog}}]

If \(\mu\) is a \(c\)-uniformly log-concave probability measure on
\(\mathbb R^n\), with covariance matrix \(\operatorname{Cov}(\mu)\), then
\begin{equation}\label{klart}
C_P(\mu)
\leq
\sqrt{\frac{\lVert\operatorname{Cov}(\mu)\rVert_{\mathrm{op}}}{c}}\leq \frac 1c.
\end{equation}
This improves the usual Lichnerowicz estimate
\(
C_P(\mu)\leq c^{-1},
\)
since uniform log-concavity implies
\(
\lVert\operatorname{Cov}(\mu)\rVert_{\mathrm{op}}\leq c^{-1}.
\)
\end{theorem}
The following result provides a free counterpart of Klartag's
inequality under the curvature condition \(CD(c,\infty)\), and strengthens Proposition~\ref{log-lich}.
\begin{theorem}[Free version of Klartag's log-Lichnerowicz inequality]
Let \(X=(X_1,\ldots,X_n)\) be a self-adjoint \(n\)-tuple satisfying \(CD(c,\infty)\)
for some \(c>0\). Let us denote by
\(
\operatorname{Cov}(X)
=
(
\tau(
(X_i-\tau(X_i).1)(X_j-\tau(X_j).1)
)
)_{i,j=1}^n
\)
the covariance matrix of \(X\), and by \(C_P(X)\) its optimal
free Poincaré constant.
\newline
Then we have
\[
\lVert\operatorname{Cov}(X)\rVert_{\mathrm{op}}
\leq
C_P(X)
\leq
\sqrt{
\frac{\lVert\operatorname{Cov}(X)\rVert_{\mathrm{op}}}{c}
}
\leq
\frac{1}{c}.
\]
\end{theorem}

\begin{proof}
By Theorem~\ref{main}, the operator \(\Delta\) has compact resolvent and
\(
\ker(\Delta)=\mathbb C.1.
\)
Consequently, its minimal non-zero spectral value is an eigenvalue. We denote it by
\(
\lambda
:=
\min(\sigma(\Delta)\setminus\{0\}).
\)
The free Poincaré constant is therefore given by
\(
C_P(X)=\lambda^{-1}.
\)
\newline
We choose a normalized centered eigenfunction \(f\in L^2(M)\):  
\[
\Delta f=\lambda f,
\qquad
\tau(f)=0,
\qquad
\lVert f\rVert_2=1.
\]
Note that $\tau(f)=0$ is immediate as $f$ belongs to the orthogonal complement of $\Ker(\Delta)=\C.1$. Moreover, since \(f\) belongs to an eigenspace of \(\Delta\), it belongs to
\(\operatorname{Dom}(\Delta^r)\) for every \(r>0\). In particular,
\(
f\in\operatorname{Dom}(\Delta^{3/2}).
\)
Dabrowski's almost-commutation formula~\ref{dab} therefore applies to \(f\), and
\(
\bar\partial_i f\in\operatorname{Dom}(\Delta^\otimes),
\)  for all \(i=1,\ldots,n.
\)
\newline
Let \(P_0\) denote the orthogonal projection of \(L^2(M)\) onto
\(\mathbb C.1\), and use the same notation for the corresponding
projection of \(L^2(M^{op})\) onto \(\mathbb C.1\). By a slight abuse of notation, the identity operators on \(L^2(M)\), on \(L^2(M^{op})\),
and on
\(
L^2(M)\bar\otimes L^2(M^{op})
\)
are all denoted by the same symbol \(\mathrm{id}\). By the definition of \(\lambda\)
and the spectral theorem, we have
\[
\Delta\geq\lambda(\mathrm{id}-P_0)
\]
in the sense of quadratic forms. Recalling that
\[
\Delta^\otimes
=
\overline{
\Delta\otimes\mathrm{id}
+
\mathrm{id}\otimes\Delta
},
\]
it follows that
\begin{align}
\Delta^\otimes
&\geq
\lambda(
(\mathrm{id}-P_0)\otimes\mathrm{id}
+
\mathrm{id}\otimes(\mathrm{id}-P_0)
)\nonumber\\
&\geq
\lambda(
\mathrm{id}-(P_0\otimes P_0)
).
\end{align}
Indeed, one can check on the core $\Dom(\Delta)\otimes \Dom(\Delta)$ (here again $\otimes$ denotes the algebraic tensor product) that 
\begin{align}
&
(
(\mathrm{id}-P_0)\otimes\mathrm{id}
+
\mathrm{id}\otimes(\mathrm{id}-P_0)
)
-
(
\mathrm{id}-(P_0\otimes P_0)
)
=
(\mathrm{id}-P_0)\otimes(\mathrm{id}-P_0)
\geq0.\nonumber
\end{align}
Since \(P_0\otimes P_0\), initially defined on simple tensors, extends to the orthogonal projection onto
\(\mathbb C.(1\otimes1)\), we conclude that, for every
\(\eta\in\operatorname{Dom}(\Delta^\otimes)\),
\begin{equation}\label{tensorP}
\langle
\Delta^\otimes\eta,\eta
\rangle_2
\geq
\lambda
(
\lVert\eta\rVert_2^2
-
\lvert
\langle
\eta,1\otimes1
\rangle_2
\rvert^2
).
\end{equation}
The same conclusion can also be seen directly by diagonalizing
\(\Delta\). Indeed, if
\(
0=\lambda_0<\lambda_1=\lambda\leq\lambda_2\leq\cdots
\)
are its eigenvalues, counted with multiplicity, then the eigenvalues
of the tensor extension \(\Delta^\otimes\) (which thus has discrete spectrum) are the sums
\(\lambda_i+\lambda_j\). The only zero sum is for
\(\lambda_0+\lambda_0\), and every other sum is at least \(\lambda\).
\bigbreak
For \(i=1,\ldots,n\), set
\(
m_i
:=
\langle
\bar\partial_i f,1\otimes1
\rangle_2
=:
(\tau\otimes\tau)(\bar\partial_i f)
\).
Applying the preceding estimate to each \(\bar\partial_i f\) and
summing over \(i\), we obtain
\begin{align}
\langle
(\Delta^\otimes\otimes I_n)\bar\partial f,
\bar\partial f
\rangle_2
&\geq
\lambda
\bigg(
\lVert\bar\partial f\rVert_2^2
-
\sum_{i=1}^n\lvert m_i\rvert^2
\bigg).
\end{align}
Since
\(
\lVert\bar\partial f\rVert_2^2
=
\langle\Delta f,f\rangle_2
=
\lambda
\),
this becomes
\[
\langle
(\Delta^\otimes\otimes I_n)\bar\partial f,
\bar\partial f
\rangle_2
\geq
\lambda^2
-
\lambda\sum_{i=1}^n\lvert m_i\rvert^2.
\]
We now apply Dabrowski's almost-commutation formula~\ref{dab}:
\[
\bar\partial\Delta f
=
(
\Delta^\otimes\otimes I_n
+
\mathcal R_{\mathscr J\xi}
)\bar\partial f.
\]
Taking the inner product with \(\bar\partial f\), and using
\(\Delta f=\lambda f\), we get
\[
\lambda^2
=
\langle
(\Delta^\otimes\otimes I_n)\bar\partial f,
\bar\partial f
\rangle_2
+
\langle
\mathcal R_{\mathscr J\xi}\bar\partial f,
\bar\partial f
\rangle_2.
\]
Combining the last two estimates gives
\[
\langle
\mathcal R_{\mathscr J\xi}\bar\partial f,
\bar\partial f
\rangle_2
\leq
\lambda\sum_{i=1}^n\lvert m_i\rvert^2.
\]
Equivalently,
\[
\sum_{i=1}^n\lvert m_i\rvert^2
\geq
\frac{1}{\lambda}
\langle
\mathcal R_{\mathscr J\xi}\bar\partial f,
\bar\partial f
\rangle_2.
\]
This is the free analogue of the ``\emph{preferred-direction}'' estimate in Klartag's terminology (see the discussion after~\cite[Proposition~2.3]{Klartaglog}).
\newline
For \(i=1,\ldots,n\), we denote
\(
X_i^\circ:=X_i-\tau(X_i).1.
\)
Since the free difference quotients are taken with respect to the
variables \(X_1,\ldots,X_n\), we have
\[
\bar\partial_jX_i^\circ
=
\delta_{ij}(1\otimes1).
\]
Integration by parts therefore gives
\begin{align}
m_i
&=
\langle
\bar\partial f,\bar\partial X_i^\circ
\rangle_2\nonumber\\
&=
\langle
\Delta f,X_i^\circ
\rangle_2\nonumber\\
&=
\lambda
\langle
f,X_i^\circ
\rangle_2.
\end{align}
Hence
\[
\sum_{i=1}^n\lvert m_i\rvert^2
=
\lambda^2
\sum_{i=1}^n
\lvert
\langle
f,X_i^\circ
\rangle_2
\rvert^2.
\]
Consider the operator
\[
T:\mathbb C^n\longrightarrow L^2(M),
\qquad
T(a)
:=
\sum_{i=1}^n a_iX_i^\circ,
\qquad
a=(a_1,\ldots,a_n).
\]
We can check for every \(a,b\in\mathbb C^n\), using the self-adjointness of the \(X_i^\circ\), that we have
\begin{align}
\langle T(a),T(b)\rangle_2
&=
\left\langle
\sum_{i=1}^n a_iX_i^\circ,
\sum_{j=1}^n b_jX_j^\circ
\right\rangle_2\nonumber\\
&=
\sum_{i,j=1}^n
a_i\overline{b_j}\,
\tau(X_j^\circ X_i^\circ)\nonumber\\
&=
\left\langle
\operatorname{Cov}(X)a,b
\right\rangle_{\mathbb C^n}.
\end{align}
Hence, the Gram matrix of \(T\) is \(\operatorname{Cov}(X)\). Therefore, we have
\(
T^*T=\operatorname{Cov}(X),
\) so that \(
\lVert T\rVert^2
=
\lVert\operatorname{Cov}(X)\rVert_{\mathrm{op}}
\).
Since \(\lVert f\rVert_2=1\), it follows that
\begin{align}\label{eq21}
\sum_{i=1}^n
\left\lvert
\left\langle
f,X_i^\circ
\right\rangle_2
\right\rvert^2
&=
\lVert T^*f\rVert_{\mathbb C^n}^2\nonumber\\
&\leq
\lVert T\rVert^2\lVert f\rVert_2^2\nonumber\\
&=
\lVert\operatorname{Cov}(X)\rVert_{\mathrm{op}}.
\end{align}
We then get
\(
\sum_{i=1}^n\lvert m_i\rvert^2
\leq
\lambda^2
\lVert\operatorname{Cov}(X)\rVert_{\mathrm{op}}.
\)
Combining this with the preferred-direction estimate gives
\[
\langle
\mathcal R_{\mathscr J\xi}\bar\partial f,
\bar\partial f
\rangle_2
\leq
\lambda^3
\lVert\operatorname{Cov}(X)\rVert_{\mathrm{op}}.
\]
Finally, the \(CD(c,\infty)\) condition implies
\(
\mathcal R_{\mathscr J\xi}
\geq
c\,\mathrm{id}.
\)
Consequently,
\begin{align}
c\lambda
&=
c\lVert\bar\partial f\rVert_2^2\nonumber\\
&\leq
\langle
\mathcal R_{\mathscr J\xi}\bar\partial f,
\bar\partial f
\rangle_2\nonumber\\
&\leq
\lambda^3
\lVert\operatorname{Cov}(X)\rVert_{\mathrm{op}}.
\end{align}
Since \(\lambda>0\), this yields
\[
\lambda^2
\geq
\frac{c}{\lVert\operatorname{Cov}(X)\rVert_{\mathrm{op}}},
\]
and therefore we deduce
\[
C_P(X)
=
\frac{1}{\lambda}
\leq
\sqrt{
\frac{\lVert\operatorname{Cov}(X)\rVert_{\mathrm{op}}}{c}
}.
\]
It remains to establish the lower bound. It is easy to see that testing the free Poincar\'e inequality on linear functions $
T(a):=a_1X_1^{\circ}+\ldots+a_nX_n^{\circ}$ and taking the supremum over all unit vectors \(a\in\mathbb C^n\) gives
\[
\lVert\operatorname{Cov}(X)\rVert_{\mathrm{op}}
\leq
C_P(X).
\]
Finally, the free Lichnerowicz inequality under \(CD(c,\infty)\)
gives
\[
C_P(X)\leq\frac{1}{c}.
\]
Together with
\(
\lVert\operatorname{Cov}(X)\rVert_{\mathrm{op}}\leq C_P(X),
\)
this implies
\[
\sqrt{
\frac{\lVert\operatorname{Cov}(X)\rVert_{\mathrm{op}}}{c}
}
\leq
\frac{1}{c},
\]
and completes the proof.
\end{proof}
The equality case in the preceding theorem admits a rigidity statement
which is the free counterpart of a recent result of Courtade and
Fathi~\cite[Theorem~1.2]{courtfa}. Equality in the improved estimate
necessarily implies that a first eigenfunction is affine in a direction of maximal
variance. It then follows that our usual free Lichnerowicz inequality is
itself saturated, so that our free Obata rigidity
theorem~\cite{DiezObata} gives the conclusion. Courtade-Fathi in fact obtained a stability statement involving the $L^1$-Wasserstein distance~\cite[Theorem~1.3]{courtfa}; we leave this question for further investigations. For convenience, we also restate Courtade-Fathi's result.
\begin{theorem}[{Rigidity for Klartag's log-Lichnerowicz inequality, Courtade-Fathi~\cite[Theorem~1.2]{courtfa}}]
Let \(\mu\) be a \(c\)-log-concave probability measure satisfying~\eqref{klart} with equality. Then \(\mu\) splits off a Gaussian factor in a direction of largest variance. More precisely, up to a rotation and a translation,
\(
\mu=\bar\mu\otimes\gamma_c,
\)
where \(\bar\mu\) is the marginal of \(\mu\) on the first \(d-1\) coordinates and \(\gamma_c\) is the centered one-dimensional Gaussian measure with variance
\(
\frac{1}{c}
=
\|\operatorname{Cov}(\mu)\|_{\mathrm{op}}.
\)
\end{theorem}
\begin{flushleft}
    We first state a proposition about the behavior of the curvature condition under affine
similarities. The proof is elementary and left to the reader; it only uses the chain rule in free differential calculus for affine changes
of variables and the analogous results for
conjugate variables (see Voiculescu~\cite{Voic} and also~\cite{DiezObata} for details and similar computations).
\end{flushleft}

\begin{proposition}[Stability of the $CD$ criterion]
\label{similarity}
Let \(X=(X_1,\ldots,X_n)\) be a self-adjoint tuple with Lipschitz
conjugate variables, let \(a>0\), \(U\in O(n)\), and let
\(b=(b_1,\ldots,b_n)\in\mathbb R^n\). Define
\[
Y_i
:=
a\sum_{j=1}^nU_{ij}X_j+b_i1,
\qquad
i=1,\ldots,n.
\]
Then 
the tuple \(Y\) has Lipschitz conjugate variables, and
\[
C_P(Y)
=
a^2C_P(X),
\qquad
\operatorname{Cov}(Y)
=
a^2U\operatorname{Cov}(X)U^{\mathsf T}.
\]
Moreover, if \(X\) satisfies \(CD(c,\infty)\), then \(Y\) satisfies
\[
CD\left(\frac{c}{a^2},\infty\right).
\]
In particular, translations and orthogonal changes of variables
preserve \(CD(c,\infty)\).
\end{proposition}

\begin{theorem}[Rigidity of the free improved Lichnerowicz inequality]
\label{rigidity}
Let \(X=(X_1,\ldots,X_n)\) be a self-adjoint \(n\)-tuple satisfying \(CD(c,\infty)\) for some
\(c>0\). Let
\(
M=W^*(X_1,\ldots,X_n),
\)
and denote by \(C_P(X)\) its optimal free Poincaré constant. Assume that
equality holds in the free improved log-Lichnerowicz inequality:
\[
C_P(X)
=
\sqrt{
\frac{\lVert\operatorname{Cov}(X)\rVert_{\mathrm{op}}}{c}
}.
\]
Then
\(
\lVert\operatorname{Cov}(X)\rVert_{\mathrm{op}}
=
C_P(X)
=
\frac1c.
\)
Moreover, there exist a unit vector \(\theta\in\mathbb R^n\) and an
orthogonal matrix \(U\in O(n)\), whose first row is \(\theta\), such
that
\(
\operatorname{Cov}(X)\theta
=
\frac1c\,\theta
\)
and, writing
\[
X_j^\circ
:=
X_j-\tau(X_j).1,
\qquad
Y_i
:=
\sum_{j=1}^nU_{ij}X_j^\circ,
\qquad
i=1,\ldots,n,
\]
the following assertions hold:
\begin{enumerate}
\item \(Y_1\) is a centered semicircular variable of variance \(1/c\).

\item \(Y_1\) is free from the family \((Y_2,\ldots,Y_n)\).

\item If $n\geq2$, the tuple \((Y_2,\ldots,Y_n)\) satisfies \(CD(c,\infty)\).
\end{enumerate}
Consequently, setting
\(
N
:=
W^*(Y_2,\ldots,Y_n),
\)
one has the trace-preserving free product decomposition
\(
(M,\tau)
\cong
(W^*(Y_1),\tau_{\mathrm{sc},c})
*
(N,\tau\lvert_N)
\),
where \(\tau_{\mathrm{sc},c}\) is the centered semicircular distribution of
variance \(1/c\).
\end{theorem}

\begin{proof}
By the spectral theorem, since \(\ker(\Delta)=\mathbb C.1\), the free
Poincaré inequality is equivalently written for $f\in\operatorname{Dom}(\Delta)$ as
\(
\lVert\bar\partial f\rVert_2^2
\leq
C_P(X)\lVert\Delta f\rVert_2^2.
\)
Using again the free Bochner formula, for every
\(f\in\operatorname{Dom}(\Delta^{3/2})\),
\begin{align}
\lVert\Delta f\rVert_2^2
&=
\left\langle
(\Delta^\otimes\otimes I_n)
\bar\partial f,\bar\partial f
\right\rangle_2
+
\left\langle
\mathcal R_{\mathscr J\xi}\bar\partial f,
\bar\partial f
\right\rangle_2
\nonumber\\
&\geq
\left\langle
(\Delta^\otimes\otimes I_n)
\bar\partial f,\bar\partial f
\right\rangle_2
+
c\lVert\bar\partial f\rVert_2^2\nonumber\\
&\geq
C_P(X)^{-1}
\left(
\lVert\bar\partial f\rVert_2^2
-
\sum_{i=1}^n
\lvert
\langle
\bar\partial_i f,1\otimes1
\rangle_2
\rvert^2
\right)
+
c\lVert\bar\partial f\rVert_2^2.
\label{bochner}
\end{align}
The last inequality is obtained by applying the spectral-gap estimate
for \(\Delta^\otimes\) to each component \(\bar\partial_i f\), that is,
by applying~\eqref{tensorP}.
\newline
Since \(\lVert\operatorname{Cov}(X)\rVert_{\mathrm{op}}>0\), we may center the
tuple and apply a dilation to normalize its maximal variance to one.
Thus, without loss of generality, we assume that
\(
\tau(X_i)=0,
\) for all \(
i=1,\ldots,n,
\) and \(
\lVert\operatorname{Cov}(X)\rVert_{\mathrm{op}}=1.
\)
\newline
In the sequel, we keep the same notation for the tuple and all the
associated objects. Under this dilation, the free Laplacian and the
curvature constant are multiplied by the original value of
\(\lVert\operatorname{Cov}(X)\rVert_{\mathrm{op}}\), while the Poincaré
constant is divided by the same quantity. We continue to denote the
(new) curvature constant by \(c\). Thus the normalized tuple
satisfies \(CD(c,\infty)\), and the equality assumption becomes
\[
C_P(X)=\frac{1}{\sqrt c}.
\]
Testing the free Poincaré inequality on a linear combination of the
variables in a direction of maximal variance gives
\(1\leq C_P(X)=c^{-1/2}\), and hence \(c\leq1\). The goal is to prove
that \(c=1\). This will show that equality in the improved free
Lichnerowicz inequality implies equality in the usual free Lichnerowicz
inequality, and the conclusion will then follow from the free Obata
rigidity theorem~\cite[Theorem~3.15]{DiezObata}.
\newline
By Theorem~\ref{main}, the spectral gap is attained.

Since
\(\Delta\) commutes with the Tomita conjugation
(see the proof of~\cite[Corollary~3.4]{DiezObata}), we can choose a
self-adjoint eigenfunction \(f\in L^2(M)\) satisfying
\(
\Delta f=\sqrt c\,f,
\) $\tau(f)=0$,
and \(
\lVert f\rVert_2=1
\).
\newline
As
\(
f\in\ker(\Delta-\sqrt c\,\mathrm{id}),
\)
it belongs to \(\operatorname{Dom}(\Delta^r)\) for every \(r>0\). We
also have
\(
\lVert\bar\partial f\rVert_2^2
=
\langle\Delta f,f\rangle_2
=
\sqrt c,
\) and \(
\lVert\Delta f\rVert_2^2=c.
\)
For every \(i=1,\ldots,n\), integration by parts gives, using the self-adjointness of
\(\Delta f\) and \(X_i\),
\[
\left\langle
\bar\partial_i f,1\otimes1
\right\rangle_2
=
\left\langle
\bar\partial f,\bar\partial X_i
\right\rangle_2
=
\left\langle
\Delta f,X_i
\right\rangle_2
\in\mathbb R.
\]
If
\(
\left(
\langle
\bar\partial_i f,1\otimes1
\rangle_2
\right)_{i=1}^n
=0,
\)
then applying~\eqref{bochner} to $f$ would yield
\(c\geq c+c^{3/2}\), which is impossible. We may therefore define
\(
\theta
:=
\frac{
(
\langle
\bar\partial_i f,1\otimes1
\rangle_2
)_{i=1}^n
}{
\left(
\sum_{i=1}^n
\lvert
\langle
\bar\partial_i f,1\otimes1
\rangle_2
\rvert^2
\right)^{1/2}
}.
\)
Then \(\theta\in\mathbb R^n\) and \(\lvert\theta\rvert=1\). Moreover,
\begin{align}
\sum_{i=1}^n
\left\lvert
\left\langle
\bar\partial_i f,1\otimes1
\right\rangle_2
\right\rvert^2
&=
\left\lvert
\sum_{i=1}^n
\theta_i
\left\langle
\bar\partial_i f,1\otimes1
\right\rangle_2
\right\rvert^2\nonumber\\
&=
\left\lvert
\left\langle
\bar\partial f,
\bar\partial
\left(
\sum_{i=1}^n\theta_iX_i
\right)
\right\rangle_2
\right\rvert^2\nonumber\\
&=
\left\lvert
\left\langle
\Delta f,
\sum_{i=1}^n\theta_iX_i
\right\rangle_2
\right\rvert^2 \nonumber\\
&\leq
\lVert\Delta f\rVert_2^2
\left\lVert
\sum_{i=1}^n\theta_iX_i
\right\rVert_2^2
\label{cauchy}\\
&\leq
\lVert\Delta f\rVert_2^2
\lVert\operatorname{Cov}(X)\rVert_{\mathrm{op}}
\lvert\theta\rvert^2
=
\lVert\Delta f\rVert_2^2
=
c.
\label{direction}
\end{align}
Combining this estimate with~\eqref{bochner}, applied to
\(f\), gives
\begin{align}
c
&\geq
c
-\sqrt c \sum_{i=1}^n \lvert \langle \bar\partial_i f,1\otimes1\rangle_2
\rvert^2
+
c^{3/2}\nonumber\\
&\geq
c-c^{3/2}+c^{3/2}
=
c.
\end{align}
Hence, equality holds at every step, and in particular
in~\eqref{cauchy}
and~\eqref{direction}.
\newline
Equality in~\eqref{direction} means that
\(\theta\) is a direction of unit variance, and hence a direction of
maximal variance under the normalization under consideration.
Equivalently,
\[
\operatorname{Cov}(X)\theta=\theta.
\]
Equality in~\eqref{cauchy} implies that
\(
\Delta f
=
\lambda
\sum_{i=1}^n\theta_iX_i
\)
for some \(\lambda\in\mathbb R\). Taking the \(L^2\)-norms of both
sides then gives \(\lvert\lambda\rvert=\sqrt c\). Since
\(\Delta f=\sqrt c\,f\), we deduce that
\(
f
=
\pm
\sum_{i=1}^n\theta_iX_i.
\)
Therefore,
\(
\bar\partial_jf
=
\pm\theta_j(1\otimes1),
\)
for all \(j=1,\ldots,n\), and hence
\[
\sqrt c
=
\lVert\bar\partial f\rVert_2^2
=
\sum_{j=1}^n\theta_j^2
=
1.
\]
Thus \(c=1\), and consequently \(C_P(X)=1\) in the normalized setting.
\newline
The normalized tuple therefore saturates the usual free Lichnerowicz
inequality. The conclusion follows from its equality case, that is, the
free Obata rigidity theorem~\cite[Theorem~3.15]{DiezObata}. Replacing
\(\theta\) by \(-\theta\), if necessary, and removing the normalization,
we obtain
\[
\lVert\operatorname{Cov}(X)\rVert_{\mathrm{op}}
=
C_P(X)
=
\frac1c,
\]
where \(c\) again denotes the curvature constant of the original tuple,
and \(\theta\) is a direction of maximal variance. After an orthogonal
change of variables whose first direction is \(\theta\), the first
centered variable is semicircular with variance \(1/c\), it is free
from the remaining variables, and, when $n\geq2$, the remaining tuple satisfies
\(CD(c,\infty)\) by the conditional-expectation argument in the proof of Corollary~\ref{obata}.
\end{proof}
\subsection{A free KLS bound}

The analogy with classical convex geometry naturally leads us to ask
whether, under $CD(0,\infty)$, the free Poincar\'e constant can be
bounded by a universal multiple of the operator norm of the covariance
matrix. This is the free counterpart of the
Kannan--Lov\'asz--Simonovits (KLS) conjecture~\cite{KLS}, a longstanding
problem in convex geometry, recently resolved affirmatively in the classical setting
by Bizeul, Klartag and Lehec~\cite{bizkla} and by Song and
Zhang~\cite{song}. Let us first recall the classical result.

\begin{theorem}[Bizeul--Klartag--Lehec~\cite{bizkla},
Song--Zhang~\cite{song}]
There exists a universal constant $C>0$ such that, for every dimension
$n\geq1$ and every log-concave random vector $X$ in $\mathbb R^n$,
\[
 \lVert\operatorname{Cov}(X)\rVert_{\mathrm{op}}
 \leq C_P(X)
 \leq C\,\lVert\operatorname{Cov}(X)\rVert_{\mathrm{op}}.
\]
\end{theorem}

\begin{flushleft}
The lower bound follows by testing the Poincar\'e inequality on linear
functions. The upper bound shows that, for log-concave laws, these
functions determine the optimal constant up to a universal factor.
Equivalently, the Poincar\'e constants of isotropic log-concave measures
are bounded independently of the dimension.

We refer to the lecture notes of Klartag and Lehec~\cite{klalehec}
and the survey of Lee and Vempala~\cite{lee} for an account of the
connections between KLS, Bourgain's slicing problem, and the thin-shell
conjecture. In particular, the KLS bound implies the thin-shell
inequality by testing the Poincar\'e inequality on the squared
Euclidean norm.

Following a long series of advances, the slicing problem was resolved
by Klartag and Lehec~\cite{Klalehec1}. Their proof combines stochastic
localization, Milman's $M$-ellipsoids and stability estimates for the
Shannon--Stam inequality. Bizeul~\cite{bizeul} also gave an alternative proof, based on small-ball estimates and avoiding information-theoretic
arguments. Both proofs build on a breakthrough bound of order
$\log\log n$ for the slicing constant, due to Guan~\cite{Guan}.
Klartag and Lehec~\cite{Klalehec2} subsequently resolved the thin-shell
conjecture by constructing a parallel coupling between exponential
tilts of a log-concave measure.

More recently, Chen and Klartag~\cite{Chenkla} proved the sharp
thin-shell inequality
\[
 \operatorname{Var}(|X|^2)\leq8n
\]
for centered isotropic log-concave random vectors. Equality holds when
the coordinates are independent copies of a standard centered
exponential random variable. Their argument uses log-concave moment
measures and the associated Hessian metric (see for example Klartag's paper for a nice exposition on the subject~\cite{Klart}).

The recent proofs of KLS build on these advances.
Bizeul, Klartag and Lehec~\cite{bizkla} start from a third-moment
estimate given by an anisotropic extension of the sharp Chen-Klartag argument due to Bizeul~\cite[Section~2]{bizeul2} and Letwin~\cite{Letwin}. They use stochastic localization to propagate this estimate
to cumulants of all orders (this is closely inspired to the higher-order Stein kernels construction of Fathi~\cite{Fathi}, note that such Stein kernels were also introduced in the free setting by the author~\cite{diez3}). These bounds then allow them to apply a
criterion introduced by Song and Zhang~\cite{song}, which turns control
of averages under exponential tilting into a dimension-free
Poincar\'e inequality.

Song and Zhang~\cite{song} obtain the universal bound by alternating
polynomial variance estimates with curvature-dependent bounds for the
spectral gap.
\end{flushleft}
\begin{flushleft}
We prove below that a universal comparison holds in the free setting
under $CD(0,\infty)$. The proof is inspired by Klartag's analysis of
eigenfunctions with a preferred direction~\cite[Section~2]{Klartaglog}
and the $H^{-1}$ point of view emphasized
in~\cite[Section~4]{Klartaglog}. For a first eigenfunction, the mean
of each of the free difference quotients measures its correlation with linear functions.
If this mean captures a fixed proportion of the row free difference quotient energy, the covariance matrix gives the KLS bound. Under non-negative curvature alone, however, the first combination of the Bochner formula and the Poincar\'e inequality cannot however by itself give this conclusion:
the two estimates can possibly balance exactly, leaving no useful lower bound for this quantity. We implement this eigenfunction heuristic
through the use of \emph{Poisson-type} equations, without assuming that the spectral gap is attained.

The free setting provides in fact an additional step: slicing or taking a partial trace of the free differences quotient, i.e. apply the conditional expectation $\mathrm{id} \otimes \tau$ gives an element of $L^2(M)$ which, under the domain conditions considered below, can itself be differentiated
using Dabrowski's Lemma~\ref{dab}. We study these elements through
the inverse Laplacian, following the $H^{-1}$ point of view mentioned
above. The free Poincar\'e inequality already gives a positive
spectral gap, so that $\Delta^{-1}$ is bounded on $L^2_0(M)$.
For a centered $u\in L^2(M)$, its free $H^{-1}$ norm satisfies
\[
 \|u\|_{H^{-1}}^2
 =
 \langle u,\Delta^{-1}u\rangle_2
 =
 \|\bar\partial\Delta^{-1}u\|_2^2.
\]
After centering the partial traces of the gradient, we solve the
Poisson equations $h_i=\Delta^{-1}a_i$.
Integrating these equations against $a_i$ expresses their squared
$L^2$ norms as pairings between the gradients of $a_i$ and $h_i$,
and hence as pairings involving the free Hessian.

Co-associativity and the real property of free difference quotients now give two representations of the same sum of squared norms through the two outer tensor factors of this Hessian.
The corresponding test gradients occupy different tensor positions,
so their mixed inner product depends only on their partial traces.
We can therefore estimate this mixed term by another application
of Poincar\'e and Bochner. This is exactly where the gain appears: the free Bochner
formula counts derivatives on both tensor factors, whereas the
estimate of the mixed term uses only one. By co-associativity,
these two contributions to the second-order energy are equal.
The means that remain in the computations are then still controlled by the covariance
matrix. Together, these estimates give a strict improvement of the initial Poincar\'e bound, which finally gives the universal comparison.
\end{flushleft}
\begin{theorem}[Free KLS bound]\label{thm:free-KLS}
Let $X=(X_1,\ldots,X_n)$ be a self-adjoint tuple satisfying
$CD(0,\infty)$. Then
\[
\norm{\operatorname{Cov}(X)}_{\mathrm{op}}
\leq C_P(X)
\leq 10\norm{\operatorname{Cov}(X)}_{\mathrm{op}}.
\]
\end{theorem}

\begin{proof}
Voiculescu's free Poincar\'e inequality~\ref{voipub} gives a positive
spectral gap on
\(
L^2_0(M):=L^2(M)\ominus\mathbb C.1.
\)
Writing $\lambda:=\inf\sigma(\Delta|_{L^2_0(M)})$, the spectral theorem gives
\begin{equation}\label{eq:kls-inverse}
\lambda=C_P(X)^{-1}>0,
\qquad
0\leq\Delta^{-1}\leq\frac1\lambda\mathrm{id},
\qquad
\|\Delta^{-1}\|_{\mathbb B(L^2_0(M))}=\frac1\lambda.
\end{equation}
Here and below, $\Delta^{-1}$ acts on $L^2_0(M)$.
We do not assume that $\lambda$ is an eigenvalue.
Choose a centered self-adjoint $f\in\Dom(\Delta^{3/2})$.
Recall that the free Poincar\'e inequality is equivalently written as
$\|\bar\partial f\|_2^2\leq C_P(X)\|\Delta f\|_2^2$.
Our aim is to prove the following estimate:
\begin{equation}\label{eq:kls-self-improvement}
\|\bar\partial f\|_2^2
\leq
\left(\frac{7}{8\lambda}
+\frac54\|\operatorname{Cov}(X)\|_{\mathrm{op}}\right)
\|\Delta f\|_2^2.
\end{equation}
Since $C_P(X)=1/\lambda$, the coefficient $7/8<1$ will allow us to
absorb the term involving $C_P(X)$ and obtain the upper bound in the theorem.
We divide the proof into five steps.

\medskip
\noindent\emph{Step 1:}
Let $P_0$ be the orthogonal projection of $L^2(M)$ onto $\mathbb C.1$.
The spectral theorem gives
\(
\Delta\geq\lambda(\mathrm{id}-P_0),
\) so that
\[
\Delta\otimes\mathrm{id}
\geq\lambda((\mathrm{id}-P_0)\otimes\mathrm{id})
\]
as quadratic forms. If $\eta\in\Dom(\Delta^\otimes)$, the joint
spectral calculus gives $\eta\in\Dom(\Delta\otimes\mathrm{id})$. Hence,
using the tensor Dirichlet form identity
(see, for instance,~\cite[Proposition~3.5]{DiezObata}),
\begin{align}
\lambda\|((\mathrm{id}-P_0)\otimes\mathrm{id})\eta\|_2^2
&\leq\langle(\Delta\otimes\mathrm{id})\eta,\eta\rangle_2\nonumber\\
&=\sum_{j=1}^n
\left\|\overline{\partial_j\otimes\mathrm{id}}\,\eta\right\|_2^2.
\label{eq:kls-first-leg-energy}
\end{align}
For a simple tensor,
\(
(P_0\otimes\mathrm{id})(x\otimes y)
=\tau(x).1\otimes y
=1\otimes(\tau\otimes\mathrm{id})(x\otimes y).
\)
Both sides extend continuously to the Hilbert tensor product.
Orthogonality therefore gives
\[
\|\eta\|_2^2
=\|(\tau\otimes\mathrm{id})\eta\|_2^2
+\|((\mathrm{id}-P_0)\otimes\mathrm{id})\eta\|_2^2.
\]
Together with~\eqref{eq:kls-first-leg-energy}, this proves
\begin{equation}\label{eq:kls-tensor-Poincare}
\|\eta\|_2^2
\leq\|(\tau\otimes\mathrm{id})\eta\|_2^2
+\frac1\lambda\sum_{j=1}^n
\left\|\overline{\partial_j\otimes\mathrm{id}}\,\eta\right\|_2^2.
\end{equation}
Dabrowski's almost commutation formula in Lemma~\ref{dab}
gives, for
$u\in\Dom(\Delta^{3/2})$,
\begin{equation}\label{eq:kls-commutation}
\bar\partial u\in\Dom(\Delta^\otimes\otimes I_n),
\qquad
\bar\partial\Delta u
=(\Delta^\otimes\otimes I_n
+\mathcal R_{\mathscr J\xi})\bar\partial u.
\end{equation}
Define the tensor Dirichlet energy by
\begin{align*}
\mathcal E_2(u)
&:=\langle(\Delta^\otimes\otimes I_n)\bar\partial u,
                   \bar\partial u\rangle_2\\
&=\sum_{i,j=1}^n
\left\|\overline{\partial_i\otimes\mathrm{id}}\,\bar\partial_j u\right\|_2^2
+\sum_{i,j=1}^n
\left\|\overline{\mathrm{id}\otimes\partial_j}\,\bar\partial_i u\right\|_2^2.
\end{align*}
Taking the inner product of~\eqref{eq:kls-commutation} with
$\bar\partial u$, and using $CD(0,\infty)$, yields
\begin{equation}\label{eq:kls-Bochner}
\mathcal E_2(u)
=\|\Delta u\|_2^2
-\langle\mathcal R_{\mathscr J\xi}\bar\partial u,
                         \bar\partial u\rangle_2
\leq\|\Delta u\|_2^2.
\end{equation}
Co-associativity also holds on this domain:
\begin{equation}\label{eq:kls-coass}
\overline{\partial_i\otimes\mathrm{id}}\,\bar\partial_j u
=\overline{\mathrm{id}\otimes\partial_j}\,\bar\partial_i u,
\qquad 1\leq i,j\leq n.
\end{equation}
Indeed, approximate $u$ by polynomials in the graph norm of
$\bar\partial$ and test the polynomial identity against
$V\in M_0^{\otimes3}$. Voiculescu's formula for the adjoint~\cite[Proposition~4.1]{Voic} and density of $M_0^{\otimes3}$ then
give~\eqref{eq:kls-coass}.
In particular, the two tensor energies are equal:
\begin{equation}\label{eq:kls-half}
\sum_{i,j=1}^n
\left\|\overline{\partial_i\otimes\mathrm{id}}\,\bar\partial_j u\right\|_2^2
=\sum_{i,j=1}^n
\left\|\overline{\mathrm{id}\otimes\partial_j}\,\bar\partial_i u\right\|_2^2
=\frac12\mathcal E_2(u).
\end{equation}

\medskip
\noindent\emph{Step 2: Reduction to the centered partial traces of the free gradients}
For the function $f$ fixed above, set
\[
m_i:=(\tau\otimes\tau)(\bar\partial_i f),
\qquad
a_i:=(\mathrm{id}\otimes\tau)(\bar\partial_i f)-m_i.1,
\qquad
A:=\sum_{i=1}^n\|a_i\|_2^2.
\]
Thus $\tau(a_i)=0$. Moreover, $a_i\in\Dom(\Delta)$.
This follows from Dabrowski's Lemma~\ref{dab}(ii).
One can also see it directly: $\mathrm{id}\otimes P_0$ commutes with
the spectral projections of
$\Delta^\otimes$ and its range is $L^2(M)\bar\otimes\mathbb C.1$.
The restriction of $\Delta^\otimes$ to this range is $\Delta\otimes\mathrm{id}$.
Applying this observation to $\bar\partial_i f\in\Dom(\Delta^\otimes)$
gives the conclusion.
\newline
The real property of free difference quotients gives
\(
\bar\partial_i(f^*)=(\bar\partial_i f)^\dagger.
\)
Since $f=f^*$, we have $m_i\in\mathbb R$. Using Mai's result~\cite[Lemma~4.6]{Mai} again, for any $u\in L^2(M)\bar\otimes L^2(M^{op})$ one has 
\(
(\mathrm{id}\otimes \tau)(u)^\ast = (\tau\otimes \mathrm{id})(u^\dagger)
\), so that
\begin{equation}\label{eq:kls-opposite-trace}
(\tau\otimes\mathrm{id})(\bar\partial_i f)=a_i^*+m_i\,1.
\end{equation}
In particular, by orthogonality, since $\tau(a_i)=0$, we have $a_i^*\perp1$ in $L^2(M)$.
Pythagoras' theorem and traciality therefore give
\[
\|(\tau\otimes\mathrm{id})(\bar\partial_i f)\|_2^2
=\|a_i\|_2^2+|m_i|^2.
\]
Apply~\eqref{eq:kls-tensor-Poincare} to $\bar\partial_i f$ and sum
over $i=1,\ldots,n$. By~\eqref{eq:kls-half},
\begin{equation}\label{eq:kls-reduction}
\|\bar\partial f\|_2^2
\leq\sum_{i=1}^n|m_i|^2+A+\frac1{2\lambda}\mathcal E_2(f).
\end{equation}
Writing $X_i^\circ:=X_i-\tau(X_i).1$, integration by parts gives
$m_i=\langle\Delta f,X_i^\circ\rangle_2$.
Since $\operatorname{Cov}(X)$ is the Gram matrix of the $X_i^\circ$,
Cauchy--Schwarz yields
\begin{equation}\label{eq:kls-covariance-means}
\sum_{i=1}^n|m_i|^2
\leq\|\operatorname{Cov}(X)\|_{\mathrm{op}}\|\Delta f\|_2^2.
\end{equation}
It remains to estimate $A$. Applying Poincar\'e directly to each
$a_i$ would give no improvement in~\eqref{eq:kls-reduction}.
We instead use the inverse Laplacian.

\medskip
\noindent\emph{Step 3: The Poisson equation and exterior Hessian traces.}
To express $\|a_i\|_2^2$ by integration by parts, we solve the Poisson
equation. For each $i=1,\ldots,n$, let
\[
h_i:=\Delta^{-1}a_i\in L^2_0(M),
\qquad \Delta h_i=a_i.
\]
Integration by parts and~\eqref{eq:kls-inverse} give
\begin{equation}\label{eq:kls-Hminus-bound}
\sum_{i=1}^n\|\bar\partial h_i\|_2^2
=\sum_{i=1}^n\langle a_i,\Delta^{-1}a_i\rangle_2
\leq\frac A\lambda.
\end{equation}
Each term $\langle a_i,\Delta^{-1}a_i\rangle_2$ is precisely the squared
$H^{-1}$ norm of $a_i$. Thus~\eqref{eq:kls-Hminus-bound} controls the
energy of the Poisson solutions, while $\Delta h_i=a_i$ gives
$\|a_i\|_2^2=\langle\bar\partial a_i,\bar\partial h_i\rangle_2$.
This is the identity we will use to pair the Hessian of $f$ with the
gradients of the $h_i$.
\newline
Since $a_i\in\Dom(\Delta)$, in particular
$a_i\in\Dom(\Delta^{1/2})$. On $L^2_0(M)$, spectral calculus gives
$\Delta^{3/2}\Delta^{-1}a_i=\Delta^{1/2}a_i$, so
\(
h_i\in\Dom(\Delta^{3/2})
\).
The free Bochner formula and co-associativity therefore apply to each $h_i$.
\newline
Define the Hessian
\begin{equation}\label{eq:kls-Hessian}
H_{ij}:=\overline{\partial_i\otimes\mathrm{id}}\,\bar\partial_j f
=\overline{\mathrm{id}\otimes\partial_j}\,\bar\partial_i f,\qquad 1\leq i,j\leq n.
\end{equation}
By Dabrowski~\cite[Lemma~18]{Dab14}, the ``sliced'' derivatives, i.e. composing the conditional expectation $\mathrm{id}\otimes\tau$ with the free difference quotients:
$(\mathrm{id}\otimes\tau)\circ \partial_j$, extend to bounded operators
on $L^2(M)$. The same holds for $(\tau\otimes\mathrm{id})\circ\partial_i$
by applying the involution.
As in~\cite[Lemma~19(ii)]{Dab14}, co-associativity on polynomials,
approximation in the graph norm of $\bar\partial$, and closedness
of $\bar\partial$ therefore give
\begin{align}
\bar\partial_i\bigl((\mathrm{id}\otimes\tau)\bar\partial_j f\bigr)
&=\bigl(\mathrm{id}\otimes
  [(\mathrm{id}\otimes\tau)\bar\partial_j]\bigr)\bar\partial_i f
\nonumber\\
&=(\mathrm{id}\otimes\mathrm{id}\otimes\tau)H_{ij},
\label{eq:kls-trace-derivative-right}\\
\bar\partial_j\bigl((\tau\otimes\mathrm{id})\bar\partial_i f\bigr)
&=\bigl([(\tau\otimes\mathrm{id})\bar\partial_i]
  \otimes\mathrm{id}\bigr)\bar\partial_j f
\nonumber\\
&=(\tau\otimes\mathrm{id}\otimes\mathrm{id})H_{ij}.
\label{eq:kls-trace-derivative-left}
\end{align}
Here the ``bracketed'' compositions of conditional expectations with
derivatives denote their bounded extensions on $L^2(M)$.
For the last equalities, the corresponding identities hold on tensors of non-commutative polynomials, which form a core for the closed operators $\overline{\mathrm{id}\otimes\partial_j}$ and
$\overline{\partial_i\otimes\mathrm{id}}$
in~\eqref{eq:kls-Hessian}. They therefore extend by approximation of
$\bar\partial_i f$ and $\bar\partial_j f$ in the respective graph norms,
using boundedness of these extensions and contractivity of
the $L^2$-extensions of the conditional expectations
$\mathrm{id}\otimes\mathrm{id}\otimes\tau$ and
$\tau\otimes\mathrm{id}\otimes\mathrm{id}$.
Since constants have zero derivative, the definition of $a_j$
and~\eqref{eq:kls-opposite-trace} yield
\begin{align}
(\mathrm{id}\otimes\mathrm{id}\otimes\tau)H_{ij}
&=\bar\partial_i a_j,
\nonumber\\
(\tau\otimes\mathrm{id}\otimes\mathrm{id})H_{ij}
&=\bar\partial_j(a_i^*).
\label{eq:kls-two-outer-traces}
\end{align}
The adjoints of these conditional expectation maps are respectively given by
$x\mapsto x\otimes1$ and $x\mapsto1\otimes x$.
Thus, by integration by parts, we see immediately that
\begin{align*}
A
&=\sum_{j=1}^n\langle a_j,\Delta h_j\rangle_2\\
&=\sum_{i,j=1}^n\langle\bar\partial_i a_j,
                              \bar\partial_i h_j\rangle_2\\
&=\sum_{i,j=1}^n\langle H_{ij},\bar\partial_i h_j\otimes1\rangle_2.
\end{align*}
Since $\Delta$
commutes with the Tomita involution, we have $\Delta(h_i^*)=a_i^*$. Consequently, using traciality,
\begin{align*}
A
&=\sum_{i=1}^n\langle a_i^*,\Delta(h_i^*)\rangle_2\\
&=\sum_{i,j=1}^n\langle\bar\partial_j(a_i^*),
                              \bar\partial_j(h_i^*)\rangle_2\\
&=\sum_{i,j=1}^n\langle H_{ij},1\otimes\bar\partial_j(h_i^*)\rangle_2.
\end{align*}
Adding the two representations gives
\begin{equation}\label{eq:kls-double-duality}
2A=\sum_{i,j=1}^n
\left\langle H_{ij},
\bar\partial_i h_j\otimes1+1\otimes\bar\partial_j(h_i^*)
\right\rangle_2.
\end{equation}
Each outer trace of the Hessian produces one copy of $A$.
Co-associativity allows us to combine them in the same scalar product.
Since $\sum_{i,j=1}^n\|H_{ij}\|_2^2=\mathcal E_2(f)/2$
by~\eqref{eq:kls-half}, Cauchy--Schwarz yields
\begin{equation}\label{eq:kls-double-CS}
4A^2\leq\frac12\mathcal E_2(f)
\sum_{i,j=1}^n
\left\|\bar\partial_i h_j\otimes1
+1\otimes\bar\partial_j(h_i^*)\right\|_2^2.
\end{equation}

\medskip
\noindent\emph{Step 4: Estimating the mixed term.}
A direct calculation on simple tensors, extended by linearity and
continuity, gives the following identity (recall that our inner products
are linear in the first variable): 
\(
\langle x\otimes1,1\otimes y\rangle_2
=\left\langle(\tau\otimes\mathrm{id})x,
                   (\mathrm{id}\otimes\tau)y\right\rangle_2.
\)

Set
\(
s_{ij}:=(\tau\otimes\mathrm{id})(\bar\partial_i h_j).
\)
Using Mai's result~\cite[Lemma~4.6]{Mai} once more gives
\(
(\mathrm{id}\otimes\tau)\bar\partial_j(h_i^*)=s_{ji}^*.
\)
Hence, Cauchy--Schwarz implies
\[
\left|
\sum_{i,j=1}^n
\left\langle\bar\partial_i h_j\otimes1,
                   1\otimes\bar\partial_j(h_i^*)\right\rangle_2
\right|
=\left|\sum_{i,j=1}^n\langle s_{ij},s_{ji}^*\rangle_2\right|
\leq\sum_{i,j=1}^n\|s_{ij}\|_2^2.
\]
Expanding the squared norms, and using
$\|\bar\partial(h_j^*)\|_2=\|\bar\partial h_j\|_2$ (since $^\dagger$ is an isometric involution), we obtain 
\begin{equation}\label{eq:kls-mixed-bound}
\sum_{i,j=1}^n
\left\|\bar\partial_i h_j\otimes1
+1\otimes\bar\partial_j(h_i^*)\right\|_2^2\leq2\sum_{j=1}^n\|\bar\partial h_j\|_2^2
+2\sum_{i,j=1}^n\|s_{ij}\|_2^2.
\end{equation}
The first sum is bounded by $A/\lambda$ by
\eqref{eq:kls-Hminus-bound}.

For the second sum, the domains in~\eqref{eq:kls-commutation} and
the trace identity~\eqref{eq:kls-trace-derivative-left}, applied to $h_j$, give
$s_{ij}\in\Dom(\bar\partial)$, with
\[
\bar\partial_k s_{ij}
=(\tau\otimes\mathrm{id}\otimes\mathrm{id})
\overline{\mathrm{id}\otimes\partial_k}\,\bar\partial_i h_j.
\]
The partial trace is contractive. Thus, using~\eqref{eq:kls-half},
\eqref{eq:kls-Bochner}, and $\Delta h_j=a_j$, we obtain
\begin{align}
\sum_{i,j,k=1}^n\|\bar\partial_k s_{ij}\|_2^2
&\leq\sum_{i,j,k=1}^n
\left\|\overline{\mathrm{id}\otimes\partial_k}\,\bar\partial_i h_j\right\|_2^2
\nonumber\\
&=\frac12\sum_{j=1}^n\mathcal E_2(h_j)
\leq\frac12\sum_{j=1}^n\|a_j\|_2^2
=\frac A2.
\label{eq:kls-s-derivative}
\end{align}
The $h_j$ need not be self-adjoint: the Bochner estimate
and~\eqref{eq:kls-half} hold for arbitrary complex $u$.
\newline
The means of the $s_{ij}$ are also controlled. Integration by parts gives
\begin{align*}
\tau(s_{ij})
&=(\tau\otimes\tau)(\bar\partial_i h_j)\\
&=\langle\bar\partial h_j,\bar\partial X_i^\circ\rangle_2\\
&=\langle\Delta h_j,X_i^\circ\rangle_2
=\langle a_j,X_i^\circ\rangle_2.
\end{align*}
As in~\eqref{eq:kls-covariance-means}, this yields
\begin{equation}\label{eq:kls-s-means}
\sum_{i,j=1}^n|\tau(s_{ij})|^2
\leq\|\operatorname{Cov}(X)\|_{\mathrm{op}}A.
\end{equation}
We now apply the free Poincar\'e inequality to the centered family $s_{ij}-\tau(s_{ij})\,1$.  Using
\eqref{eq:kls-s-derivative} and~\eqref{eq:kls-s-means} gives
\begin{align}
\sum_{i,j=1}^n\|s_{ij}\|_2^2
&\leq\frac1\lambda
\sum_{i,j,k=1}^n\|\bar\partial_k s_{ij}\|_2^2
+\sum_{i,j=1}^n|\tau(s_{ij})|^2\nonumber\\
&\leq\left(\frac1{2\lambda}+\|\operatorname{Cov}(X)\|_{\mathrm{op}}\right)A.
\label{eq:kls-s-bound}
\end{align}
Combining this with~\eqref{eq:kls-Hminus-bound} and
\eqref{eq:kls-mixed-bound} yields
\begin{equation}\label{eq:kls-second-factor}
\sum_{i,j=1}^n
\left\|\bar\partial_i h_j\otimes1
+1\otimes\bar\partial_j(h_i^*)\right\|_2^2
\leq\left(\frac3\lambda+2\|\operatorname{Cov}(X)\|_{\mathrm{op}}\right)A.
\end{equation}
Substituting into~\eqref{eq:kls-double-CS}, we obtain
\[
4A^2\leq\frac12\mathcal E_2(f)
\left(\frac3\lambda+2\|\operatorname{Cov}(X)\|_{\mathrm{op}}\right)A.
\]
Dividing by $4A$ when $A>0$ gives the following bound, which also
holds trivially when $A=0$:
\begin{equation}\label{eq:kls-A-bound}
A\leq\left(\frac{3}{8\lambda}
+\frac14\|\operatorname{Cov}(X)\|_{\mathrm{op}}\right)\mathcal E_2(f).
\end{equation}

\medskip
\noindent\emph{Step 5: Self-improvement of the Poincar\'e constant.}
Combining~\eqref{eq:kls-reduction},
\eqref{eq:kls-covariance-means}, and~\eqref{eq:kls-A-bound} gives
\[
\|\bar\partial f\|_2^2
\leq\|\operatorname{Cov}(X)\|_{\mathrm{op}}\|\Delta f\|_2^2+\left(\frac{3}{8\lambda}
+\frac14\|\operatorname{Cov}(X)\|_{\mathrm{op}}\right)\mathcal E_2(f)
+\frac1{2\lambda}\mathcal E_2(f).
\]
By~\eqref{eq:kls-Bochner}, we get
\[
\|\bar\partial f\|_2^2
\leq\left(\frac{7}{8\lambda}
+\frac54\|\operatorname{Cov}(X)\|_{\mathrm{op}}\right)\|\Delta f\|_2^2,
\]
which is~\eqref{eq:kls-self-improvement}.

We have proved this for centered self-adjoint
$f\in\Dom(\Delta^{3/2})$, and both sides are unchanged by adding a scalar.
Since $\Delta$ commutes with the Tomita involution $f^*:=Jf$
(see~\cite[Corollary~3.4]{DiezObata}), its spectral calculus does as well.
For $f=\Re f+i\Im f$, with $\Re f=(f+f^*)/2$ and
$\Im f=(f-f^*)/(2i)$, the real property of the free difference quotient, stability under adjunction~\ref{adjunc} and the reality of the corresponding inner products give
$\|\bar\partial f\|_2^2=\|\bar\partial\Re f\|_2^2
+\|\bar\partial\Im f\|_2^2$, and the same decomposition holds for
$\|\Delta f\|_2^2$. Thus~\eqref{eq:kls-self-improvement} holds for every
$f\in\Dom(\Delta^{3/2})$.
For $f\in\Dom(\Delta)$, consider the spectral truncations
\[
f_N:=\mathbf1_{[0,N]}(\Delta)f.
\]
Then $f_N\in\Dom(\Delta^r)$ for every $r>0$, and
$f_N\to f$, $\Delta f_N\to\Delta f$ in $L^2(M)$.
Moreover,
\[
\|\bar\partial(f_N-f)\|_2^2
=\langle\Delta(f_N-f),f_N-f\rangle_2\longrightarrow0.
\]
Thus~\eqref{eq:kls-self-improvement} extends to all
$f\in\Dom(\Delta)$.
\newline
Now take $g\in L^2_0(M)$ and apply this estimate to
$f=\Delta^{-1}g$. Integration by parts gives
\begin{align*}
\langle g,\Delta^{-1}g\rangle_2
&=\|\bar\partial\Delta^{-1}g\|_2^2\\
&\leq\left(\frac{7}{8\lambda}
+\frac54\|\operatorname{Cov}(X)\|_{\mathrm{op}}\right)\|g\|_2^2.
\end{align*}
Taking the supremum over unit vectors $g$ and using the positivity of
$\Delta^{-1}$ yields
\[
\frac1\lambda
=\|\Delta^{-1}\|
\leq\frac{7}{8\lambda}+\frac54\|\operatorname{Cov}(X)\|_{\mathrm{op}}.
\]
Rearranging gives
\(
\frac1{8\lambda}\leq\frac54\|\operatorname{Cov}(X)\|_{\mathrm{op}},
\)
and hence
\[
C_P(X)=\frac1\lambda\leq10\,\|\operatorname{Cov}(X)\|_{\mathrm{op}},
\]
which concludes the proof.
\end{proof}

The free KLS bound also yields a counterpart of the thin-shell
theorem of Klartag--Lehec~\cite{Klalehec2} (see again Chen--Klartag~\cite{Chenkla} for the optimal bound). In the classical
setting, this theorem states that the Euclidean norm of a centered
isotropic log-concave random vector remains in a thin spherical shell at a universally bounded
$L^2$ distance from $\sqrt n$. In our setting, the Euclidean norm
is replaced by the positive square root of the sum of squares.
The corresponding estimate follows by applying the free
Poincar\'e inequality to this sum.
\begin{corollary}[A free thin-shell estimate]
\label{cor:free-thin-shell}
Let $X=(X_1,\ldots,X_n)$ be a centered isotropic self-adjoint
tuple satisfying $CD(0,\infty)$.
Let $\mathbb{X}:=\sum_{i=1}^n X_i^2$. Then we have
\[
\tau\bigl((\mathbb{X}-n.1)^2\bigr)\leq2Cn,
\qquad
\tau\bigl((\mathbb{X}^{1/2}-\sqrt n.1)^2\bigr)\leq2C,
\]
where $C$ denotes the optimal universal constant in the free
KLS inequality.
\end{corollary}

\begin{proof}
Since $X$ is centered and isotropic, $\tau(\mathbb{X})=n$ and
\[
\partial_j\mathbb{X}=X_j\otimes1+1\otimes X_j,
\qquad
\|\partial_j\mathbb{X}\|_2^2=2.
\]
The free Poincar\'e inequality and the definition of $C$ give
\[
\tau\bigl((\mathbb{X}-n.1)^2\bigr)
\leq C_P(X)\sum_{j=1}^n\|\partial_j\mathbb{X}\|_2^2
=2nC_P(X)\leq2Cn.
\]
Moreover, functional calculus yields
\[
\begin{aligned}
(\mathbb{X}^{1/2}-\sqrt n.1)^2
&=(\mathbb{X}-n.1)^2
  (\mathbb{X}^{1/2}+\sqrt n.1)^{-2}\\
&\leq\frac1n(\mathbb{X}-n.1)^2.
\end{aligned}
\]
Taking the trace proves the second assertion and concludes the proof.
\end{proof}
\begin{remark}
We do not address the optimal constants in the free KLS and thin-shell inequalities in this paper. Various computations
nevertheless suggest intriguing parallels between the classical and free theories, even at the level of sharp constants and
extremal behaviour. We leave these questions for further investigation.
\end{remark}
\subsection{Further questions on spectral comparison}\label{milman}

The compactness of the resolvent naturally raises the question of a
comparison for the entire spectrum. In the classical Euclidean setting,
E. Milman~\cite[Theorem~1.1]{MilmanSpectral} proved that the usual $CD(c,\infty)$ criterion on Euclidean spaces implies a very strong spectral comparison using Caffarelli's contraction theorem~\cite{cafa}
and the min--max principle.

\begin{theorem}[{E. Milman's spectral comparison, \cite[Theorem~1.1]{MilmanSpectral}}]
    Let $(\mathbb R^n,\lvert\cdot\rvert,\mu)$ satisfy $CD(c,\infty)$ for some \(c>0\), where $\lvert\cdot\rvert$ is the usual Euclidean norm, i.e., \[
d\mu=e^{-V}dx
\qquad\text{and}\qquad
\nabla^2V\geq cI_n.
\] Then, for the eigenvalues arranged in non-decreasing order, we have
\[
\lambda_k(\mathbb R^n,\lvert\cdot\rvert,\mu)
\geq
\lambda_k(\mathbb R^n,\lvert\cdot\rvert,\gamma_c^n),
\qquad
k\geq0,
\]
where \(\gamma_c^n\) is the centered Gaussian measure with covariance
\(\frac{1}{c}\mathrm {I}_n\), and the eigenvalues are repeated according to
multiplicity.

\end{theorem}
\begin{flushleft}
Let \(S_c=(S_{c}^1,\ldots,S_{c}^n)\) be a centered semicircular family
with covariance \(c^{-1}I_n\), and let \(\mathcal L_n\) denote the
number operator associated with the usual Wigner-chaos decomposition, normalized
so that
\[
\mathcal L_n\lvert_{\mathcal H_k}=k\:\mathrm{id},
\qquad
k\geq0,
\]
where \(\mathcal H_k\) denotes the \(k\)-th \emph{Wigner chaos} (see Biane-Speicher~\cite{BS} for details). The free
Laplacian associated with \(S_c\) is then
\[
\Delta_{S_c}=c\mathcal L_n.
\]
Consequently,
\[
\sigma(c\mathcal L_n)=c\mathbb N_0,
\qquad
\dim\Ker(c\mathcal L_n-mc\:\mathrm{id})=n^m,
\qquad
m\geq0.
\]
By Theorem~\ref{main}, for every tuple
\(X=(X_1,\ldots,X_n)\) satisfying \(CD(c,\infty)\), the corresponding
free Laplacian \(\Delta\) has compact resolvent. Its eigenvalues may
therefore be arranged in non-decreasing order and repeated according
to multiplicity.
\end{flushleft}
\begin{conjecture}[Free version of E. Milman's comparison]
Let \(X=(X_1,\ldots,X_n)\)
satisfy \(CD(c,\infty)\) for some $c>0$ and denote its free Laplacian by $\Delta=\partial^*\bar\partial$. Then we have
\[
\lambda_k(\Delta)
\geq
\lambda_k(c\mathcal L_n),
\qquad
k\geq0.
\]
\end{conjecture}
\begin{flushleft}
Equivalently, this question can be reformulated in this way: does one have
\begin{equation}\label{milcount}
N_{\Delta}(t)
\leq
N_{c\mathcal L_n}(t)
=
\sum_{m=0}^{\lfloor t/c\rfloor}n^m,
\qquad
t\geq0?
\end{equation}
\end{flushleft}
\begin{flushleft}
The preceding Corollary~\ref{expon} gives the same comparison up to Catalan
multiplicities.
A positive answer would in particular improve our heat-trace estimate.
\end{flushleft}
\begin{flushleft}
    There is at least a positive answer to the above question in dimension one, that is, the commutative case (under some mild
regularity assumptions). This is possible thanks to the free Caffarelli contraction theorem
proved in our previous work~\cite{Di4}. Note below that both measures $\nu$ and $\eta$, the latter being the semicircular law, are absolutely continuous with compact support, so the Brenier map exists.
\end{flushleft}

\begin{proposition}[One-dimensional free spectral comparison]
\label{milman1}
Let \(\nu_u\) be the free Gibbs measure associated with a smooth
potential \(u:\mathbb R\to\mathbb R\) satisfying
\(
u''\geq c>0
\).
Assume that the optimal transport map \(T\) sending the standard
semicircular distribution \(\eta\) onto \(\nu_u\) is of class
\(\mathcal C^3\). Then
\[
\lambda_k(\Delta_{\nu_u})
\geq
ck,
\qquad
k\geq0.
\]
\end{proposition}
\begin{proof}
If $X$ has law $\nu_u$, then its conjugate variable is $\xi=u'(X)$ and
\[
\partial u'(x,y)
=
\frac{u'(x)-u'(y)}{x-y}
\geq c,
\]
and on the diagonal $\partial u'(x,x)=u''(x)$. Thus $X$ satisfies $CD(c,\infty)$, and Theorem~\ref{main} gives compact resolvent. Moreover, the free Caffarelli contraction theorem~\cite[Theorem~7]{Di4} gives $\operatorname{Lip}(T)\leq c^{-1/2}$. The map $f\mapsto f\circ T$ is an isometry from $L^2(\nu_u)$ into $L^2(\eta)$, and the chain rule for the free difference quotient yields, first for a polynomial $f$, the following comparison of the Dirichlet forms $\mathcal E_\eta$ and $\mathcal E_{\nu_u}$:
\[
\mathcal E_\eta(f\circ T)
\leq
\frac1c\mathcal E_{\nu_u}(f).
\]
Standard approximation for $f$ and closedness then extend this inequality to the full form domain.
The min--max principle therefore yields the following comparison of eigenvalues: 
\(
\lambda_k(\Delta_{\nu_u})\geq c\lambda_k(\Delta_\eta)=ck
\)
for every $k\geq0$. The reader can consult similar details in E. Milman~\cite[Section~3]{MilmanSpectral} and \cite{Di4} for precise explicit computations of the free Laplacian in this one-dimensional case.
\end{proof}
The preceding proposition lies exactly in our curvature--dimension setting. What remains open, even in dimension one, is whether the same spectral comparison follows directly from the curvature condition, without requiring a free Caffarelli contraction theorem. This problem in the classical case is already mentioned by E.~Milman~\cite[Section~1.1]{MilmanSpectral}, as it is not clear how a \(CD(c,\infty)\) criterion, i.e. information at the level of \(\Gamma_2\), is enough by itself to encode the entire eigenvalue structure. For the moment, the work of Ledoux~\cite{Ledoux3} shows that higher-order iterated \emph{carr\'e du champ operators} are useful to obtain higher-order eigenvalue estimates for the generator, suggesting that information at the level of \(\Gamma_2\) alone is not enough, at least with the presently known methods.

\section{Acknowledgments}

This work was supported by ISF grant 2574/24 and NSF--BSF grant 2022707. The author is grateful to Prof.~David Jekel for many helpful comments and discussions on an early version of the manuscript, and in particular for suggesting the argument that yields an exponential
upper bound for the spectral counting function.

\subsection{AI Declaration}

The overall mathematical framework of this work was developed by the author. GPT-6 Pro (OpenAI) was used as a research assistant, with mathematical assistance focused on two specific parts of the work: the construction of an embedding into a free group factor under the curvature--dimension criterion (Theorem~\ref{thm:curvature-embedding}) and the proof of the free KLS inequality (Theorem~\ref{thm:free-KLS}). The arguments suggested through these interactions were subsequently checked, revised, and substantially reorganized by the author before being incorporated into the manuscript.
GPT-6 Pro was also used for language editing and to improve the clarity and presentation of several parts of the paper. The author takes full responsibility for the mathematical statements, the verification of the proofs, and the final
content of the manuscript.


\begin{thebibliography}{99}

\bibitem{AnanPopa}
\textsc{C. Anantharaman and S. Popa},
\textit{An introduction to $\mathrm{II}_1$ factors}, \url{https://www.math.ucla.edu/~popa/Books/IIunV15.pdf}.

\bibitem{jek3}
\textsc{T. Amrutam, D. Jekel and M. Wasilewski},
\textit{Simplicity of $q$-Gaussian $C^*$-algebras},
arXiv:2606.10723, 2026.

\bibitem{Bakry}
\textsc{D. Bakry and M. \'Emery},
\textit{Diffusions hypercontractives},
in \textit{S\'eminaire de Probabilit\'es XIX, 1983/84},
Lecture Notes in Math., vol.~1123, Springer, Berlin, 1985, 177--206.

\bibitem{BS}
\textsc{P. Biane and R. Speicher},
\textit{Stochastic calculus with respect to free Brownian motion and analysis on Wigner space},
Probab. Theory Related Fields \textbf{112} (1998), 373--409.

\bibitem{bizeul}
\textsc{P. Bizeul},
\textit{The slicing conjecture via small ball estimates},
to appear in Ann. Probab.; arXiv:2501.06854, 2025.

\bibitem{bizeul2}
\textsc{P. Bizeul},
\textit{Optimal $MM^*$ bounds for convex bodies},
arXiv:2607.29458, 2026.

\bibitem{bizkla}
\textsc{P. Bizeul, B. Klartag and J. Lehec},
\textit{Presenting a proof of the
Kannan--Lov\'asz--Simonovits conjecture},
arXiv:2610.05474, 2026.

\bibitem{Bouyrie}
\textsc{R. Bouyrie},
\textit{Rigidity phenomenons for an infinite dimension diffusion operator and cases of near equality in the Bakry--Ledoux isoperimetric comparison theorem},
arXiv:1708.07203, 2017.

\bibitem{Bozejko}
\textsc{M. Bo\.{z}ejko},
\textit{Ultracontractivity and strong Sobolev inequality for
$q$-Ornstein--Uhlenbeck semigroup ($-1<q<1$)},
Infin. Dimens. Anal. Quantum Probab. Relat. Top.
\textbf{2} (1999), no.~2, 203--220.

\bibitem{BSK97}
\textsc{M. Bo\.{z}ejko, B. K\"ummerer and R. Speicher},
\textit{$q$-Gaussian processes: non-commutative and classical aspects},
Comm. Math. Phys. \textbf{185} (1997), no.~1, 129--154.

\bibitem{cafa}
\textsc{L. A. Caffarelli},
\textit{Monotonicity properties of optimal transportation and the FKG
and related inequalities},
Comm. Math. Phys. \textbf{214} (2000), no.~3, 547--563;
erratum, Comm. Math. Phys. \textbf{225} (2002), no.~2, 449--450.

\bibitem{CGP}
\textsc{A. Chandra, M. Gordina and M. Peev},
\textit{Free Bakry--\'Emery calculus},
arXiv:2608.23948, 2026.

\bibitem{Chenkla}
\textsc{Y. Chen and B. Klartag},
\textit{Digesting the proof of the sharp thin-shell inequality},
arXiv:2607.23307, 2026.

\bibitem{Choda}
\textsc{M. Choda},
\textit{Group factors of the Haagerup type},
Proc. Japan Acad. Ser. A Math. Sci. \textbf{59} (1983), no.~5, 174--177.

\bibitem{Cianchi}
\textsc{A. Cianchi and V. Maz'ya},
\textit{On the discreteness of the spectrum of the Laplacian on noncompact Riemannian manifolds},
J. Differential Geom. \textbf{87} (2011), no.~3, 469--491.

\bibitem{CiprianiSauvageot2003}
\textsc{F. Cipriani and J.-L. Sauvageot},
\textit{Derivations as square roots of Dirichlet forms},
J. Funct. Anal. \textbf{201} (2003), no.~1, 78--120.

\bibitem{Connes}
\textsc{A. Connes},
\textit{Classification of injective factors. Cases $\mathrm{II}_1$, $\mathrm{II}_{\infty}$, $\mathrm{III}_{\lambda}$, $\lambda\neq1$},
Ann. of Math. (2) \textbf{104} (1976), 73--115.

\bibitem{courtfa}
\textsc{T. A. Courtade and M. Fathi},
\textit{Stability of Klartag's improved Lichnerowicz inequality},
in \textit{High Dimensional Probability X},
Progr. Probab., vol.~82, Birkh\"auser, Cham, 2026, 91--99.

\bibitem{haagcowl}
\textsc{M. Cowling and U. Haagerup},
\textit{Completely bounded multipliers of the Fourier algebra of a simple Lie group of real rank one},
Invent. Math. \textbf{96} (1989), 507--549.

\bibitem{Dab10}
\textsc{Y. Dabrowski},
\textit{A note about proving non-$\Gamma$ under a finite non-microstates free Fisher information assumption},
J. Funct. Anal. \textbf{258} (2010), no.~11, 3662--3674.

\bibitem{Dab14}
\textsc{Y. Dabrowski},
\textit{A free stochastic partial differential equation},
Ann. Inst. Henri Poincar\'e Probab. Stat. \textbf{50} (2014),
no.~4, 1404--1455.

\bibitem{Dabioa}
\textsc{Y. Dabrowski and A. Ioana},
\textit{Unbounded derivations, free dilations, and indecomposability
results for $\mathrm{II}_1$ factors},
Trans. Amer. Math. Soc. \textbf{368} (2016), no.~7, 4525--4560.

\bibitem{Davies}
\textsc{E. B. Davies},
\textit{Heat Kernels and Spectral Theory},
Cambridge Tracts in Mathematics, vol.~92, Cambridge University Press, Cambridge, 1989.

\bibitem{Diez2}
\textsc{C.-P. Diez},
\textit{Free Malliavin--Stein--Dirichlet method: multidimensional semicircular approximations and chaos of a quantum Markov operator},
arXiv:2211.07595, 2022.

\bibitem{diez3}
\textsc{C.-P. Diez},
\textit{Free Stein kernels, free moment maps, and higher order derivatives},
arXiv:2306.00967, 2023.

\bibitem{Di4}
\textsc{C.-P. Diez},
\textit{Free Stein kernels and moment maps},
J. Anal. Math., in press; arXiv:2410.02470v2, 2024.

\bibitem{DiezObata}
\textsc{C.-P. Diez},
\textit{Obata's rigidity theorem in free probability},
arXiv:2603.05466, 2026.

\bibitem{DKEP}
\textsc{C. Ding, S. Kunnawalkam Elayavalli and J. Peterson},
\textit{Properly proximal von Neumann algebras},
Duke Math. J. \textbf{172} (2023), no.~15, 2821--2894.

\bibitem{DP}
\textsc{C. Ding and J. Peterson},
\textit{Biexact von Neumann algebras},
arXiv:2309.10161, 2023.

\bibitem{Fathi}
\textsc{M. Fathi},
\textit{Higher-order Stein kernels for Gaussian approximation},
Studia Math. \textbf{256} (2021), no.~3, 241--258.

\bibitem{gross}
\textsc{L. Gross},
\textit{Logarithmic Sobolev inequalities},
Amer. J. Math. \textbf{97} (1975), no.~4, 1061--1083.

\bibitem{Guan}
\textsc{Q. Guan},
\textit{A note on Bourgain's slicing problem},
arXiv:2412.09075, 2024.

\bibitem{GS}
\textsc{A. Guionnet and D. Shlyakhtenko},
\textit{Free diffusions and matrix models with strictly convex interaction},
Geom. Funct. Anal. \textbf{18} (2009), no.~6, 1875--1916.

\bibitem{GS2}
\textsc{A. Guionnet and D. Shlyakhtenko},
\textit{Free monotone transport},
Invent. Math. \textbf{197} (2014), no.~3, 613--661.

\bibitem{Haag}
\textsc{U. Haagerup},
\textit{An example of a non nuclear $C^*$-algebra, which has the metric
approximation property},
Invent. Math. \textbf{50} (1979), no.~3, 279--293.

\bibitem{Isono}
\textsc{Y. Isono},
\textit{Examples of factors which have no Cartan subalgebras},
Trans. Amer. Math. Soc. \textbf{367} (2015), no.~11, 7917--7937.

\bibitem{JekelLandry}
\textsc{D. Jekel and T. Basa Landry},
\textit{Compact quantum metric spaces from free probability},
arXiv:2606.19578, 2026.

\bibitem{Judge}
\textsc{C. Judge and S. Mondal},
\textit{Euclidean triangles have no hot spots},
Ann. of Math. (2) \textbf{191} (2020), no.~1, 167--211;
erratum, Ann. of Math. (2) \textbf{195} (2022), no.~1, 337--362.

\bibitem{KLS}
\textsc{R. Kannan, L. Lov\'asz and M. Simonovits},
\textit{Isoperimetric problems for convex bodies and a localization lemma},
Discrete Comput. Geom. \textbf{13} (1995), no.~3--4, 541--559.

\bibitem{Kato}
\textsc{T. Kato},
\textit{Perturbation Theory for Linear Operators},
2nd ed., corrected printing, Grundlehren der mathematischen
Wissenschaften, vol.~132, Springer-Verlag, Berlin, 1980.

\bibitem{Klart}
\textsc{B. Klartag},
\textit{Logarithmically-concave moment measures I},
in \textit{Geometric Aspects of Functional Analysis},
B. Klartag and E. Milman (eds.), Lecture Notes in Math., vol.~2116,
Springer, Cham, 2014, 231--260.

\bibitem{KlartagPutterman}
\textsc{B. Klartag and E. Putterman},
\textit{Spectral monotonicity under Gaussian convolution},
Ann. Fac. Sci. Toulouse Math. (6) \textbf{32} (2023), no.~5, 939--967.

\bibitem{Klartaglog}
\textsc{B. Klartag},
\textit{Logarithmic bounds for isoperimetry and slices of convex sets},
Ars Inveniendi Analytica (2023), Paper No.~4, 17 pp.

\bibitem{klalehec}
\textsc{B. Klartag and J. Lehec},
\textit{Isoperimetric inequalities in high-dimensional convex sets},
Bull. Amer. Math. Soc. (N.S.) \textbf{62} (2025), no.~4, 575--642.

\bibitem{Klalehec1}
\textsc{B. Klartag and J. Lehec},
\textit{Affirmative resolution of Bourgain's slicing problem using
Guan's bound},
Geom. Funct. Anal. \textbf{35} (2025), 1147--1168.

\bibitem{Klalehec2}
\textsc{B. Klartag and J. Lehec},
\textit{Thin-shell bounds via parallel coupling},
arXiv:2507.15495, 2025.

\bibitem{Ledoux3}
\textsc{M. Ledoux},
\textit{L'alg\`ebre de Lie des gradients it\'er\'es d'un g\'en\'erateur
markovien---d\'eveloppements de moyennes et entropies},
Ann. Sci. \'Ecole Norm. Sup. (4) \textbf{28} (1995), no.~4, 435--460.

\bibitem{Ledoux2}
\textsc{M. Ledoux},
\textit{Chaos of a Markov operator and the fourth moment condition},
Ann. Probab. \textbf{40} (2012), no.~6, 2439--2459.

\bibitem{Ledoux}
\textsc{D. Bakry, I. Gentil and M. Ledoux},
\textit{Analysis and Geometry of Markov Diffusion Operators},
Grundlehren der mathematischen Wissenschaften, vol.~348,
Springer, Cham, 2014.

\bibitem{Letwin}
\textsc{B. Letwin},
\textit{The KLS constant is $O(\log^{1/4}n)$},
arXiv:2607.24164, 2026.

\bibitem{lee}
\textsc{Y. T. Lee and S. S. Vempala},
\textit{The Kannan--Lov\'asz--Simonovits conjecture},
in \textit{Current Developments in Mathematics 2017},
International Press, Somerville, MA, 2019, 1--36.

\bibitem{Ma}
\textsc{Z.-M. Ma and M. R\"ockner},
\textit{Introduction to the Theory of (Non-Symmetric) Dirichlet Forms},
Universitext, Springer-Verlag, Berlin, 1992.

\bibitem{Mai}
\textsc{T. Mai},
\textit{Regularity of distributions of Wigner integrals},
arXiv:1512.07593, 2015; revised version arXiv:1512.07593v2, 2016.

\bibitem{MS}
\textsc{T. Mai and R. Speicher},
\textit{A note on the free and cyclic differential calculus},
J. Operator Theory \textbf{85} (2021), no.~1, 183--215.

\bibitem{Mikulin}
\textsc{W. Dudarov and D. Mikulincer},
\textit{Geometric obstructions to Lipschitz transport between weighted
Hessian $\mathrm{CD}(\kappa,\infty)$ manifolds},
arXiv:2606.11085, 2026.

\bibitem{MilmanSpectral}
\textsc{E. Milman},
\textit{Spectral estimates, contractions and hypercontractivity},
J. Spectr. Theory \textbf{8} (2018), no.~2, 669--714.
\bibitem{openai}
\textsc{OpenAI}, 
\textit{An isomorphism of the free group factors}, OpenAI Math Release preprint
\href{https://github.com/openai/math/blob/main/preprints/An-isomorphism-of-the-free-group-factors-September-23-2026/An-isomorphism-of-the-free-group-factors-September-23-2026.pdf}{OAI:An-isomorphism-of-the-free-group-factors-September-23-2026, 2026}.
\bibitem{openai2}
\textsc{OpenAI}
\textit{A finite-entropy separation of microstates and nonmicrostates free entropy},
\href{https://github.com/openai/math/blob/main/preprints/A-finite-entropy-separation-of-microstates-and-nonmicrostates-free-entropy-September-25-2026/paper.pdf}{OAI:A-finite-entropy-separation of-microstates-and-nonmicrostates-free-entropy-September-23-2026, 2026}.

\bibitem{Ozawa}
\textsc{N. Ozawa},
\textit{Solid von Neumann algebras},
Acta Math. \textbf{192} (2004), no.~1, 111--117.

\bibitem{OzawaWeak2}
\textsc{N. Ozawa},
\textit{Weakly exact von Neumann algebras},
J. Math. Soc. Japan \textbf{59} (2007), no.~4, 985--991.

\bibitem{OP1}
\textsc{N. Ozawa and S. Popa},
\textit{On a class of $\mathrm{II}_1$ factors with at most one Cartan subalgebra},
Ann. of Math. (2) \textbf{172} (2010), 713--749.

\bibitem{OP2}
\textsc{N. Ozawa and S. Popa},
\textit{On a class of $\mathrm{II}_1$ factors with at most one Cartan subalgebra, II},
Amer. J. Math. \textbf{132} (2010), 841--866.

\bibitem{OzawaWeak}
\textsc{N. Ozawa},
\textit{Examples of groups which are not weakly amenable},
Kyoto J. Math. \textbf{52} (2012), 333--344.

\bibitem{Peterson}
\textsc{J. Peterson},
\textit{$L^2$-rigidity in von Neumann algebras},
Invent. Math. \textbf{175} (2009), no.~2, 417--433.

\bibitem{ReedSimonIV}
\textsc{M. Reed and B. Simon},
\textit{Methods of Modern Mathematical Physics IV: Analysis of Operators},
Academic Press, New York, 1978.

\bibitem{Sauvage}
\textsc{J.-L. Sauvageot},
\textit{Strong Feller semigroups on $C^*$-algebras},
J. Operator Theory \textbf{42} (1999), no.~1, 83--102.

\bibitem{song}
\textsc{Z. Song and X. Zhang},
\textit{An $O(1)$ bound for the KLS constant},
arXiv:2610.01447v2, 2026.

\bibitem{Voic2}
\textsc{D.-V. Voiculescu},
\textit{The analogues of entropy and of Fisher's information measure
in free probability theory, II},
Invent. Math. \textbf{118} (1994), 411--440.

\bibitem{Voic}
\textsc{D.-V. Voiculescu},
\textit{The analogues of entropy and of Fisher's information measure
in free probability theory V. Noncommutative Hilbert transforms},
Invent. Math. \textbf{132} (1998), no.~1, 189--227.

\bibitem{Voipub}
\textsc{D.-V. Voiculescu},
\textit{A free Poincar\'e inequality},
unpublished manuscript.

\bibitem{Yang}
\textsc{F.-Y. Wang},
\textit{Functional inequalities for empty essential spectrum},
J. Funct. Anal. \textbf{170} (2000), no.~1, 219--245.

\bibitem{Aim}
\textsc{AIM06},
\textit{Problems posed in workshop “Free Analysis” at American Institute of Mathematics, 2006}, \url{https://aimath.org/WWN/freeanalysis/freeanalysis.pdf}.

\end{thebibliography}
\end{document}